# Università degli Studi di Trento

DIPARTIMENTO DI MATEMATICA
Doctoral school in Mathematics

On a *co-tutelle* agreement with the University of Stuttgart

# Gluing silting objects along recollements of well generated triangulated categories

Fabiano Bonometti

2019

# Abstract

We provide an explicit procedure to glue (not necessarily compact) silting objects along recollements of triangulated categories with coproducts having a 'nice' set of generators, namely, well generated triangulated categories. This procedure is compatible with gluing co-t-structures and it generalizes a result by Liu, Vitória and Yang. We provide conditions for our procedure to restrict to tilting objects and to silting and tilting modules. As applications, we retrieve the classification of silting modules over the Kronecker algebra and the classification of non-compact tilting sheaves over a weighted noncommutative regular projective curve of genus 0.

# Acknowledgments

First of all I want to express my gratitude to Lidia Angeleri, my supervisor in Verona and to Steffen König, in Stuttgart, for their guidance, expertise, and all the opportunities they have given me.

I also wish to thank all my current and former colleagues both at the University of Verona and at the University of Stuttgart for the scientific and personal exchange, and for providing happy distractions from work.

Special thanks to my parents for their unconditional trust and to my brothers and sister for their unique encouragement. Thanks to the other family members and to my friends that have always been there for me.

I also want to acknowledge the joint doctoral program in Mathematics of the Universities of Trento and Verona for the financial support, and the Short-Term Research Grant program of DAAD for the financial support during my stay in Stuttgart.

# Contents

# Introduction

Representation theory aims to explore algebraic structures such as algebras and rings through their actions on simpler objects (such as vector spaces). Classifying these actions, called *representations*, is the main aim of representation theory.

To achieve this goal there are two main ways: the first one is to explicitly compute representations, the second one is to compare the categories of representations of different rings, or equivalently, their categories of modules. This second way was initiated in 1958 by Morita, who considered two rings to be equivalent (now *Morita equivalent*) if, roughly speaking, they share the same representation theory. An example is any ring $R$ being Morita equivalent to the ring of matrices $M_n(R)$ with entries in $R$. He proved that two rings are Morita equivalent if and only if one ring is the endomorphism ring of a special module (a *projective generator*) over the other. Indeed $M_n(R) \cong \mathrm{End}_R(R^{\oplus n})$.

Tilting theory was developed in the seventies as a generalization of the theory of Morita. A central theorem of Brenner and Butler (1979) established the equivalence of certain subcategories of $\mathrm{mod}\, R$ and $\mathrm{mod}\, \mathrm{End}\, T$, for a special $R$-module $T$ that they called *tilting*, which is a generalization of a projective generator. Some years later, Happel introduced into representation theory the concept of derived module categories of algebras, a formalism that had recently been developed by Grothendieck and Verdier. He proved [31] that the bounded derived categories of $\mathrm{mod}\, R$ and $\mathrm{mod}\, \mathrm{End}\, T$ are equivalent for a tilting module $T$. Shortly afterwards, Rickard [61] and Keller [33] established a *derived Morita theory* by means of tilting complexes.

Tilting theory can now be understood as a tool to study derived categories and their equivalences. It provides links between different areas of mathematics: a classical example is the derived equivalence between the category of quasicoherent sheaves over the projective line and the category of modules over the Kronecker algebra, linking algebraic geometry and representation theory.

Silting theory was initiated by Keller and Vossieck [36] as an extension of tilting theory: the conditions of tilting objects being exceptional is weakened to the condition of having no positive self-extensions, required for silting objects. Silting theory was rediscovered years later by several authors [2, 19, 32, 49, 56] as a way to overcome a problem in tilting theory, namely, that mutations of (compact) tilting objects are not always possible. Meanwhile, in [35, 37, 49], correspondences relating silting complexes to t-structures and co-t-structures were studied. These are particular torsion pairs, which are decompositions of triangulated categories and represent another tool to compare categories. For any finite dimensional $k$-algebra $A$, Koenig and Yang [37] established one-to-one correspondences between equivalence classes of compact silting objects in $\mathbf{K}^b(\operatorname{proj} A)$, bounded t-structures on $\mathbf{D}^b(\operatorname{mod} A)$ with length heart, and bounded co-t-structures on $\mathbf{K}^b(\operatorname{proj} A)$.

Initially, the definition of a silting object $T$ in a triangulated category $\mathcal{D}$ included a strong generation property and the study of silting objects was mainly concentrated on categories of compact objects. As happened for tilting modules and tilting complexes, the notion of silting complex was extended to the unbounded derived category of a ring (see [69] and [11]). A more general version of the Koenig-Yang bijections, which included these new 'large' silting objects was provided in [11].

A further step in this direction was made independently in [56] and [59]; the authors introduced a new and more abstract concept of silting object in a triangulated category with coproducts, defined here as an object $T$ such that the pair of subcategories $(T^{\perp_{>0}}, T^{\perp_{\leq 0}})$ is a t-structure, where $T^{\perp_{>0}}$ (resp. $T^{\perp_{\leq 0}}$) denotes the set of objects $X$ such that $\operatorname{Hom}(T, X[i]) = 0$ for $i > 0$ (resp. for $i \leq 0$). For a large class of triangulated categories, it was proved in [13] that silting objects correspond bijectively to certain TTF triples, that are triples of subcategories consisting of a t-structure and a left-adjacent co-t-structure.

TTF triples are closely related to recollements, which consist of a diagram of three triangulated categories together with six additive functors

$$\mathcal{Y} \xleftarrow{\overset{i^*}{\underset{i^!}{\longleftarrow}}}_{\overset{i_*}{\longrightarrow}} \mathcal{D} \xleftarrow{\overset{j_!}{\underset{j_*}{\longleftarrow}}}_{\overset{j^*}{\longrightarrow}} \mathcal{X} \ ,$$

satisfying a number of axioms. Each row of functors resembles a short exact sequence, and the category $\mathcal{Y}$ can be seen as a (categorical) quotient $\mathcal{D}/\mathcal{X}$. A recollement can be thought as a way to glue together triangulated categories and obtain a bigger one. The properties of $\mathcal{X}$, $\mathcal{Y}$ and $\mathcal{D}$ are closely related. In particular, given t-structures, respectively co-t-structures, in $\mathcal{X}$ and $\mathcal{Y}$, one can produce a t-structure (resp. co-t-structure) in the category $\mathcal{D}$ in a quite 'natural' way (the *glued* t-structure).

In the context of comparing different categories, several authors have studied the relations existing between silting or tilting objects in the three categories composing a recollement. In particular, given a recollement as above and silting objects $T_\mathcal{X}$ in $\mathcal{X}$ and $T_\mathcal{Y}$ in $\mathcal{Y}$, Liu, Vitória and Yang [46] and Saorín and Zvonareva [65] have proposed procedures to explicitly compute the silting object in $\mathcal{D}$ associated to the t-structure or co-t-structure obtained by gluing the t-structures or co-t-structures corresponding to $T_\mathcal{X}$ and $T_\mathcal{Y}$.

Our work fits into this context, but builds instead on a paper by Angeleri, Koenig and Liu [6], who construct a (compact) tilting object in $\mathcal{D}$ from two exceptional objects in $\mathcal{X}$ and $\mathcal{Y}$ describing the recollement. Following their ideas, we consider two objects $\sigma_1$ and $\sigma_2$ with no positive self-extensions in a triangulated category $\mathcal{D}$ with coproducts and a map $\alpha \colon \sigma_2 \to \sigma_1[1]$. We show that, when $\sigma_1$ and $\sigma_2$ satisfy some orthogonality conditions, the cocone $\tilde\sigma$ of $\alpha$ have no positive self-extensions if and only if the map $\alpha$ satisfies some approximation properties. Now, if $\sigma_1$ satisfies the additional hypothesis that the class $\sigma_1^{\perp_{>0}}$ is closed under coproducts (i.e. $\sigma_1$ is *partial silting*) and the category $\mathcal{D}$ is well-generated, then $\sigma_1$ gives rise to a recollement of triangulated categories

$$\sigma_1^{\perp_\mathbb{Z}} \xleftarrow{\overset{i^*}{\underset{i^!}{\longleftarrow}}} \mathcal{D} \xleftarrow{\overset{j_!}{\underset{j_*}{\longleftarrow}}} \mathrm{Loc}(\sigma_1) \ ,$$

and it is silting in $\mathrm{Loc}(\sigma_1)$ ([13, 55]). If $\sigma_2$ is silting in $\sigma_1^{\perp_\mathbb{Z}}$, then we use the result above to show (see Theorem 2.10) that, under some hypotheses, $\tilde{\sigma} \oplus i_*\sigma_2$ is silting in $\mathcal{D}$, where $\tilde{\sigma}$ is the cocone of the $\mathrm{Add}(i_*\sigma_2)$-precover $\alpha\colon i_*\sigma_2^{(I)} \to j_!\sigma_1$, for $I = \mathrm{Hom}_\mathcal{D}(i_*\sigma_2, j_!\sigma_1)$. We show that this procedure is compatible with gluing co-t-structures, in the sense that, whenever we can associate a co-t-structure to the silting objects, the glued silting object is associated to the glued co-t-structure. This means that the glued silting object is equivalent to the one constructed in [46] if both methods can be applied. However, our result applies also to large silting objects, in contrast with [46].

In Theorem 2.16, we prove a 'dual' version of the previous theorem: if $\alpha\colon i_*\sigma_2 \to j_*\sigma_1^I$ is an $\mathrm{Add}(j_*\sigma_1)$-preenvelope, then $\tilde{\sigma} \oplus i_*\sigma_2$ is silting. This result requires however some additional and much stronger hypotheses, in particular, $j_*\sigma_1$ is required to be compact and the set $I$ finite. This procedure is compatible with gluing t-structures along a lower adjacent recollement, if it exists. Thus, this variant gives a silting object that is equivalent to the one constructed in [65], when both procedures apply.

Both variants of our procedure restrict to tilting objects and to tilting modules when some natural conditions are satisfied. We also study the particular case in which the recollement is induced by a homological ring epimorphism $\lambda$, for example a homological universal localization of rings. If $\lambda\colon R \to S$ is a universal localization, or more generally a silting ring epimorphism, then we can associate a 2-term partial silting complex $\omega$ to it. In some particular case, for instance if $S$ has projective dimension at most 1 as $R$-module, the partial silting complex $\omega$ can be explicitly computed, and the induced recollement gets a nice form. Finally, we study the case in which the recollement is induced by a triangular matrix ring (that also comprises one-point extensions).

To conclude, we apply our results to two examples. In the first one we consider the category of modules over the Kronecker algebra $R = k(\ \bullet \rightrightarrows \bullet\ )$. The classification of silting $R$-modules was completed in [12] and consists of the compact tilting modules, the compact silting non-tilting ones (0, the simple projective and the simple injective), the non-compact tilting modules, which are parametrized by universal localizations of $R$ at a nonempty set of simple regular modules, and the Lukas tilting module.

A universal localization $R \to R_S$ at a simple regular module $S$ induces an embedding of the module category $\mathrm{Mod}\,k[x]$ over the polynomial ring into the category of $R$-modules and a recollement of derived categories. If we let $S$ run through the simple regular $R$-modules, knowing the classification of silting $k[x]$-modules, we can get all non-compact silting $R$-modules. The compact silting $R$-modules can be obtained by gluing along a recollement induced by a localization at a preprojective or preinjective $R$-module.

For the second example we consider the category $\mathrm{Qcoh}\,\mathbb{X}$ of quasicoherent sheaves over a weighted noncommutative regular projective curve of genus zero.

The classification of non-compact tilting sheaves over $\mathbb{X}$ was completed in [9]. To talk about objects in nonhomogeneous tubes (or even their closure under direct limits), we make use of a geometric model introduced in [16], where a tube of rank $n$ is represented as an annulus with $n$ marked points on the outer boundary, finite length indecomposable objects are oriented arcs joining two such points, and Prüfer objects are infinite arcs starting from a point and spiraling around the inner boundary (see Figure 1).

**Figure 1:** An example of a Prüfer and a length 2 object in a tube of rank $n$, in the geometric model of [16].

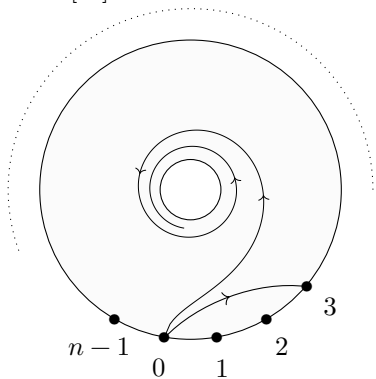

Since extensions correspond to negative crossings between arcs ([16, Theorem 4.3]), the model is very handy to visualize exceptional objects in tubes. Assuming that there is at least one exceptional tube in $\operatorname{Qcoh}\mathbb{X}$, we can get a recollement of triangulated categories by means of *expansions* of abelian categories, a construction by Chen and Krause [24] which, in this case, provides an alternative point of view for reduction of weights [30, §9]. More precisely, we can arbitrarily choose a simple exceptional sheaf $S$ at a point of weight $r \geq 2$: the right perpendicular category $S^{\perp_{0,1}}$ is again a category $\operatorname{Qcoh}\mathbb{X}'$ of quasicoherent sheaves, where the point $x$ in $\mathbb{X}'$ has weight $r-1$ and there is a recollement of triangulated categories ([24, Theorem 4.2.1])

$$\mathbf{D}^b(\operatorname{coh}\mathbb{X}') \xrightarrow[\substack{\longleftarrow i^* \\ \longrightarrow i_* \\ \longleftarrow i^!}]{} \mathbf{D}^b(\operatorname{coh}\mathbb{X}) \xrightarrow[\substack{\longleftarrow j_! \\ \longrightarrow j^* \\ \longleftarrow j_*}]{} \mathbf{D}^b(\operatorname{mod} k).$$

It turns out that gluing tilting objects with respect to Theorem 2.10, always yields a glued object with a new summand from the 'inflated' tube. In contrast, gluing with respect to Theorem 2.16, in some cases, keeps the component in this tube 'as it is' while modifying the torsionfree part.

By appropriately choosing the simple object $S$ and a non-compact tilting sheaf in $\operatorname{Qcoh}\mathbb{X}'$, we recover the whole classification of non-compact tilting sheaves in $\operatorname{Qcoh}\mathbb{X}$, including the Lukas sheaf.

The work is structured as follows. In Chapter 1 we introduce the notation and basic definitions. We present the concepts of recollements, TTF triple, t-structure and recollements of triangulated categories and see how t-structures and co-t-structures are related in the three categories forming a recollement. We motivate why well-generated categories are a suitable environment to work in, and present the notions of silting object, (2-term) silting complex and silting module. Finally, we introduce ring epimorphisms, universal localizations and their connections with (partial) silting complexes.

Chapter 2 is devoted to the proof of the main theorem and its dual version. In the first two sections we present some general results about universal maps and prove some preliminary lemmas. After proving the main theorem, we explore the compatibility of the gluing procedure with the gluing of t-structures and co-t-structures. We also restrict the gluing procedures to tilting objects and tilting modules, and then explore the particular case of a recollement induced by a silting ring epimorphism and a triangular matrix ring. Finally, we compare our result to the gluing procedures of other authors, more precisely, the ones appearing in [46] and [65].

In Chapter 3 we recover the classification of silting modules over the Kronecker ring $R$. To this aim, we consider different recollements involving the derived category of the Kronecker algebra over a ring $k$, induced by universal localization. Actually, since $R$ is a hereditary ring, all recollements of $\mathbf{D}(R)$ are induced by a universal localization (see Proposition 1.9 and [40, Theorem 8.1]). To get the large silting modules, we consider, in turn, the universal localizations at sets $\{S\}$ consisting of a single simple regular module, which induces a recollement

$$\mathbf{D}(k[x]) \cong \mathbf{D}(\mathcal{R}_S) \xleftarrow{\overset{i^*}{\longleftarrow}}_{\overset{i_*}{\longrightarrow}}_{\overset{i^!}{\longleftarrow}} \mathbf{D}(\mathcal{R}) \xleftarrow{\overset{j_!}{\longleftarrow}}_{\overset{j^*}{\longrightarrow}}_{\overset{j_*}{\longleftarrow}} \mathrm{Loc}(S_\infty)$$

involving the derived category of the polynomial ring and the localizing subcategory generated by the partial silting object $S_\infty$, that is the (not finitely generated) Prüfer module with simple socle $S$. To get the compact silting modules, we consider instead the localization at a preprojective or preinjective module. In this case, the induced recollement involves the localizing subcategories generated by two different finitely presented modules.

Finally, in Chapter 4 we start by studying the two tools that are needed: the expansion of abelian categories ([24]), and the geometric model of objects in the tube ([16]). We present the classification of non-compact tilting objects in categories of coherent sheaves over a weighted noncommutative regular projective curve of genus zero that had been completed in [8] and [9]. Finally, we consider the recollement of derived categories of quasicoherent sheaves induced by an expansion of categories of coherent sheaves and recover the classification presented above by means of gluing.

# 1

# Definitions and preliminaries

We introduce here the notation and some preliminary results. After defining TTF triples and recollements, we present the concepts of t-structures and co-t-structures and see how they naturally glue along recollements. We compare different definitions of silting and tilting objects, (2-term) silting complexes, silting and tilting modules. Finally, we consider ring epimorphisms, universal localizations and their connections to (partial) silting complexes.

## 1.1 Preliminaries and notation

Throughout, $\mathcal{D}$ will be a triangulated category with shift functor $[1]$. All subcategories considered will be full and closed under isomorphisms. Given a triangle $X \to Y \to Z \to X[1]$ in $\mathcal{D}$, the object $Y$ is said to be an *extension* of $X$ and $Z$. For a class of objects $\mathcal{X}$ in $\mathcal{D}$ and a set of integers $I \in \mathbb{Z}$, we define the following orthogonal classes

$$\mathcal{X}^{\perp_I} = \{Y \in \mathcal{D} \mid \mathrm{Hom}_\mathcal{D}(Y, X[i]) = 0 \text{ for all } X \in \mathcal{X} \text{ and } i \in I\}$$
$$^{\perp_I}\mathcal{X} = \{Y \in \mathcal{D} \mid \mathrm{Hom}_\mathcal{D}(X, Y[i]) = 0 \text{ for all } X \in \mathcal{X} \text{ and } i \in I\}.$$

We will often replace the set $I$ by notation such as $n, > n, \geq n, < n, \leq n, \neq n$, with the obvious associated meaning.

We denote by $\mathrm{add}(\mathcal{X})$ the smallest subcategory of $\mathcal{D}$ containing $\mathcal{X}$ and closed under finite coproducts and direct summands. If $\mathcal{D}$ has coproducts (respectively, products), we denote by $\mathrm{Add}(\mathcal{X})$ (respectively, $\mathrm{Prod}(\mathcal{X})$) the smallest subcategory of $\mathcal{D}$ containing $\mathcal{X}$ and closed under coproducts (respectively, products) and summands. If $\mathcal{X}$ consists of a single object $X$, we write $X^{\perp_I}$, $^{\perp_I}X$, $\mathrm{add}(X)$, $\mathrm{Add}(X)$ and $\mathrm{Prod}(X)$.

An additive functor from $\mathcal{D}$ to and abelian category $\mathcal{A}$ is *cohomological* if it takes triangles in $\mathcal{D}$ to long exact sequences in $\mathcal{A}$. A *triangulated subcategory* of $\mathcal{D}$ is a subcategory which is closed under extensions and shifts; if it is also closed under summands, then it is called *thick*. A subcategory of $\mathcal{D}$ is *suspended* if it is closed under extensions and positive shifts, *cosuspended* if it is closed under extensions and negative shifts. Given a class of objects $\mathcal{X}$ in $\mathcal{D}$, we denote by thick($\mathcal{X}$) the smallest thick subcategory of $\mathcal{D}$ containing $\mathcal{X}$. We denote by Susp($\mathcal{X}$) the smallest suspended subcategory of $\mathcal{D}$ containing $\mathcal{X}$ and closed under all existing coproducts. We say that a class of objects $\mathcal{X}$ *(weakly) generates* $\mathcal{D}$ if $\mathcal{X}^{\perp_{\mathbb{Z}}} = 0$; when $\mathcal{X}$ consists of a single object $X$, we call $X$ a *(weak) generator*. Of course, $\mathcal{X}$ generates $\mathcal{D}$ whenever thick($\mathcal{X}$) = $\mathcal{D}$.

Given an abelian category $\mathcal{A}$, we denote by $\mathbf{D}(\mathcal{A})$ its unbounded derived category. For a ring $R$, we denote by Mod $R$ the category of all right $R$-modules and write $\mathbf{D}(R) = \mathbf{D}(\text{Mod}\, R)$. The subcategories of projective and of finitely generated projective $R$-modules are denoted, respectively, by Proj $R$ and proj $R$. Their bounded homotopy categories are denoted by $\mathbf{K}^b(\text{Proj}\, R)$ and $\mathbf{K}^b(\text{proj}\, R)$, respectively. The category of all finitely presented right $R$-modules is denoted by mod $R$, and $\mathbf{D}^b(\text{mod}\, R)$ is its bounded derived category. Given an object $X$ in an abelian category $\mathcal{A}$ possessing arbitrary coproducts, we denote by $\text{Gen}(X)$ the class of all modules that are epimorphic images of objects in $\text{Add}(X)$.

## 1.2 TTF triples and recollements of triangulated categories

A *torsion pair* in a triangulated category $\mathcal{D}$ is a pair of subcategories $(\mathcal{U}, \mathcal{V})$ of $\mathcal{D}$ such that

(1) $\mathcal{U}$ and $\mathcal{V}$ are closed under direct summands;

(2) $\text{Hom}_{\mathcal{D}}(\mathcal{U}, \mathcal{V}) = 0$;

(3) every object $X$ of $\mathcal{D}$ fits into a triangle $U \to X \to V \to U[1]$ with $U \in \mathcal{U}$ and $V \in \mathcal{V}$.

A torsion pair $(\mathcal{V}, \mathcal{W})$ is a *t-structure* [17] if $\mathcal{V}$ is suspended, or equivalently $\mathcal{W}$ is cosuspended. It is a *co-t-structure* [19, 57] if $\mathcal{W}$ is suspended, or equivalently if $\mathcal{V}$ is cosuspended. A torsion pair $(\mathcal{U}, \mathcal{V})$ is *left nondegenerate* if $\bigcap_{n \in \mathbb{Z}} \mathcal{U}[n] = 0$, *right nondegenerate* if $\bigcap_{n \in \mathbb{Z}} \mathcal{V}[n] = 0$, *nondegenerate* if it is left and right nondegenerate.

**Example 1.1.** In the derived category $\mathcal{D} = \mathbf{D}(R)$ of a module category consider the subcategories

$$D^{\leq n} = \{X \in \mathbf{D}(R) \mid H^0(X) = 0 \text{ for any } i > n\},$$
$$D^{>n} = \{X \in \mathbf{D}(R) \mid H^0(X) = 0 \text{ for any } i \leq n\}.$$

The pair of subcategories $(D^{\leq 0}, D^{>0})$ is a t-structure, called the *standard t-structure*.

A *TTF triple* is a triple $(\mathcal{U}, \mathcal{V}, \mathcal{W})$ of subcategories such that $(\mathcal{U}, \mathcal{V})$ and $(\mathcal{V}, \mathcal{W})$ are torsion pairs. A TTF triple is

- *stable* if $\mathcal{V}$ is a triangulated subcategory of $\mathcal{D}$,

- *suspended* if $\mathcal{V}$ is suspended,

- *generated by a set* of objects $\mathcal{S}$ if $\mathcal{V} = \mathcal{S}^{\perp_0}$.

A suspended TTF triple $(\mathcal{U}, \mathcal{V}, \mathcal{W})$ is *nondegenerate* if $(\mathcal{V}, \mathcal{W})$ is nondegenerate. Note that a TTF triple $(\mathcal{U}, \mathcal{V}, \mathcal{W})$ is suspended if and only if $(\mathcal{U}, \mathcal{V})$ is a co-t-structure, if and only if $(\mathcal{V}, \mathcal{W})$ is a t-structure.

**Definition 1.2.** [17, §1.4.1] A *recollement (of $\mathcal{D}$)* is a diagram

$$\mathcal{Y} \xleftarrow[i_*]{\overset{i^*}{\underset{i^!}{\longleftarrow}}} \mathcal{D} \xleftarrow[j_*]{\overset{j_!}{\underset{j^*}{\longleftarrow}}} \mathcal{X} \qquad (1.1)$$

of three triangulated categories and six triangulated functors satisfying

(1) $(i^*, i_*)$, $(i_*, i^!)$, $(j_!, j^*)$, $(j^*, j_*)$ are adjoint pairs;

(2) $i^! j_* = 0$ (and thus $j^* i_* = 0$ and $i^* j_! = 0$);

(3) $i_*$, $j_*$, $j_!$ are full embeddings (and thus $i^* i_* \cong i^! i_* \cong 1_\mathcal{Y}$);

It follows from the definition of recollement ([58]) that for each object $C$ in $\mathcal{D}$ there are triangles

$$i_* i^! C \to C \to j_* j^* C \to i_* i^! C[1] \text{ and}$$
$$j_! j^* C \to C \to i_* i^* C \to j_! j^* C[1].$$

Two recollements

$$\mathcal{Y} \xleftarrow[i_*]{\overset{i^*}{\underset{i^!}{\longleftarrow}}} \mathcal{D} \xleftarrow[j_*]{\overset{j_!}{\underset{j^*}{\longleftarrow}}} \mathcal{X}$$

and

$$\mathcal{Y}' \xleftarrow[i'_*]{\overset{i'^*}{\underset{i'^!}{\longleftarrow}}} \mathcal{D} \xleftarrow[j'_*]{\overset{j'_!}{\underset{j'^*}{\longleftarrow}}} \mathcal{X}'$$

are *equivalent* [55, Definition 4.2.2] if $\operatorname{Im}(j_!) = \operatorname{Im}(j'_!)$, $\operatorname{Im}(i_*) = \operatorname{Im}(i'_*)$ and $\operatorname{Im}(j_*) = \operatorname{Im}(j'_*)$, where by $\operatorname{Im}(i_*)$ we mean the essential image of $i_*$, and analogously with the other functors.

This equivalence relation is motivated by the relation between TTF triples and recollements. Indeed, stable TTF triples in a triangulated category $\mathcal{D}$ are in bijection with equivalence classes of recollements of $\mathcal{D}$ ([17, p. 1.4.1], [53, Subsection 9.2]). The bijection associates to a recollement like (1.1) the TTF triple $(j_!\mathcal{X}, i_*\mathcal{Y}, j_*\mathcal{X})$, where by $j_!\mathcal{X}$ we mean the essential image of $j_!$, and analogously with the other functors. Conversely, given a TTF triple $(\mathcal{U}, \mathcal{V}, \mathcal{W})$ in $\mathcal{D}$, there exists a recollement of $\mathcal{D}$ as follows

$$\mathcal{V} \xleftarrow[i_*]{\overset{}{\underset{}{\longleftarrow}}} \mathcal{D} \xleftarrow[]{\overset{j_!}{\underset{}{\longleftarrow}}} \mathcal{U},$$

where $i_*$ and $j_!$ are the inclusion functors of $\mathcal{V}$ and $\mathcal{U}$ into $\mathcal{D}$, respectively. See [55, §4.2] for more details.

Given torsion pairs in the two outer terms of a recollement, the *glued torsion pair* is defined in [19] along the lines of [17]; see also [46].

**Proposition 1.3.** *[19, Theorem 8.2.3],[17, Theorem 1.4.10] Consider a recollement*

$$\mathcal{Y} \underset{\underset{\leftarrow i^!}{\leftarrow i^*}}{\overset{\leftarrow i^*}{\underset{i_*}{\longrightarrow}}} \mathcal{D} \underset{\underset{\leftarrow j_*}{\leftarrow j_!}}{\overset{\leftarrow j_!}{\underset{j^*}{\longrightarrow}}} \mathcal{X} \ ,$$

*and assume torsion pairs $(\mathcal{X}', \mathcal{X}'')$ and $(\mathcal{Y}', \mathcal{Y}'')$ are given, respectively in $\mathcal{X}$ and $\mathcal{Y}$. Then there is a torsion pair $(\mathcal{D}', \mathcal{D}'')$ in $\mathcal{D}$, defined by*

$$\mathcal{D}' = \{Z \in \mathcal{D} \mid i^*Z \in \mathcal{Y}', j^*Z \in \mathcal{X}'\},$$
$$\mathcal{D}'' = \{Z \in \mathcal{D} \mid i^!Z \in \mathcal{Y}'', j^*Z \in \mathcal{X}''\}.$$

*If the torsion pairs in the outer categories are t-structures (resp. co-t-structures), then so is the glued torsion pair.*

## 1.3 Triangulated categories with coproducts

Let now $\mathcal{D}$ be a triangulated category with arbitrary (set-indexed) coproducts.

A subcategory $\mathcal{L}$ of $\mathcal{D}$ is

- *complete* if it is closed under products,
- *cocomplete* if it is closed under coproducts,
- *localizing* if it is triangulated and cocomplete,
- *smashing* if it a localizing subcategory such that the inclusion functor admits a right adjoint and $\mathcal{L}^{\perp_0}$ is cocomplete.

For a set $\mathcal{S}$ of objects in $\mathcal{D}$, denote by $\mathrm{Loc}(\mathcal{S})$ the smallest localizing subcategory of $\mathcal{D}$ containing $\mathcal{S}$. If $\mathcal{L}$ is a localizing subcategory and $\mathcal{S}$ is a set of objects in $\mathcal{D}$, we say that $\mathcal{S}$ *(strongly) generates* $\mathcal{L}$ if $\mathrm{Loc}(\mathcal{S}) = \mathcal{L}$. If $\mathcal{S}$ consists of a single object $X$, we call $X$ a *(strong) generator* of $\mathcal{L}$.

An object $X \in \mathcal{D}$ is said to be *compact* if the functor $\mathrm{Hom}_\mathcal{D}(X, -)$ commutes with coproducts. The subcategory of compact objects will be denoted by $\mathcal{D}^c$; if $\mathcal{D}^c$ is skeletally small and generates $\mathcal{D}$, then $\mathcal{D}$ is said to be *compactly generated*. By [53, Proposition 8.4.6 and Theorem 8.3.3], in this case $\mathcal{D}$ also admits products. For a ring $R$, it is well known that $\mathbf{K}^b(\mathrm{proj}\,R) = \mathbf{D}(R)^c$ and $\mathbf{D}(R)$ is compactly generated.

We will often consider a generalization of compactly generated categories, namely well generated categories, where things work well enough.

**Definition 1.4.** [38, 53] Given a regular cardinal $\alpha$ and a triangulated category $\mathcal{D}$, we say that

- an object $X$ in $\mathcal{D}$ is $\alpha$-*small* if given any map $h \colon X \to \coprod_{\lambda \in \Lambda} Y_\lambda$ for some family of objects $(Y_\lambda)_{\lambda \in \Lambda}$ in $\mathcal{D}$, the map $h$ factors through a subcoproduct $\coprod_{\omega \in \Omega} Y_\omega$, where $\Omega$ is a subset of $\Lambda$ of cardinality strictly less then $\alpha$;
- $\mathcal{D}$ is $\alpha$-*well generated* if it has set-indexed coproducts and it has a set of objects $\mathcal{S}$ such that
  - $\mathcal{S}^{\perp_\mathbb{Z}} = 0$;

- for every set of maps $(g_\lambda\colon X_\lambda \to Y_\lambda)_{\lambda\in\Lambda}$ in $\mathcal{D}$, if $\operatorname{Hom}_\mathcal{D}(S, g_\lambda)$ is surjective for all $\lambda \in \Lambda$ and all $S \in \mathcal{S}$, then $\operatorname{Hom}_\mathcal{D}(S, \coprod_{\lambda\in\Lambda} g_\lambda)$ is surjective for all $S \in \mathcal{S}$;
- every object $S$ in $\mathcal{S}$ is $\alpha$-small.

- $\mathcal{D}$ is *well generated* if it is $\alpha$-well generated for some regular cardinal $\alpha$.

Of course, compactly generated triangulated categories with coproducts are well generated categories. Other examples of well generated triangulated categories include

- the homotopy category of projective modules over a ring is $\aleph_1$-well generated ([52, Theorem 1.1]);

- the derived category of a Grothendieck abelian category $\mathcal{A}$ is $\alpha$-well generated, with $\alpha$ a regular cardinal that depends on $\mathcal{A}$ ([50, Theorem 0.2]).

One of the advantages of working in well generated triangulated categories is given by the following proposition.

**Proposition 1.5.** *[22, Proposition 3.8][54, Theorem 2.3] In a well generated triangulated category $\mathcal{D}$ the pairs of subcategories $(\operatorname{Loc}(\omega), \omega^{\perp_\mathbb{Z}})$ and $(\operatorname{Susp}(\omega), \omega^{\perp_{\leq 0}})$ are torsion pairs, for any object $\omega \in \mathcal{D}$; thus, if $\omega^{\perp_\mathbb{Z}} = 0$, then $\operatorname{Loc}(\omega) = \mathcal{D}$.*

In other words, in a well generated category the following properties are equivalent for an object $\omega$:

(i) $\omega^{\perp_\mathbb{Z}} = 0$ (*weak generating property*) and

(ii) $\operatorname{Loc}(\omega) = \mathcal{D}$ (*strong generating property*).

A further interesting property of well generated triangulated categories is related to representability of cohomological functors.

**Definition 1.6.** We say that the triangulated category $\mathcal{D}$ *satisfies Brown representability* if every (contravariant) cohomological functor $\mathcal{D}^{\operatorname{op}} \to \operatorname{Mod}\mathbb{Z}$ sending coproducts to products is representable. Dually, $\mathcal{D}$ satisfies *dual Brown representability* if it is complete and every (covariant) cohomological functor $\mathcal{D} \to \operatorname{Mod}\mathbb{Z}$ sending products to products is representable.

By [53, Proposition 8.4.6], if a triangulated category $\mathcal{D}$ satisfies Brown representability, then it is complete.

**Theorem 1.7.** *[53, Theorem 8.3.3 and Proposition 8.4.2] Any well generated triangulated category satisfies Brown representability. Any compactly generated triangulated category satisfies, in addition, dual Brown representability.*

The following proposition is a useful characterization of smashing subcategory in well generated triangulated categories ([55, Propositon 4.2.4 and 4.4.14], see also [13, Proposition 2.9]).

**Proposition 1.8.** *If $\mathcal{D}$ is a well generated triangulated category, then the assignment $\mathcal{L} \mapsto (\mathcal{L}, \mathcal{L}^{\perp_0}, (\mathcal{L}^{\perp_0})^{\perp_0})$ yields a bijection between smashing subcategories of $\mathcal{D}$ and stable TTF triples in $\mathcal{D}$.*

Taking into account the bijection between TTF triples and recollement, Proposition 1.8 yields the following.

**Proposition 1.9.** *For a well generated triangulated category $\mathcal{D}$, there is a bijections between*

(1) *smashing subcategories of $\mathcal{D}$ and*

(2) *recollements of the category $\mathcal{D}$.*

*The bijection assigns to any smashing subcategory $\mathcal{L}$ of $\mathcal{D}$ the recollement*

$$\mathcal{L}^{\perp_0} \xrightarrow[\xleftarrow{i^*}]{\xleftarrow{i_*}} \mathcal{D} \xrightarrow[\xleftarrow{j_!}]{\xleftarrow{j^*}} \mathcal{L}, \tag{1.2}$$

*where $j_!$ and $i_*$ are the respective embeddings of the subcategories $\mathcal{L}$ and $\mathcal{L}^{\perp_0}$.*

## 1.4 Silting and tilting objects

We collect here some of the definitions of silting and tilting objects in different settings and outline the differences and connections between them. Let $\mathcal{D}$ be a triangulated category with (set-indexed) coproducts (later on we will assume $\mathcal{D}$ to be well generated). Let $R$ be any ring. When not otherwise stated, $R$-modules will be *right* $R$-modules.

**Definition 1.10.** The following definition is due to Psaroudakis and Vitória [59] and it is equivalent to the one given by Nicolás, Saorín and Zvonareva [56]. Let $\mathcal{D}$ be a triangulated category with coproducts. We say an object $\sigma$ of $\mathcal{D}$ is *silting* if the pair $(\sigma^{\perp_{>0}}, \sigma^{\perp_{\leq 0}})$ is a t-structure. It is *tilting* if it is silting and $\mathrm{Add}(\sigma) \subseteq \sigma^{\perp_{\neq 0}}$.

Note that a silting object according to this definition is contained in $\sigma^{\perp_{>0}}$ and generates $\mathcal{D}$, i.e., $\sigma^{\perp_{\mathbb{Z}}} = 0$. Following Angeleri, Marks and Vitória [13] we can give the following equivalent definition of silting object in a well generated triangulated category:

**Definition 1.11.** [13, Definition 3.1 and Remark 3.2] An object $\sigma$ in a well generated category $\mathcal{D}$ is partial silting if

(S1) $\sigma$ lies in $\sigma^{\perp_{>0}}$ and

(S2) $\sigma^{\perp_{>0}}$ is closed under coproducts.

*Remark* 1.12. [13, 56] An object $\sigma$ in a well generated category $\mathcal{D}$ is silting if and only if it is partial silting and

(S3) $\sigma^{\perp_{\mathbb{Z}}} = 0$.

It is is *(partial) tilting* if and only if it is (partial) silting and

(T) $\sigma^{(J)} \in \sigma^{\perp_{\neq 0}}$ for any set $J$.

If $\mathcal{D}$ is a well generated triangulated category satisfying dual Brown representability, to a silting object $T$ in $\mathcal{D}$ we can associate a co-t-structure $(^{\perp_0}(T^{\perp_{>0}}), T^{\perp_{>0}})$, which is called the *co-t-structure associated to $T$*. The next result makes this more precise.

**Theorem 1.13.** *[13, Theorem 3.9] Let $\mathcal{D}$ be a well generated triangulated category satisfying dual Brown representability. The assignment sending an object $T$ to the triple $(^{\perp_0}(T^{\perp_{>0}}), T^{\perp_{>0}}, (T^{\perp_{>0}})^{\perp_0})$ yields a bijection between*

- *equivalence classes of partial silting objects in $\mathcal{D}$ and*

- *suspended TTF triples $(\mathcal{U}, \mathcal{V}, \mathcal{W})$ in $\mathcal{D}$ such that $(\mathcal{U}, \mathcal{V})$ is generated by a set of objects from $\mathcal{D}$ and $(\mathcal{V}, \mathcal{W})$ is right nondegenerate.*

From now on, $\mathcal{D}$ will be a well generated triangulated category. The following results will also be useful.

**Proposition 1.14.** *[13, Proposition 3.5] The following are equivalent for an object $\omega$ in $\mathcal{D}$:*

*(i) $\omega$ is partial silting,*

*(ii) $\mathrm{Loc}(\omega)$ is smashing and $\omega$ is silting in $\mathrm{Loc}(\omega)$.*

Combining Proposition 1.14 and Proposition 1.9, and noting that $\mathrm{Loc}(\omega)^{\perp_0} = \omega^{\perp_{\mathbb{Z}}}$, one gets the following recollement associated to a partial silting object $\omega$ in $\mathcal{D}$.

$$\omega^{\perp_{\mathbb{Z}}} \xrightarrow[\substack{\leftarrow i^* \\ \xrightarrow{i_*} \\ \leftarrow i^!}]{} \mathcal{D} \xrightarrow[\substack{\leftarrow j_! \\ \xrightarrow{j^*} \\ \leftarrow j_*}]{} \mathrm{Loc}(\omega) \ , \tag{1.3}$$

Working with recollements of well generated categories we will find useful the following result of Krause.

**Theorem 1.15.** *[39, Theorem 7.2.1] Given a well generated category $\mathcal{D}$ and a localizing subcategory $\mathcal{X}$ that is generated by a set of objects, assume there exists a recollement*

$$\mathcal{Y} \xrightarrow[\substack{\leftarrow i^* \\ \xrightarrow{i_*} \\ \leftarrow i^!}]{} \mathcal{D} \xrightarrow[\substack{\leftarrow j_! \\ \xrightarrow{j^*} \\ \leftarrow j_*}]{} \mathcal{X} \ , \tag{1.4}$$

*Then $\mathcal{X}$ and $\mathcal{Y}$ are well generated.*

From Theorem 1.15, it follows that $\mathrm{Loc}(\omega)$ and $\omega^{\perp_{\mathbb{Z}}}$ are well generated for any partial silting object $\omega$ in a well generated triangulated category $\mathcal{D}$. Next proposition states that recollements of the form (1.3) are indeed quite general.

**Proposition 1.16.** *Let $\mathcal{D}$ be a well generated triangulated category and assume there is a recollement like (1.4). If there exists a silting object $\omega$ in $\mathcal{X}$, then the recollement is of the form (1.3).*

*Proof.* By Krause's Theorem 1.15, the category $\mathcal{X}$ is well generated. Thus, by Proposition 1.5, $\mathcal{X} = \mathrm{Loc}(\omega)$ and consequently $\mathcal{Y} = \omega^{\perp_{\mathbb{Z}}}$. □

## 1.5 Silting complexes

Let now $R$ be a ring. A complex $\sigma$ in $\mathbf{K}^b(\operatorname{Proj} R)$ is *n-term* if it is concentrated between degrees $-n+1$ and $0$ ([11, §4]). Silting complexes in $\mathbf{D}(A)$ were introduced in [69] by Wei, who named them "big semi-tilting". See also [5, Definition 5.1].

**Definition 1.17.** A bounded complex of projective $R$-modules $\sigma$ in $\mathbf{K}^b(\operatorname{Proj} R)$ is a *silting complex* if it satisfies the following conditions:

(1) $\operatorname{Hom}_{\mathbf{D}(R)}(\sigma, \sigma^{(I)}[i]) = 0$ for all sets $I$ and $i > 0$;

(2) the smallest triangulated category of $\mathbf{D}(R)$ containing $\operatorname{Add}(\sigma)$ is $\mathbf{K}^b(\operatorname{Proj} R)$.

Moreover, $\sigma$ is a *tilting complex* if, in addition, it is compact and exceptional, that is, it belongs to $\mathbf{K}^b(\operatorname{proj} R)$ and satisfies $\operatorname{Hom}_{\mathbf{D}(R)}(\sigma, \sigma[i]) = 0$ for all $i \neq 0$.

By condition (1), the class $\sigma^{\perp_{>0}}$ given by a silting complex $\sigma$ contains $\operatorname{Susp}(\sigma)$ ([69, Corollary 2.6]). Conversely, given an object $X$ in $\sigma^{\perp_{>0}}$, we can consider the canonical triangle $V \to X \to W \to V[1]$ with respect to the t-structure $(\operatorname{Susp}(\sigma), \sigma^{\perp_{\leq 0}})$ from Proposition 1.5. Then $W$ belongs to $\sigma^{\perp_{\leq 0}}$, but also to $\sigma^{\perp_{\geq 0}}$ since so do $X$ and $V[1]$. Since $\sigma$ is a generator, we infer that $W = 0$ and thus $X$ belongs to $\operatorname{Susp}(\sigma)$. Then, the pair $(\sigma^{\perp_{>0}}, \sigma^{\perp_{\leq 0}})$ coincides with $(\operatorname{Susp}(\sigma), \sigma^{\perp_{\leq 0}})$ and hence it is a t-structure, that means $\sigma$ is a silting object in $\mathbf{D}(R)$.

The following proposition states the compatibility of the definitions of silting object and silting complex.

**Proposition 1.18.** *[11, Proposition 4.2] (see also [5, Proposition 5.3]) A bounded complex of projective $R$-modules is a silting object in $\mathbf{D}(R)$ if and only if it is a silting complex.*

Note, however, that a silting object in $\mathbf{D}(R)$ needs not be bounded ([5, Example 7.9]). A bounded complex of projectives $\sigma$ is called a *partial silting complex* if it is a partial silting object in $\mathbf{D}(R)$.

## 1.6 Silting and tilting modules

For a morphism $\sigma$ between projective $R$-modules, consider the following classes of $R$-modules

$$\mathcal{D}_\sigma = \{X \in \operatorname{Mod} R \mid \operatorname{Hom}_R(\sigma, X) \text{ is surjective}\};$$
$$\mathcal{Y}_\sigma = \{X \in \operatorname{Mod} R \mid \operatorname{Hom}_R(\sigma, X) \text{ is a bijection}\}.$$

Note that if $\omega$ is a complex in $\operatorname{Proj}(R)$ concentrated in degrees $0$ and $-1$ with $H^0(\omega) = X$, the following lemma holds true.

**Lemma 1.19.** (1) *The following conditions are equivalent for two 2-term complexes $\sigma$ and $\omega$ in $\mathbf{K}^b(\operatorname{Proj} R)$ with $H^0(\omega) = X$:*

*(i) the map $\operatorname{Hom}_R(\sigma, X)$ is surjective,*

*(ii) $\operatorname{Hom}_{\mathbf{D}(R)}(\sigma, X[1]) = 0$,*

*(iii)* $\mathrm{Hom}_{\mathbf{K}^b(\mathrm{Proj}\,R)}(\sigma, \omega[1]) = 0$.

*(2) In the same setting, the following are equivalent:*

*(i) the map* $\mathrm{Hom}_R(\sigma, X)$ *is bijective,*

*(ii)* $\mathrm{Hom}_{\mathbf{D}(R)}(\sigma, X[i]) = 0$ *for* $i \in \{0, 1\}$,

*Proof.* We just prove (2); the ideas for (1) are similar (see also [11, Lemma 3.9(3)], [1, Lemma 3.4]).

Identify $\sigma$ with its cone:
$$Q_{-1} \xrightarrow{\sigma} Q_0 \to Z \to Q_{-1}[1]$$
and apply the functor $\mathrm{Hom}_{\mathbf{D}(R)}(-, X)$ to get the long exact sequence

$$(Q_{-1}[1], X) \to (Z, X) \to (Q_0, X) \xrightarrow{(\sigma, X)} (Q_{-1}, X) \to (Z[-1], X) \to (Q_0[-1], X).$$

The first item is zero since it is a negative extension of modules, and the last item is zero because $Q_0$ is projective. Thus, $\mathrm{Hom}_{\mathbf{D}(R)}(\sigma, X)$ is an isomorphism if and only if $\mathrm{Hom}_{\mathbf{D}(R)}(Z, X) = 0 = \mathrm{Hom}_{\mathbf{D}(R)}(Z[-1], X)$, and since $Z$ is quasi-isomorphic to $\sigma$, this is in turn equivalent to $\mathrm{Hom}_{\mathbf{D}(R)}(\sigma, X) = 0 = \mathrm{Hom}_{\mathbf{D}(R)}(\sigma, X[1])$. □

**Definition 1.20.** An $R$-module $T$ is

- *silting* (with respect to $\sigma$) if there exists a projective presentation $\sigma$ such that $\mathrm{Gen}(T) = \mathcal{D}_\sigma$ [11, Definition 3.10];

- *partial silting* (with respect to $\sigma$) if there exists a projective presentation $\sigma$ such that $T \in \mathcal{D}_\sigma$ and $\mathcal{D}_\sigma$ is closed under coproducts;

- *(partial) tilting* if it is (partial) silting with respect to a monomorphic projective presentation $\sigma$.

A module $T$ is tilting if and only if $\mathrm{Gen}(T) = \mathrm{Ker}\,\mathrm{Ext}^1_R(T, -)$. This is equivalent to the following set of conditions [26, Proposition 2.2]:

(1) $\mathrm{proj.}\dim T \leq 1$;

(2) $\mathrm{Ext}^1_R(T, T^{(J)}) = 0$ for any set $J$;

(3) for all $R$-modules $M$, if $\mathrm{Hom}_R(T, M) = 0 = \mathrm{Ext}^1_R(T, M)$, then $M = 0$.

Moreover, a tilting module is always a tilting object in $\mathbf{D}(R)$ according to Definition 1.10 (see [59, Example 4.2(iv)]). Conversely, we can prove the following.

**Lemma 1.21.** *Assume $T$ is a tilting object in the derived category $\mathbf{D}(R)$ of a module category over a ring $R$. If $T$ is concentrated in degree $0$ and $\mathrm{proj.}\dim T \leq 1$, then $T$ is a tilting $R$-module.*

*Proof.* We have to prove conditions (2) and (3) of the definition above. Indeed, $\mathrm{Ext}^1_R(T, T^{(J)}) = 0$ for any set $J$ by condition (T) of Definition 1.11; thus, (2) is satisfied. If $X$ is a module with $\mathrm{Hom}_R(T, X) = 0 = \mathrm{Ext}^1_R(T, X)$, from $\mathrm{proj.}\dim T \leq 1$ it follows that $\mathrm{Ext}^k_R(T, X)$ vanishes also for all $k \geq 2$. Thus, $\mathrm{Hom}_{\mathbf{D}(R)}(T, X[k]) = 0$ for all $k \in \mathbb{Z}$ and then $X = 0$ by (S3). □

Colpi and Trlifaj defined a partial tilting module as a module $T$ such that $\operatorname{Gen}(T) \subseteq \operatorname{Ker} \operatorname{Ext}^1_R(T,-)$ and $\operatorname{Gen}(T) \subseteq \operatorname{Ker} \operatorname{Ext}^1_R(T,-)$ is a torsion class [27, Definition 1.4]. It is easy to see that this definition is equivalent to ours. Indeed, if $T$ is partial silting with respect to a projective presentation $\sigma$, then $\mathcal{D}_\sigma$ is a torsion class if and only if it is closed under coproducts, and if $\sigma$ is injective, then $\mathcal{D}_\sigma = \operatorname{Ker} \operatorname{Ext}^1_R(T,-)$. Conversely, if $T$ is partial silting in the definition of Colpi and Trlifaj, then $\mathcal{D}_\sigma = \operatorname{Ker} \operatorname{Ext}^1_R(T,-)$ for a projective resolution $\sigma$ of $T$ (which has projective dimension at most 1) and thus $T$ is partial silting with respect to $\sigma$.

**Lemma 1.22.** *[11, Lemma 4.8] The following hold for a 2-term complex $\sigma$ in $\mathbf{K}^b(\operatorname{Proj} R)$ with $H^0(\sigma) = T$:*

(1) *An object $X$ in $D^{\leq 0}$ belongs to $\sigma^{\perp_{>0}}$ if and only if $H^0(X)$ lies in $\mathcal{D}_\sigma$.*

(2) *An object $X$ in $D^{\geq 0}$ belongs to $\sigma^{\perp_{\leq 0}}$ if and only if $H^0(X)$ lies in $T^{\perp_0}$.*

The following proposition generalizes Lemma 1.22 and describes the complexes lying in $\sigma^{\perp_{>0}}$ in terms of their cohomologies.

**Proposition 1.23.** *Let $\sigma$ be a two-term complex in $\mathbf{K}^b(\operatorname{Proj} R)$ with $T = H^0(\sigma)$. Then*
$$\sigma^{\perp_{>0}} = \{X \in \mathbf{D}(R) \mid H^0(X) \in \mathcal{D}_\sigma \text{ and } H^i(X) \in \mathcal{Y}_\sigma \text{ for all } i \geq 1\}.$$

*Proof.* In view of the equality $\mathcal{Y}_\sigma = \mathcal{D}_\sigma \cap T^{\perp_0}$, the characteristic property can be rewritten as
$$H^i(X) \in \mathcal{D}_\sigma \quad \text{and} \quad H^{i+1}(X) \in T^{\perp_0} \text{ for } i \geq 0.$$

For $i \geq 0$, consider the triangle
$$U_i \to X[i] \to V_i \to U_i[1], \tag{$\dagger$}$$

where $U_i \in D^{\leq 0}$ and $V_i \in D^{>0}$. Then $H^i(X) = H^0(U_i)$, $H^{i+1}(X) = H^0(V_i[1])$ and by Lemma 1.22, the claim amounts to showing that $X \in \sigma^{\perp_{>0}}$ if and only if
$$U_i \in \sigma^{\perp_{>0}} \text{ and } V_i[1] \in \sigma^{\perp_{\leq 0}}, \text{ for all } i \geq 0.$$

We first prove the implication "$\Leftarrow$".
Consider the triangle
$$U_i[1] \to X[i+1] \to V_i[1] \to U_i[2] \tag{$\star$}$$

for $i \geq 0$ and apply the functor $\operatorname{Hom}_{\mathbf{D}(R)}(\sigma,-)$. One gets
$$\operatorname{Hom}_{\mathbf{D}(R)}(\sigma, U_i[1]) \to \operatorname{Hom}_{\mathbf{D}(R)}(\sigma, X[i+1]) \to \operatorname{Hom}_{\mathbf{D}(R)}(\sigma, V_i[1]) \to .$$

Since the two outer terms are zero for $i \geq 0$, we have that $X \in \sigma^{\perp_{>0}}$.

"$\Rightarrow$". Assume now that $X \in \sigma^{\perp_{>0}}$. Since $\sigma$ is concentrated in degrees $-1$ and 0, we have $D^{\leq -2} \subseteq \sigma^{\perp_0}$ and $D^{>0} \subseteq \sigma^{\perp_0}$.

Consider the triangle $(\star)$ for $i \geq 0$; an application of $\operatorname{Hom}_{\mathbf{D}(R)}(\sigma,-)$ yields
$$\operatorname{Hom}_{\mathbf{D}(R)}(\sigma, X[i+1]) \to \operatorname{Hom}_{\mathbf{D}(R)}(\sigma, V_i[1]) \to \operatorname{Hom}_{\mathbf{D}(R)}(\sigma, U_i[2]) \to .$$

For $i \geq 0$, since $U_i[2] \in D^{\leq -2}$, the two outer terms are zero, hence the middle term is zero and $V_i[1] \in \sigma^{\perp_0}$. Moreover, since $V_i[1][-j] \in D^{\geq j} \subseteq \mathcal{D}^{>0}$ for $j > 0$, we conclude $V_i[1] \in \sigma^{\perp_{\leq 0}}$.

The second condition follows by analogous reasoning: in the short exact sequence

$$\mathrm{Hom}_{\mathbf{D}(R)}(\sigma, V_i) \to \mathrm{Hom}_{\mathbf{D}(R)}(\sigma, U_i[1]) \to \mathrm{Hom}_{\mathbf{D}(R)}(\sigma, X[i+1]) \to$$

the two outer terms vanish for $i \geq 0$. Thus, $U_i \in \sigma^{\perp_1}$; but since $U_i[j] \in D^{\leq -j} \subseteq D^{\leq -2}$ for $j \geq 2$, we also have $\mathrm{Hom}_{\mathbf{D}(R)}(\sigma, U_i[j]) = 0$ for $j \geq 2$. Hence, $U_i \in \sigma^{\perp_{>0}}$. □

**Corollary 1.24.** *Let $\sigma$ be a 2-term complex of projectives and $T = H^0(\sigma)$. If $\mathcal{D}_\sigma$ is closed under coproducts, then so are $\mathcal{Y}_\sigma$ and $\sigma^{\perp_{>0}}$.*

*Proof.* The fact that $\mathcal{Y}_\sigma$ is closed under coproducts follows from the equality $\mathcal{Y}_\sigma = \mathcal{D}_\sigma \cap T^{\perp_0}$. Thus, $\sigma^{\perp_{>0}}$ is closed under coproducts as a consequence of Proposition 1.23. □

**Proposition 1.25.** *The following hold for an $R$-module $T$.*

(1) *$T$ is silting if and only if there exists a 2-term silting complex $\sigma$ with $H^0(\sigma) = T$.*

(2) *$T$ is partial silting if and only if there exists a 2-term partial silting complex $\sigma$ with $H^0(\sigma) = T$.*

*Proof.* (1) follows from [11, Proposition 4.2, Theorem 4.11].

(2) Assume $\sigma$ is a 2-term partial silting complex of projectives with $T = H^0(\sigma)$. By the characterization of Lemma 1.22 and since $\sigma \in \sigma^{\perp_{>0}}$, it is immediate that $T \in \mathcal{D}_\sigma$. Moreover, again Lemma 1.22 implies that $\mathcal{D}_\sigma = \sigma^{\perp_{>0}} \cap \mathrm{Mod}\, R$ and as such it is closed under coproducts.

Conversely, assume $\mathcal{D}_\sigma$ is closed under coproducts and $T \in \mathcal{D}_\sigma$. By Corollary 1.24, $\sigma^{\perp_{>0}}$ is closed under coproducts. Finally, $H^0(\sigma) = T \in \mathcal{D}_\sigma$ by assumption and $H^i(\sigma) = 0 \in \mathcal{Y}_\sigma$ for $i \geq 1$, and hence $\sigma \in \sigma^{\perp_{>0}}$. □

## 1.7 Ring epimorphisms and universal localizations

Let now $\lambda \colon R \to S$ be a *ring epimorphism*, that is, an epimorphism in the category of rings. We say that two ring epimorphisms $f \colon R \to S$ and $g \colon R \to S'$ are *equivalent* if there is a ring isomorphism $h \colon S \to S'$ such that $g = h \circ f$. In this case we say that $f$ and $g$ lie in the same *epiclass* of $R$.

A particular case of ring epimorphism is given by the following theorem of Schofield's.

**Theorem 1.26.** *[66, Theorem 4.1] Let $R$ be a ring and $\Sigma$ a set of morphisms in $\mathrm{proj}\, R$. Then there exists a ring $R_\Sigma$ and ring homomorphism $f \colon R \to R_\Sigma$ such that*

(1) *$R_\Sigma \otimes_R \sigma$ is an isomorphism for all $\sigma \in \Sigma$ and*

(2) *every ring homomorphism $g \colon R \to S$ such that $S \otimes_R \sigma$ is an isomorphism for all $\sigma \in \Sigma$ factors uniquely through $f$.*

The ring $R_\Sigma$ of Theorem 1.26 is called the *universal localization* of $R$ at $\Sigma$. The homomorphism $f\colon R \to R_\Sigma$ is a ring epimorphism with $\operatorname{Tor}_1^R(R_\Sigma, R_\Sigma) = 0$ (see [66, Theorem 4.7]).

A subcategory $\mathcal{X}$ of $\operatorname{Mod} R$ is called *bireflective* if the inclusion functor admits both a left and a right adjoint. Bireflective subcategories are precisely those that are closed under products, coproducts, kernels and cokernels.

Define the following full subcategory of $\operatorname{Mod}(R)$:

$$\mathcal{X}_S = \text{the essential image of the restriction of}$$
$$\text{scalars functor } \lambda_*\colon \operatorname{Mod} S \to \operatorname{Mod} R.$$

In case the ring epimorphism $\lambda$ is a universal localization $R \to R_\Sigma$ we will write $\mathcal{X}_\Sigma$ instead of $\mathcal{X}_{R_\Sigma}$.

**Theorem 1.27.** (1) *[28, 30] The assignment $\Phi\colon \lambda \mapsto \mathcal{X}_S$ yields a bijection between epiclasses of ring epimorphisms $\lambda\colon R \to S$ and bireflective subcategories of $\operatorname{Mod} R$.*

(2) *[12, Proposition 3.3] If $T$ is a partial silting module with respect to $\sigma$, then the full subcategory $\mathcal{Y}_\sigma$ of $\operatorname{Mod} R$ is bireflective and extension-closed.*

As a consequence of Theorem 1.27 and Proposition 1.25, to any 2-term partial silting complex $\sigma$ in $\operatorname{Mod}(R)$ one can associate a ring epimorphism $\lambda\colon R \to S$ such that $\mathcal{Y}_\sigma = \mathcal{X}_S$. A ring epimorphism $\lambda\colon R \to S$ such that $\mathcal{Y}_\sigma = \mathcal{X}_S$ for some two-term partial silting complex $\sigma$ is called a *silting ring epimorphism*.

A ring epimorphism $\lambda\colon R \to S$ is *homological* if $\operatorname{Tor}_i^R(S, S) = 0$ for all $i > 0$. In this case there exists a recollement of triangulated categories

$$\mathbf{D}(S) \xrightarrow{\lambda_*} \mathbf{D}(R) \xrightarrow{j^*} \operatorname{Loc}(X) ,\qquad(1.5)$$

with $\lambda^* = -\otimes_R^{\mathbf{L}} S$, $\lambda^! = \mathbf{R}\operatorname{Hom}_R(S,-)$, and $j_!, j_*$,

where $\lambda_*$ is the derived functor of the restriction of scalars functor $\lambda_*\colon \operatorname{Mod} S \to \operatorname{Mod} R$, $X$ is the cone of the morphism $\lambda$ in $\mathbf{D}(R)$

$$R \xrightarrow{\lambda} S \longrightarrow X \longrightarrow R[1]$$

and $j_!$ is the embedding of $\operatorname{Loc}(X)$ into $\mathbf{D}(R)$ [6, §1.7].

# 2 Gluing silting objects

This chapter is devoted to the proof of the main theorem (Theorem 2.10) and its dual version (Theorem 2.16). We will study the compatibility of the gluing procedure with the gluing of t-structures and co-t-structures. We also restrict the gluing procedures to tilting objects and tilting modules, and then explore the particular case of a recollement induced by a silting ring epimorphism and a triangular matrix ring. Finally, we will compare our result to other gluing procedures appearing in [46] and [65].

Throughout this chapter, $\mathcal{D}$ will be a well generated triangulated category.

## 2.1 Universal maps

**Definition 2.1.** Let $\mathcal{X}$ be an additive subcategory of $\mathcal{D}$. For objects $M$ and $N$ in $\mathcal{D}$, a map $\alpha \in \operatorname{Hom}_{\mathcal{D}}(M, N)$ is called

(1) *left universal* if the map $\operatorname{Hom}_{\mathcal{D}}(M, \alpha)\colon \operatorname{End}_{\mathcal{D}}(M) \to \operatorname{Hom}_{\mathcal{D}}(M, N)$ is an epimorphism, that is, for any morphism $f\colon M \to N$ there exists a $g\colon M \to M$ such that the following diagram

$$\begin{array}{ccc} M & & \\ g\downarrow & \searrow f & \\ M & \xrightarrow{\alpha} & N \end{array}$$

commutes;

(2) an *$\mathcal{X}$-precover*, or a *right $\mathcal{X}$-approximation* if $M \in \mathcal{X}$ and for any object $X$ in $\mathcal{X}$, the map $\operatorname{Hom}_{\mathcal{D}}(X, \alpha)$ is an epimorphism, that is, for any morphism $f\colon X \to N$ there exists a $g\colon X \to M$ such that the following diagram

$$\begin{array}{ccc} X & & \\ g\downarrow & \searrow f & \\ M & \xrightarrow{\alpha} & N \end{array}$$

commutes.

Dually the map $\alpha$ is called

(1) *right universal* if the map $\mathrm{Hom}_{\mathcal{D}}(\alpha, N)\colon \mathrm{End}_{\mathcal{D}}(N) \to \mathrm{Hom}_{\mathcal{A}}(M, N)$ is an epimorphism;

(2) an *$\mathcal{X}$-preenvelope* or a *left $\mathcal{X}$-approximation* if $N \in \mathcal{X}$ and for any object $X$ in $\mathcal{X}$, the map $\mathrm{Hom}_{\mathcal{D}}(\alpha, X)$ is an epimorphism.

Note that a map $\alpha \in \mathrm{Hom}_{\mathcal{D}}(M, N)$ is left universal if and only if it is an $\mathrm{Add}(M)$-precover; it is right universal if and only if it is a $\mathrm{Prod}\, N$-preenvelope.

## 2.2   Some preliminary results

Inspired by [6], we consider the following conditions, for two objects $\sigma_1$ and $\sigma_2$ in $\mathcal{D}$ with $\sigma_1 \in \sigma_1^{\perp_{>0}}$ and $\sigma_2 \in \sigma_2^{\perp_{>0}}$:

(A1) $\mathrm{Hom}_{\mathcal{D}}(\sigma_1, \sigma_2[k]) = 0$ for all $k \geq 0$;

(A2) $\mathrm{Hom}_{\mathcal{D}}(\sigma_2, \sigma_1[k]) = 0$ for all $k \geq 2$.

The following results are adapted from [6]. Let $\sigma_1$ and $\sigma_2$ be two objects in $\mathcal{D}$ such that $\sigma_1 \in \sigma_1^{\perp_{>0}}$, $\sigma_2 \in \sigma_2^{\perp_{>0}}$ and conditions (A1) and (A2) are satisfied. Consider a morphism $\alpha\colon \sigma_2 \to \sigma_1[1]$ and let $\tilde{\sigma}$ be the object defined by the triangle

$$\sigma_1 \to \tilde{\sigma} \xrightarrow{\gamma} \sigma_2 \xrightarrow{\alpha} \sigma_1[1]. \qquad (2.1)$$

**Lemma 2.2.** *In the setting above, the following hold:*

- $\mathrm{Hom}_{\mathcal{D}}(\tilde{\sigma}, \sigma_2[k]) = 0$ *for any $k > 0$,*

- $\mathrm{Hom}_{\mathcal{D}}(\tilde{\sigma}, \sigma_1[k]) = 0$ *for any $k \geq 2$,*

- $\mathrm{Hom}_{\mathcal{D}}(\sigma_2, \tilde{\sigma}[k]) = 0$ *for any $k \geq 2$,*

- $\mathrm{Hom}_{\mathcal{D}}(\sigma_1, \tilde{\sigma}[k]) = 0$ *for any $k > 0$.*

*Proof.* For an integer number $k$, apply the functor $\mathrm{Hom}_{\mathcal{D}}(-, \sigma_2[k])$ to the triangle (2.1). The condition $\sigma_2 \in \sigma_2^{\perp_{>0}}$ and condition (A1) imply that $\mathrm{Hom}_{\mathcal{D}}(\tilde{\sigma}, \sigma_2[k]) = 0$ for any $k > 0$.

Next apply the functor $\mathrm{Hom}_{\mathcal{D}}(-, \sigma_1[k])$ to the same triangle. Condition (A2) and $\sigma_1 \in \sigma_1^{\perp_{>0}}$ imply that $\mathrm{Hom}_{\mathcal{D}}(\tilde{\sigma}, \sigma_1[k]) = 0$ for any $k \geq 2$.

Apply then $\mathrm{Hom}_{\mathcal{D}}(-, \sigma_2[k])$. By assumption, $\mathrm{Hom}_{\mathcal{D}}(\sigma_2, \sigma_2[k]) = 0$ for all positive integers $k$ and $\mathrm{Hom}_{\mathcal{D}}(\sigma_1, \sigma_2[k]) = 0$ for all $k \geq 0$. This implies that $\mathrm{Hom}_{\mathcal{D}}(\tilde{\sigma}, \sigma_2[k]) = 0$ for any $k > 0$.

Finally, apply $\mathrm{Hom}_{\mathcal{D}}(\sigma_2, -)$. Now $\mathrm{Hom}_{\mathcal{D}}(\sigma_2, \sigma_2[k]) = 0$ for any $k > 0$, while $\mathrm{Hom}_{\mathcal{D}}(\sigma_2, \sigma_1[k]) = 0$ for any $k \geq 2$. This implies $\mathrm{Hom}_{\mathcal{D}}(\sigma_2, \tilde{\sigma}[k]) = 0$ for any $k \geq 2$. $\square$

**Proposition 2.3.** *With the notation above, $\tilde{\sigma} \in \tilde{\sigma}^{\perp_{>0}}$ if and only if the homomorphism*

$$\varphi \colon \operatorname{Hom}_{\mathcal{D}}(\sigma_2, \sigma_2) \oplus \operatorname{Hom}_{\mathcal{D}}(\sigma_1[1], \sigma_1[1]) \to \operatorname{Hom}_{\mathcal{D}}(\sigma_2, \sigma_1[1])$$
$$(f, g) \mapsto \alpha \circ f + g \circ \alpha$$

*is surjective.*

*Remark* 2.4. Note that the homomorphism $\varphi$ above is surjective in particular if either

(i) the map $\operatorname{Hom}_{\mathcal{D}}(\sigma_2, \alpha)$ is surjective, or equivalently the map $\alpha$ is left universal, or

(ii) $\operatorname{Hom}_{\mathcal{D}}(\tilde{\sigma}, \sigma_1[1]) = 0$, equivalently, the map $\operatorname{Hom}_{\mathcal{D}}(\alpha, \sigma_1[1])$ is surjective, that is, the map $\alpha$ is right universal.

*Proof of Proposition 2.3.* By Lemma 2.2, $\operatorname{Hom}_{\mathcal{D}}(\tilde{\sigma}, \sigma_2[k]) = 0$ for $k > 0$ and $\operatorname{Hom}_{\mathcal{D}}(\tilde{\sigma}, \sigma_1[k]) = 0$ for $k \geq 2$.

Applying the functor $\operatorname{Hom}_{\mathcal{D}}(\tilde{\sigma}, -)$ to the triangle (2.1), we can see that $\operatorname{Hom}_{\mathcal{D}}(\tilde{\sigma}, \tilde{\sigma}[k]) = 0$ for any $k \geq 2$, and $\operatorname{Hom}_{\mathcal{D}}(\tilde{\sigma}, \tilde{\sigma}[1]) = 0$ if and only if the map

$$\operatorname{Hom}_{\mathcal{D}}(\tilde{\sigma}, \alpha) \colon \operatorname{Hom}_{\mathcal{D}}(\tilde{\sigma}, \sigma_2) \to \operatorname{Hom}_{\mathcal{D}}(\tilde{\sigma}, \sigma_1[1])$$

is surjective.

Finally, consider the following commutative diagram with exact columns:

$$\begin{array}{ccc}
0 & \longrightarrow & \operatorname{Hom}_{\mathcal{D}}(\sigma_1[1], \sigma_1[1]) \\
\downarrow & & \downarrow {\scriptstyle (\alpha, \sigma_1[1])} \\
\operatorname{Hom}_{\mathcal{D}}(\sigma_2, \sigma_2) & \xrightarrow{(\sigma_2, \alpha)} & \operatorname{Hom}_{\mathcal{D}}(\sigma_2, \sigma_1[1]) \\
\downarrow {\scriptstyle (\gamma, \sigma_2)} & & \downarrow {\scriptstyle (\gamma, \sigma_1[1])} \\
\operatorname{Hom}_{\mathcal{D}}(\tilde{\sigma}, \sigma_2) & \xrightarrow{(\tilde{\sigma}, \alpha)} & \operatorname{Hom}_{\mathcal{D}}(\tilde{\sigma}, \sigma_1[1]) \\
\downarrow & & \downarrow \\
0 & \longrightarrow & 0.
\end{array}$$

One proves now easily by diagram chasing that the map $(\tilde{\sigma}, \alpha)$ being surjective is equivalent to the surjectivity of the map

$$\operatorname{Hom}_{\mathcal{D}}(\sigma_2, \sigma_2) \oplus \operatorname{Hom}_{\mathcal{D}}(\sigma_1[1], \sigma_1[1]) \longrightarrow \operatorname{Hom}_{\mathcal{D}}(\sigma_2, \sigma_1[1])$$
$$(f, g) \longmapsto \alpha \circ f + g \circ \alpha.$$

$\square$

**Proposition 2.5.** *Let $\sigma_1$ and $\sigma_2$ be two objects in $\mathcal{D}$ satisfying $\sigma_1 \in \sigma_1^{\perp_{>0}}$, $\sigma_2 \in \sigma_2^{\perp_{>0}}$ and conditions (A1) and (A2) and let $\tilde{\sigma}$ be constructed as in Proposition 2.3. Then the following hold*

(1) *the object $\tilde{\sigma} \oplus \sigma_2 \in (\tilde{\sigma} \oplus \sigma_2)^{\perp_{>0}}$ if and only if the morphism $\alpha \colon \sigma_2 \to \sigma_1[1]$ is left universal;*

(2) *the object $\tilde{\sigma} \oplus \sigma_1 \in (\tilde{\sigma} \oplus \sigma_1)^{\perp_{>0}}$ if and only if the morphism $\alpha \colon \sigma_2 \to \sigma_1[1]$ is right universal.*

*Proof.* (1) By Lemma 2.2, $\mathrm{Hom}_{\mathcal{D}}(\sigma_2, \tilde{\sigma}[k]) = 0$ for $k \geq 2$ and $\mathrm{Hom}_{\mathcal{D}}(\tilde{\sigma}, \sigma_2[k]) = 0$, for $k > 0$. Applying $\mathrm{Hom}_{\mathcal{D}}(\sigma_2, -)$ to the triangle (2.1), we obtain

$$\mathrm{Hom}_{\mathcal{D}}(\sigma_2, \sigma_2) \overset{(\sigma_2, \alpha)}{\longrightarrow} \mathrm{Hom}_{\mathcal{D}}(\sigma_2, \sigma_1[1]) \to \mathrm{Hom}_{\mathcal{D}}(\sigma_2, \tilde{\sigma}[1]) \to 0.$$

Then $\mathrm{Hom}_{\mathcal{D}}(\sigma_2, \tilde{\sigma}[1]) = 0$ is equivalent to the map $(\sigma_2, \alpha)$ being surjective, i.e., $\alpha$ being left universal.

(2) is proved in a similar way. By Lemma 2.2, $\mathrm{Hom}_{\mathcal{D}}(\sigma_1, \tilde{\sigma}[k]) = 0$ for $k > 0$ and $\mathrm{Hom}_{\mathcal{D}}(\tilde{\sigma}, \sigma_1[k]) = 0$, for $k \geq 2$.

Applying $\mathrm{Hom}_{\mathcal{D}}(-, \sigma_1)$ to the triangle (2.1), we obtain

$$\mathrm{Hom}_{\mathcal{D}}(\sigma_1, \sigma_1) \overset{(\alpha, \sigma_1)}{\longrightarrow} \mathrm{Hom}_{\mathcal{D}}(\sigma_2, \sigma_1[1]) \to \mathrm{Hom}_{\mathcal{D}}(\tilde{\sigma}, \sigma_1[1]) \to 0.$$

Then $\mathrm{Hom}_{\mathcal{D}}(\tilde{\sigma}, \sigma_1[1]) = 0$ is equivalent to the map $(\alpha, \sigma_1)$ being surjective, i.e., $\alpha$ being right universal.

Finally, use Remark 2.4 and Proposition 2.3 to see that in both cases (1) and (2), $\mathrm{Hom}_{\mathcal{D}}(\tilde{\sigma}, \tilde{\sigma}[k]) = 0$ for all $k > 0$. □

## 2.3 Adjacent recollements

Consider a recollement

$$\mathcal{Y} \underset{\underset{i^!}{\longleftarrow}}{\overset{\overset{i^*}{\longleftarrow}}{\underset{i_*}{\longrightarrow}}} \mathcal{D} \underset{\underset{j_*}{\longleftarrow}}{\overset{\overset{j_!}{\longleftarrow}}{\underset{j^*}{\longrightarrow}}} \mathcal{X} \tag{2.2}$$

of triangulated categories, let $\mathcal{D}$ be well generated and let $\sigma$ be a silting object in $\mathcal{Y}$. Notice that in this case $\mathcal{X}$ and $\mathcal{Y}$ are also well generated and all three categories satisfy Brown representability. We investigate equivalent conditions for the class $(i_*\sigma)^{\perp_{>0}}$ to be closed under coproducts in $\mathcal{D}$. This property will be important in the following; however, it is quite a restrictive condition to assume, in that it turns out to be equivalent to the given recollement having a lower adjacent recollement[1].

We say a recollement of triangulated categories is *lower adjacent* to (2.2) if it is of the form

$$\mathcal{X} \underset{\underset{j^\#}{\longleftarrow}}{\overset{\overset{j^*}{\longleftarrow}}{\underset{j_*}{\longrightarrow}}} \mathcal{D} \underset{\underset{i_\#}{\longleftarrow}}{\overset{\overset{i_*}{\longleftarrow}}{\underset{i^!}{\longrightarrow}}} \mathcal{Y}.$$

**Lemma 2.6.** *Consider a recollement of well generated triangulated categories of the form (2.2) and let $\sigma$ be a silting object in $\mathcal{Y}$. If the class $(i_*\sigma)^{\perp_{>0}}$ is closed under coproducts in $\mathcal{D}$, then so is also the class $\mathrm{Im}(j_*)$.*

*Proof.* The TTF triple $(\mathrm{Ker}\, i^* = \mathrm{Im}\, j_!, \mathrm{Im}\, i_*, \mathrm{Im}\, j_* = \mathrm{Ker}\, i^!)$ associated to the given recollement shows that the subcategory $\mathrm{Im}\, j_*$ of $\mathcal{D}$ coincides with the perpendicular subcategory $(i_*\mathcal{Y})^{\perp_0}$, which in turn can be written as $(i_*\sigma)^{\perp_{\mathbb{Z}}}$, since $\sigma$ generates $\mathcal{Y}$. By noticing that

$$(i_*\sigma)^{\perp_{\mathbb{Z}}} = \bigcap_{k \in \mathbb{Z}} (i_*\sigma)^{\perp_{>0}}[k],$$

we have the claim. $\square$

---

[1] We are grateful to Jorge Vitória for pointing out these equivalences.

**Lemma 2.7.** *[53, Theorem 8.4.4] Let $F\colon \mathcal{Y} \to \mathcal{D}$ be a (covariant) functor between triangulated categories with coproducts. Assume $\mathcal{Y}$ satisfies Brown representability. Then $F$ has a right adjoint if and only if it preserves coproducts.*

**Lemma 2.8.** *If $\operatorname{Im} j_*$ is closed under coproducts, then $j_*$ preserves coproducts.*

*Proof.* Consider a recollement of triangulated categories with coproducts of the form (2.2). Assume $\operatorname{Im} j_*$ is closed under coproducts and consider the torsion pair $(\operatorname{Im} i_*, \operatorname{Im} j_*)$ in $\mathcal{D}$ induced by the recollement. By [18, Lemma III.1.2], the functor $i^!$ preserves coproducts. We have then the following isomorphism of triangles, for any family $(Z_j)_{j \in J}$ of objects of $\mathcal{D}$:

$$\begin{array}{ccccccc}
i_*i^!\bigoplus_{j\in J} Z_j & \longrightarrow & \bigoplus_{j\in J} Z_j & \longrightarrow & j_*j^*\bigoplus_{j\in J} Z_j & \longrightarrow & i_*i^!\bigoplus_{j\in J} Z_j[1] \\
\downarrow{\scriptstyle \mathbb{R}} & & \parallel & & \downarrow{\scriptstyle \mathbb{R}} & & \downarrow{\scriptstyle \mathbb{R}} \\
\bigoplus_{j\in J} i_*i^! Z_j & \longrightarrow & \bigoplus_{j\in J} Z_j & \longrightarrow & \bigoplus_{j\in J} j_*j^* Z_j & \longrightarrow & \bigoplus_{j\in J} i_*i^! Z_j[1].
\end{array}$$

Thus,

$$\bigoplus_{j\in J} j_*j^* Z_j \cong j_*j^* \bigoplus_{j\in J} Z_j \cong j_* \bigoplus_{j\in J} j^* Z_j,$$

where the second isomorphism follows from the fact that $j^*$ has a right adjoint and thus it commutes with coproducts. Finally, the identity

$$\bigoplus_{j\in J} j_* X_j \cong j_* \bigoplus_{j\in J} X_j,$$

holds true for any family $(X_j)_{j \in J}$ of objects in $\mathcal{X}$ since $j^*$ is dense. $\square$

**Proposition 2.9.** *Consider a recollement of well generated triangulated categories of the form (2.2) and let $\sigma$ be a silting object in $\mathcal{Y}$. The following statements are equivalent:*

(1) *The class $(i_*\sigma)^{\perp_{>0}}$ is closed under coproducts;*

(2) *The functor $j_*$ has a right adjoint;*

(3) *The functor $i^!$ has a right adjoint;*

(4) *The recollement (2.2) has a lower adjacent recollement.*

*Proof.* The equivalence between (2), (3) and (4) follows from [25, Theorem 2.1] and its dual. We prove that (1) implies (2) and (3) implies (1).

(1) $\Rightarrow$ (2). Assume that $(i_*\sigma)^{\perp_{>0}}$ is closed under coproducts. By Lemma 2.6, the class $\operatorname{Im} j_*$ is closed under coproducts. Since $\mathcal{Y}$ is well generated (Theorem 1.15) and well generated triangulated categories satisfy Brown representability, Lemma 2.7 and Lemma 2.8 give the claim.

(3) $\Rightarrow$ (1). Assume the functor $i^!$ has a right adjoint and consider a family $(X_j)_{j \in J}$ of objects in $\mathcal{D}$ such that for each $j \in J$, $X_j \in (i_*\sigma)^{\perp_{>0}}$. By adjunction, this implies $i^! X_j \in \sigma^{\perp_{>0}}$ for any $j \in J$. Then, for each $k > 0$,

$$\operatorname{Hom}_{\mathcal{D}}(i_*\sigma, \bigoplus_{j \in J} X_j[k]) \cong \operatorname{Hom}_{\mathcal{Y}}(\sigma, i^! \bigoplus_{j \in J} X_j[k])$$
$$\cong \operatorname{Hom}_{\mathcal{Y}}(\sigma, \bigoplus_{j \in J} i^! X_j[k]) = 0,$$

where the last isomorphism holds true since $i^!$ has a right adjoint and the last term is zero since $\sigma^{\perp_{>0}}$ is closed under coproducts. $\square$

## 2.4 Recollements generated by partial silting objects

Throughout this section, let $\mathcal{D}$ be a well generated triangulated category. Suppose there is a partial silting object $\omega$ in $\mathcal{D}$ and a silting object $\sigma$ in $\omega^{\perp_{\mathbb{Z}}}$. We can construct a recollement associated to $\omega$, as in Section 1.4:

$$\omega^{\perp_{\mathbb{Z}}} \xrightarrow[\substack{\leftarrow i^* \\ \leftarrow i^!}]{i_*} \mathcal{D} \xrightarrow[\substack{\leftarrow j_! \\ \leftarrow j_*}]{j^*} \operatorname{Loc}(\omega) \ , \tag{2.3}$$

where $\omega$ and $\sigma$ are silting respectively in $\operatorname{Loc}(\omega)$ and in $\omega^{\perp_{\mathbb{Z}}}$. In view of the recollement (2.3), we will write $j_!\omega$ and $i_*\sigma$ when these objects are considered as objects of $\mathcal{D}$. We are now ready to state the main result.

**Theorem 2.10.** *Consider a recollement of the form (2.3) associated to a partial silting object $\omega$ and let $\sigma$ be a silting object in $\omega^{\perp_{\mathbb{Z}}}$. Assume the following additional hypotheses are satisfied:*

(1) *the class $(i_*\sigma)^{\perp_{>0}} \cap (j_!\omega)^{\perp_{>0}}$ is closed under coproducts in $\mathcal{D}$,*

(2) $\operatorname{Hom}_{\mathcal{D}}(i_*\sigma, j_!\omega[k]) = 0$ *for all integers $k \geq 2$.*

*Consider the set $I = \operatorname{Hom}_{\mathcal{D}}(i_*\sigma, j_!\omega[1])$, let $\alpha \colon i_*\sigma^{(I)} \to j_!\omega[1]$ be the canonically induced map and $\tilde{\sigma}$ be the cocone of $\alpha$:*

$$j_!\omega \to \tilde{\sigma} \to i_*\sigma^{(I)} \xrightarrow{\alpha} j_!\omega[1]. \tag{2.4}$$

*Then $\tilde{\sigma} \oplus i_*\sigma$ is a silting object in $\mathcal{D}$.*

*Proof.* We prove that the object $\tilde{\sigma} \oplus i_*\sigma$ is silting. First we prove (S1), that is, $\operatorname{Hom}_{\mathcal{D}}(\tilde{\sigma} \oplus i_*\sigma, (\tilde{\sigma} \oplus i_*\sigma)[k]) = 0$ for $k > 0$. We show that the hypotheses of Proposition 2.3 are satisfied for the objects $\sigma_1 = j_!\omega$ and $\sigma_2 = i_*\sigma^{(I)}$ in $\mathcal{D}$. The objects $j_!\omega$ and $i_*\sigma^{(I)}$ clearly satisfy $j_!\omega \in (j_!\omega)^{\perp_{>0}}$ and $i_*\sigma^{(I)} \in (i_*\sigma^{(I)})^{\perp_{>0}}$, and condition (A2) holds by hypothesis (2). We are left to prove (A1), that is, $\operatorname{Hom}_{\mathcal{D}}(j_!\omega, i_*\sigma^{(I)}[k]) = 0$ for all $k \geq 0$. Indeed, for any $k \in \mathbb{Z}$, $\mathrm{Hom}_{\mathcal{D}}(j_!\omega, i_*\sigma^{(I)}[k]) \cong \mathrm{Hom}_{\mathrm{Loc}(\omega)}(\omega, j^*i_*\sigma^{(I)}[k]) = 0$ since $j^*i_* = 0$. Finally, the map $\alpha$ is left universal by construction, and thus Proposition 2.5 gives the claim.

Next, we prove that $\tilde{\sigma} \oplus i_*\sigma$ is partial silting, which amounts to prove that $(\tilde{\sigma} \oplus i_*\sigma)^{\perp_{>0}}$ is closed under coproducts. Consider the long exact sequence

$$\mathrm{Hom}_{\mathcal{D}}(i_*\sigma^{(I)}, X[k]) \longrightarrow \mathrm{Hom}_{\mathcal{D}}(\tilde{\sigma}, X[k]) \longrightarrow \mathrm{Hom}_{\mathcal{D}}(j_!\omega, X[k]) \longrightarrow \mathrm{Hom}_{\mathcal{D}}(i_*\sigma^{(I)}, X[k+1])$$

for some $X$ in $\mathcal{D}$ and $k > 0$: it is easy to see that $X$ belongs to $(i_*\sigma \oplus \tilde{\sigma})^{\perp_{>0}}$ if and only if it belongs to $(i_*\sigma \oplus j_!\omega)^{\perp_{>0}}$, which coincides with $(i_*\sigma)^{\perp_{>0}} \cap (j_!\omega)^{\perp_{>0}}$. The latter class is closed under coproducts by hypothesis (2).

Finally, we prove the generating property (S3). Notice that $\tilde{\sigma} \oplus i_*\sigma$ is a generator if and only if so is $i_*\sigma \oplus j_!\omega$. Assume then $\mathrm{Hom}_{\mathcal{D}}(i_*\sigma \oplus j_!\omega, X[k]) = 0$ for all $k \in \mathbb{Z}$. In particular, $0 = \mathrm{Hom}_{\mathcal{D}}(j_!\omega, X[k]) \cong \mathrm{Hom}_{\mathrm{Loc}(\omega)}(\omega, j^*X[k])$ for all $k \in \mathbb{Z}$, but since $\omega$ is silting in $\mathrm{Loc}(\omega)$ (Proposition 1.14), it follows $j^*X = 0$. Consider then the triangle

$$i_*i^*X \to X \to j_!j^*X \to i_*i^*X[1]$$

induced by the recollement. Since $j^*X = 0$, the object $X$ is isomorphic to $i_*i^*X$. Thus

$$0 = \mathrm{Hom}_{\mathcal{D}}(i_*\sigma, X[k]) \cong \mathrm{Hom}_{\mathcal{D}}(i_*\sigma, i_*i^*X[k]) \cong \mathrm{Hom}_{\omega^{\perp_{\mathbb{Z}}}}(\sigma, i^*X[k])$$

and since $\sigma$ is silting in $\omega^{\perp_{\mathbb{Z}}}$, it follows $i^*X = 0$. Thus $X = 0$. □

*Remark* 2.11. (1) A sufficient for hypothesis (1) of Theorem 2.10 to hold is that given recollement has a lower adjacent recollement. To see this, recall that it is equivalent to the class $(i_*\sigma)^{\perp_{>0}}$ being closed under coproducts. Moreover, $(j_!\omega)^{\perp_{>0}}$ is always closed under coproducts: consider indeed a class $(X_j)_{j \in J}$ of objects lying in $(j_!\omega)^{\perp_{>0}}$. By adjunction, this implies $j^*X_j \in \omega^{\perp_{>0}}$ for every $j \in J$. Then, for each $k > 0$,

$$\mathrm{Hom}_{\mathcal{D}}(j_!\omega, \bigoplus_{j \in J} X_j[k]) \cong \mathrm{Hom}_{\mathcal{X}}(\omega, j^* \bigoplus_{j \in J} X_j[k])$$
$$\cong \mathrm{Hom}_{\mathcal{X}}(\omega, \bigoplus_{j \in J} j^*X_j[k]) = 0,$$

where the last isomorphism holds true since $j^*$ has a right adjoint, and the last term is zero since $\omega^{\perp_{>0}}$ is closed under coproducts.

Thus, the intersection $(i_*\sigma)^{\perp_{>0}} \cap (j_!\omega)^{\perp_{>0}}$ is closed under coproducts and hypothesis (1) is satisfied.

(2) With reference to point (1), combining [7, Proposition 3.2 and Lemma 2.5] we get a handy criterion to tell if a recollement can be extended one step downwards. Namely, if $A$, $B$, $R$ are $k$-algebras and

$$\mathbf{D}(B) \xrightarrow[\substack{\longleftarrow i^* \\ \longleftarrow i_* \longrightarrow \\ \longleftarrow i^!}]{} \mathbf{D}(R) \xrightarrow[\substack{\longleftarrow j_! \\ \longleftarrow j^* \longrightarrow \\ \longleftarrow j_*}]{} \mathbf{D}(A)$$

is a recollement of triangulated categories, then the recollement can be extended one step downwards if and only if $i_*B \in \mathbf{K}^b(\mathrm{proj}\,R)$, equivalently if $j^!A \in \mathbf{K}^b(\mathrm{proj}\,R)$. Dually, there is a similar criterion for such a recollement to be extended one step upwards. A generalized version of this criterion is available in [65].

Next we determine conditions under which Theorem 2.10 restricts to tilting objects.

**Corollary 2.12.** *In the hypotheses of Theorem 2.10, assume $\sigma$ and $\omega$ are tilting objects respectively in $\omega^{\perp_{\mathbb{Z}}}$ and in $\mathrm{Loc}(\omega)$. Assume moreover the following additional condition is satisfied:*

(3) $\mathrm{Hom}_{\mathcal{D}}(i_*\sigma, j_!\omega^{(J)}[k]) = 0$ *for $k < 0$, for any set $J$.*

*Then the object $\tilde{\sigma} \oplus i_*\sigma$ is tilting in $\mathcal{D}$.*

*Assume, in addition, that $\mathcal{D} = \mathbf{D}(R)$ for a ring $R$, and that $\sigma$ and $\omega$ have cohomologies concentrated in degree $0$ and have projective dimension, as $R$-modules, at most $1$. Then $\tilde{\sigma} \oplus i_*\sigma$ is a tilting $R$-module.*

*Proof.* We have to check the following conditions hold for any set $J$ and $k < 0$:

(i) $\mathrm{Hom}_{\mathcal{D}}(i_*\sigma, i_*\sigma^{(J)}[k]) = 0$,

(ii) $\mathrm{Hom}_{\mathcal{D}}(\tilde{\sigma}, i_*\sigma^{(J)}[k]) = 0$,

(iii) $\mathrm{Hom}_{\mathcal{D}}(i_*\sigma, \tilde{\sigma}^{(J)}[k]) = 0$,

(iv) $\mathrm{Hom}_{\mathcal{D}}(\tilde{\sigma}, \tilde{\sigma}^{(J)}[k]) = 0$.

Condition (i) holds because $i_*$ is fully faithful and since $\sigma$ is tilting in $\omega^{\perp_{\mathbb{Z}}}$. For (ii), apply the functor $\mathrm{Hom}_{\mathcal{D}}(-, i_*\sigma^{(J)})$, to the triangle

$$j_!\omega \to \tilde{\sigma} \to i_*\sigma^{(I)} \to j_!\omega[1].$$

Moreover, $\mathrm{Hom}_{\mathcal{D}}(\tilde{\sigma}, j_!\omega^{(J)}[k]) = 0$, for $k < 0$, follows by an application of $\mathrm{Hom}_{\mathbf{D}(R)}(-, j_!\omega^{(J)})$ to the same triangle.

For conditions (iii) and (iv) apply respectively the functors $\mathrm{Hom}_{\mathcal{D}}(i_*\sigma, -)$ and $\mathrm{Hom}_{\mathcal{D}}(\tilde{\sigma}, -)$ to the triangle

$$j_!\omega^{(J)} \to \tilde{\sigma}^{(J)} \to i_*\sigma^{(I\times J)} \to j_!\omega^{(J)}[1].$$

This proves that $\tilde{\sigma} \oplus i_*\sigma$ is a tilting object.

In case $\mathcal{D} = \mathbf{D}(R)$ and $\sigma$ and $\omega$ have cohomologies concentrated in degree $0$, then the triangle (2.4) is induced by a short exact sequence

$$0 \to \omega \to \tilde{\sigma} \to \sigma^{(I)} \to 0$$

in $\mathrm{Mod}\, R$. In particular, $\tilde{\sigma}$ is an $R$-module and the short exact sequence implies $\tilde{\sigma}$, and hence $\tilde{\sigma} \oplus i_*\sigma$, has projective dimension at most $1$. By Lemma 1.21 we conclude that $\tilde{\sigma} \oplus i_*\sigma$ is a tilting module. □

*Remark* 2.13. (1) Corollary 2.12 generalizes [6, Theorem 2.4]. Indeed, let $T_1 = j_!\omega$ and $T_2 = i_*\sigma$; if $T_2$ is compact, then hypothesis (1) of Theorem 2.10 is satisfied since $T_2^{\perp_{>0}}$ is closed under coproducts.

(2) In case $\mathcal{D} = \mathbf{D}(R)$ and $\sigma$ and $\omega$ have cohomologies concentrated in degree $0$, then the triangle (2.4) is induced from a short exact sequence

$$0 \to U \to \tilde{T} \to T^{(I)} \to 0$$

in $\mathrm{Mod}\, R$, where $U = H^0(j_!\omega)$, $\tilde{T} = H^0(\tilde{\sigma})$ and $T = H^0(i_*\sigma)$. This sequence coincides with the Bongartz sequence as generalized by Trlifaj ([68, Lemma 6.8]).

## 2.5 The gluing preserves co-t-structures

We show now that the method presented above is a gluing of silting objects with respect to co-t-structures, meaning that the co-t-structure generated by $\tilde{\sigma} \oplus i_*\sigma$ coincides with the co-t-structure obtained by gluing the co-t-structures generated by $\sigma$ and $\omega$.

Consider the co-t-structures $(\mathcal{Y}', \mathcal{Y}'') = (^{\perp_0}(\sigma^{\perp_{\geq 0}}), \sigma^{\perp_{\geq 0}})$ and $(\mathcal{X}', \mathcal{X}'') = (^{\perp_0}(\omega^{\perp_{\geq 0}}), \omega^{\perp_{\geq 0}})$ respectively in $\omega^{\perp_{\mathbb{Z}}}$ and in $\operatorname{Loc}(\omega)$, associated to the silting objects $\sigma$ and $\omega$ respectively. Denote by $\rho = \tilde{\sigma} \oplus i_*\sigma$ the silting object constructed as in Theorem 2.10. Consider a recollement of the form (2.3) and assume $\mathcal{D}$ is well generated and $\mathcal{D}$, $\operatorname{Loc}(\omega)$ and $\omega^{\perp_{\mathbb{Z}}}$ satisfy dual Brown representability (for example this holds if $\mathcal{D}$ is compactly generated and there is a lower adjacent recollement to the given one). In this case one can associate a co-t-structure to any silting objects in the categories forming the recollement. Consider the following co-t-structures in $\mathcal{D}$:

- $(^{\perp_0}(\rho^{\perp_{\geq 0}}), \rho^{\perp_{\geq 0}})$ associated to the silting object $\rho$ and

- $(\mathcal{D}', \mathcal{D}'')$ obtained by gluing the two co-t-structures $(\mathcal{Y}', \mathcal{Y}'')$ and $(\mathcal{X}', \mathcal{X}'')$ as in Proposition 1.3.

**Theorem 2.14.** *In the situation above, the co-t-structures* $(^{\perp_0}(\rho^{\perp_{\geq 0}}), \rho^{\perp_{\geq 0}})$ *and* $(\mathcal{D}', \mathcal{D}'')$ *coincide.*

*Proof.* Consider the coaisle $\mathcal{D}''$ of the glued co-t-structure:

$$\mathcal{D}'' = \left\{ Z \in \mathcal{D} \mid i^! Z \in \sigma^{\perp_{\geq 0}}, j^* Z \in \omega^{\perp_{\geq 0}} \right\}.$$

The first condition is equivalent to $\operatorname{Hom}_{\mathcal{D}}(i_*\sigma, Z[k]) = 0$ for $k \geq 0$, the second to $\operatorname{Hom}_{\mathcal{D}}(j_!\omega, Z[k]) = 0$ for $k \geq 0$. From the long exact sequence

$$\operatorname{Hom}_{\mathcal{D}}(i_*\sigma^{(I)}, Z[k]) \to \operatorname{Hom}_{\mathcal{D}}(\tilde{\sigma}, Z[k]) \to \operatorname{Hom}_{\mathcal{D}}(j_!\omega, Z[k]) \to \operatorname{Hom}_{\mathcal{D}}(i_*\sigma^{(I)}, Z[k+1])$$

it is clear that the condition $\operatorname{Hom}_{\mathcal{D}}(j_!\omega, Z[k]) = 0 = \operatorname{Hom}_{\mathcal{D}}(i_*\sigma, Z[k])$ for all $k \geq 0$ is equivalent to $\operatorname{Hom}_{\mathcal{D}}(\tilde{\sigma}, Z[k]) = 0 = \operatorname{Hom}_{\mathcal{D}}(i_*\sigma, Z[k])$ for all $k \geq 0$. Thus $\mathcal{D}''$ coincides with

$$\rho^{\perp_{\geq 0}} = \{Z \in \mathcal{D} \mid \operatorname{Hom}_{\mathcal{D}}(i_*\sigma, Z[k]) = 0, \operatorname{Hom}_{\mathcal{D}}(\tilde{\sigma}, Z[k]) = 0 \text{ for } k \geq 0\}.$$

Now the aisles of the two co-t-structures coincide since they are left Hom-perpendicular to the coaisles, and these coincide. □

## 2.6 A 'dual' version

In this section we prove a 'dual' version of Theorem 2.10, in the sense that we consider *right* universal maps and exploit the second statement of Proposition 2.5. Such a version requires however some quite strong additional hypotheses.

**Lemma 2.15.** *Let $M$ be a compact object in a triangulated category $\mathcal{D}$ with coproducts and $\alpha \colon M \to N$ be a morphism in $\mathcal{D}$. Let $\{X_j\}_{j \in J}$ be a set of objects in $\mathcal{D}$ and assume $\operatorname{Hom}_{\mathcal{D}}(\alpha, X_j)$ is surjective for all $j \in J$. Then $\operatorname{Hom}_{\mathcal{D}}(\alpha, \bigoplus_{j \in J} X_j)$ is surjective.*

*Proof.* Consider the following commutative diagram of abelian groups:

$$\begin{array}{ccc} \bigoplus_{j\in J}\operatorname{Hom}_{\mathcal{D}}(N,X_j) & \xrightarrow{\oplus(\alpha,X_j)} & \bigoplus_{j\in J}\operatorname{Hom}_{\mathcal{D}}(M,X_j) \\ \downarrow & & \downarrow \wr \\ \operatorname{Hom}_{\mathcal{D}}(N,\bigoplus_{j\in J}X_j) & \xrightarrow{(\alpha,\oplus X_j)} & \operatorname{Hom}_{\mathcal{D}}(M,\bigoplus_{j\in J}X_j), \end{array}$$

where the vertical arrows are the universal morphisms induced by the respective coproducts of Hom groups. The vertical map on the right is an isomorphism since $M$ is compact, and the map $\bigoplus(\alpha, X_j)$ is surjective since the coproduct of epimorphisms is an epimorphism. Thus, the map $(\alpha, \bigoplus X_j)$ is surjective. $\square$

**Theorem 2.16.** *Consider a recollement of the form (2.3) associated to the partial silting object $\omega$ and let $\sigma$ be a silting object in $\omega^{\perp_\mathbb{Z}}$. Assume $j_*\omega$ has an $\operatorname{add}(i_*\sigma)$-preenvelope $\alpha\colon j_*\omega \to i_*\sigma^I[1]$, for a finite index set $I$, and assume the following additional hypotheses are satisfied:*

(1) *the class $(i_*\sigma)^{\perp_{>0}}$ is closed under coproducts in $\mathcal{D}$,*

(2) $\operatorname{Hom}_{\mathcal{D}}(j_*\omega, i_*\sigma[k]) = 0$ *for all integers $k \geq 2$,*

(3) $j_*\omega$ *is a compact object in $\mathcal{D}$.*

*Let $\tilde{\sigma}$ be the cocone of the map $\alpha$:*

$$i_*\sigma^I \to \tilde{\sigma} \to j_*\omega \xrightarrow{\alpha} i_*\sigma^I[1]. \tag{2.4'}$$

*Then $\tilde{\sigma} \oplus i_*\sigma$ is a silting object in $\mathcal{D}$.*

*Proof.* We prove that the hypotheses (S1)–(S3) of Definition 1.11 are satisfied for the object $\tilde{\sigma} \oplus i_*\sigma$. To prove (S1), that is, $\operatorname{Hom}_{\mathcal{D}}(\tilde{\sigma} \oplus i_*\sigma, (\tilde{\sigma} \oplus i_*\sigma)[k]) = 0$ for $k > 0$, we show that the hypotheses of Proposition 2.5 are satisfied for the objects $\sigma_1 = (i_*\sigma)^I$ and $\sigma_2 = j_*\omega$ in $\mathcal{D}$. This means we have to prove the following four conditions:

- $i_*\sigma^I \in (i_*\sigma^I)^{\perp_{>0}}$,

- $j_*\omega \in (j_*\omega)^{\perp_{>0}}$,

(A1) $\operatorname{Hom}_{\mathcal{D}}(i_*\sigma^I, j_*\omega[k]) = 0$ for all $k \geq 0$,

(A2) $\operatorname{Hom}_{\mathcal{D}}(j_*\omega, i_*\sigma^I[k]) = 0$ for all $k \geq 2$.

The object $j_*\omega$ clearly satisfies $j_*\omega \in (j_*\omega)^{\perp_{>0}}$; moreover,

$$\operatorname{Hom}_{\mathcal{D}}(i_*\sigma^I, i_*\sigma^I[k]) \cong \operatorname{Hom}_{\mathcal{D}}(i_*\sigma, i_*\sigma)^{I \times I}$$
$$\cong \operatorname{Hom}_{\mathcal{D}}(\sigma, \sigma[k]) = 0$$

for $k > 0$, where the first isomorphism follows since $I$ is finite. Next we prove condition (A2), that is,

$$\operatorname{Hom}_{\mathcal{D}}(j_*\omega, i_*\sigma^I[k]) = 0 \text{ for } k \geq 2.$$

This holds by hypothesis (2) since the covariant Hom functor commutes with products. We are left to prove (A1), that is,

$$\mathrm{Hom}_{\mathcal{D}}(i_*\sigma^I, j_*\omega[k]) = 0 \text{ for all } k \geq 0.$$

Indeed, for any $k \in \mathbb{Z}$, $\mathrm{Hom}_{\mathcal{D}}(i_*\sigma^I, j_*\omega[k]) \cong \mathrm{Hom}_{\mathrm{Loc}(\omega)}(j^*i_*\sigma^I, \omega[k]) = 0$ since $j^*i_* = 0$. Finally, the map $\alpha$ is right universal by construction, and thus Proposition 2.5 gives the claim.

Next, we prove that $\tilde{\sigma} \oplus i_*\sigma$ is partial silting, which amounts to prove

(S2) $(\tilde{\sigma} \oplus i_*\sigma)^{\perp > 0}$ is closed under coproducts.

Consider a family of objects $X_j$, $j$ in some index set $J$, with $X_j \in (\tilde{\sigma} \oplus i_*\sigma)^{\perp > 0}$ for each $j$. Then, from the long exact sequence

$$\mathrm{Hom}_{\mathcal{D}}(i_*\sigma^I[1], X_j[k]) \xrightarrow{(\alpha, X_j[k])} \mathrm{Hom}_{\mathcal{D}}(j_*\omega, X_j[k]) \to \mathrm{Hom}_{\mathcal{D}}(\tilde{\sigma}, X_j[k]) \to \mathrm{Hom}_{\mathcal{D}}(i_*\sigma^I, X_j[k])$$

we infer that for all $j \in J$ and $k \geq 2$, $\mathrm{Hom}_{\mathcal{D}}(j_*\omega, X_j[k]) = 0$ and the map $\mathrm{Hom}_{\mathcal{D}}(\alpha, X_j[1])$ is surjective.

By hypotheses (1) and (3), $i_*\sigma$ and $j_*\omega$ are partial silting, and thus the two outer terms in the following exact sequence are zero for $k \geq 1$

$$\mathrm{Hom}_{\mathcal{D}}(j_*\omega, \bigoplus_{j \in J} X_j[k]) \to \mathrm{Hom}_{\mathcal{D}}(\tilde{\sigma}, \bigoplus_{j \in J} X_j[k]) \to \mathrm{Hom}_{\mathcal{D}}(i_*\sigma^I, \bigoplus_{j \in J} X_j[k]).$$

Moreover, by hypothesis (3) and Lemma 2.15, the map $\mathrm{Hom}_{\mathcal{D}}(\alpha, \bigoplus X_j)$ is surjective. It follows that $\bigoplus_{j \in J} X_j$ belongs to $\tilde{\sigma}^{\perp > 0}$ and thus to $(\tilde{\sigma} \oplus i_*\sigma)^{\perp > 0}$.

Finally, we prove the generating property (S3). Notice that $\tilde{\sigma} \oplus i_*\sigma$ is a generator if and only if so is $i_*\sigma \oplus j_*\omega$. Assume then $\mathrm{Hom}_{\mathcal{D}}(i_*\sigma \oplus j_*\omega, X[k]) = 0$ for all $k \in \mathbb{Z}$. In particular, $0 = \mathrm{Hom}_{\mathcal{D}}(i_*\sigma^I, X[k]) \cong \mathrm{Hom}_{\omega^{\perp_\mathbb{Z}}}(\sigma^I, i^!X[k])$ for all $k \in \mathbb{Z}$, but since $\sigma$ is silting in $\omega^{\perp_\mathbb{Z}}$, it follows $i^!X = 0$. Consider then the triangle

$$i_*i^!X \to X \to j_*j^*X \to i_*i^!X[1]$$

induced by the recollement. Since $i^!X = 0$, the object $X$ is isomorphic to $j_*j^*X$. Thus

$$0 = \mathrm{Hom}_{\mathcal{D}}(j_*\omega, X[k]) \cong \mathrm{Hom}_{\mathcal{D}}(j_*\omega, j_*j^*X[k]) \cong \mathrm{Hom}_{\mathrm{Loc}(\omega)}(\omega, j^*X[k])$$

for any $k \in \mathbb{Z}$ and since $\omega$ is silting in $\mathrm{Loc}(\omega)$, it follows $j^*X = 0$. Thus $X = 0$. □

In a similar fashion to Theorem 2.14, we present a connection of this version of the algorithm to t-structures.

Consider the t-structures $(\mathcal{Y}', \mathcal{Y}'') = (\sigma^{\perp > 0}, \sigma^{\perp \leq 0})$ and $(\mathcal{X}', \mathcal{X}'') = (\omega^{\perp > 0}, \omega^{\perp \leq 0})$ respectively in $\omega^{\perp_\mathbb{Z}}$ and in $\mathrm{Loc}(\omega)$, associated to the silting objects $\sigma$ and $\omega$ respectively. Denote by $\rho = \tilde{\sigma} \oplus i_*\sigma$ the silting object constructed as in Theorem 2.16. By hypothesis (1) and Proposition 2.9, the recollement has a lower adjacent recollement

$$\mathcal{X} \xrightarrow[\substack{\leftarrow j^* \\ \longrightarrow j_* \\ \leftarrow j^\# }]{} \mathcal{D} \xrightarrow[\substack{\leftarrow i_* \\ \longrightarrow i^! \\ \leftarrow i_\# }]{} \mathcal{Y}.$$

Consider the following t-structures in $\mathcal{D}$:

- $(\rho^{\perp_{>0}}, \rho^{\perp_{\leq 0}})$ associated to the silting object $\rho$ and

- $(\mathcal{D}', \mathcal{D}'')$ obtained by gluing the two t-structures $(\mathcal{Y}', \mathcal{Y}'')$ and $(\mathcal{X}', \mathcal{X}'')$ along the lower adjacent recollement, as in Proposition 1.3.

**Theorem 2.17.** *In the situation above, the t-structures $(\rho^{\perp_{>0}}, \rho^{\perp_{\leq 0}})$ and $(\mathcal{D}', \mathcal{D}'')$ coincide.*

*Proof.* Consider the coaisle $\mathcal{D}''$ of the glued t-structure:
$$\mathcal{D}'' = \left\{ Z \in \mathcal{D} \mid j^{\#}Z \in \omega^{\perp_{\leq 0}}, i^{!}Z \in \sigma^{\perp_{\leq 0}} \right\}.$$

The condition $j^{\#}Z \in \omega^{\perp_{\leq 0}}$ is equivalent to $\mathrm{Hom}_{\mathcal{D}}(j_*\omega, Z[k]) = 0$ for $k \leq 0$, and the condition $i^{!}Z \in \sigma^{\perp_{\leq 0}}$ to $\mathrm{Hom}_{\mathcal{D}}(i_*\sigma, Z[k]) = 0$ for $k \leq 0$. From the long exact sequence

$$\mathrm{Hom}_{\mathcal{D}}(j_*\omega, Z[k]) \to \mathrm{Hom}_{\mathcal{D}}(\tilde{\sigma}, Z[k]) \to \mathrm{Hom}_{\mathcal{D}}(i_*\sigma^I, Z[k]) \to \mathrm{Hom}_{\mathcal{D}}(j_*\omega, Z[k+1])$$

obtained applying the functor $\mathrm{Hom}_{\mathcal{D}}(-, Z[k])$ to the triangle (2.4'), it is clear that the condition $\mathrm{Hom}_{\mathcal{D}}(i_*\sigma, Z[k]) = 0 = \mathrm{Hom}_{\mathcal{D}}(j_*\omega, Z[k])$ for all $k \leq 0$ is equivalent to $\mathrm{Hom}_{\mathcal{D}}(\tilde{\sigma}, Z[k]) = 0 = \mathrm{Hom}_{\mathcal{D}}(i_*\sigma, Z[k])$ for all $k \leq 0$ (consider that $I$ is a finite set). Thus $\mathcal{D}''$ coincides with

$$\rho^{\perp_{\leq 0}} = \{Z \in \mathcal{D} \mid \mathrm{Hom}_{\mathcal{D}}(i_*\sigma, Z[k]) = 0, \mathrm{Hom}_{\mathcal{D}}(\tilde{\sigma}, Z[k]) = 0 \text{ for } k \leq 0\}.$$

Now the aisles of the two co-t-structures coincide since they are left Hom-perpendicular to the coaisles, and these coincide. $\square$

Also Theorem 2.16 restricts to tilting objects.

**Corollary 2.18.** *In the hypotheses of Theorem 2.16, assume $\sigma$ and $\omega$ are tilting objects respectively in $\omega^{\perp_{\mathbb{Z}}}$ and in $\mathrm{Loc}(\omega)$. Assume moreover the following additional condition is satisfied:*

(4) $\mathrm{Hom}_{\mathcal{D}}(j_*\omega, i_*\sigma^{(J)}[k]) = 0$ *for $k < 0$, for any set $J$.*

*Then the object $\tilde{\sigma} \oplus i_*\sigma$ is tilting in $\mathcal{D}$. In particular, if $\mathcal{D} = \mathbf{D}(R)$ and $\sigma$ and $\omega$ are isomorphic to $R$-modules of projective dimension at most 1, then $\tilde{\sigma} \oplus i_*\sigma$ is a tilting $R$-module.*

*Proof.* The proof is dual to the one of Corollary 2.12. We have to check that the following conditions hold for any set $J$ and for $k < 0$:

(i) $\mathrm{Hom}_{\mathcal{D}}(i_*\sigma, i_*\sigma^{(J)}[k]) = 0$,

(ii) $\mathrm{Hom}_{\mathcal{D}}(\tilde{\sigma}, i_*\sigma^{(J)}[k]) = 0$,

(iii) $\mathrm{Hom}_{\mathcal{D}}(i_*\sigma, \tilde{\sigma}^{(J)}[k]) = 0$,

(iv) $\mathrm{Hom}_{\mathcal{D}}(\tilde{\sigma}, \tilde{\sigma}^{(J)}[k]) = 0$.

Condition (i) holds because $i_*$ is fully faithful. For (iii), apply the functor $\operatorname{Hom}_{\mathcal{D}}(i_*\sigma, -)$, to the triangle

$$(i_*\sigma^I)^{(J)} \to \tilde{\sigma}^{(J)} \to j_*\omega^{(J)} \to (i_*\sigma^I)^{(J)}[1].$$

An application of the functor $\operatorname{Hom}_{\mathbf{D}(R)}(j_*\omega, -)$ to the same triangle yields also $\operatorname{Hom}_{\mathcal{D}}(j_*\omega, \tilde{\sigma}^{(J)}[k]) = 0$, for $k < 0$.

For conditions (ii) and (iv) apply respectively the functors $\operatorname{Hom}_{\mathcal{D}}(-, i_*\sigma^{(J)}[k])$ and $\operatorname{Hom}_{\mathcal{D}}(-, \tilde{\sigma}^{(J)}[k])$ to the triangle

$$i_*\sigma^I \to \tilde{\sigma} \to j_*\omega \to i_*\sigma^I[1]. \qquad \square$$

In case $\mathcal{D} = \mathbf{D}(R)$ and $\sigma$ and $\omega$ have cohomologies concentrated in degree 0, then the triangle (2.4') is induced by a short exact sequence

$$0 \to \sigma^I \to \tilde{\sigma} \to \omega \to 0$$

in $\operatorname{Mod} R$. In particular, $\tilde{\sigma}$ is an $R$-module and the short exact sequence implies $\tilde{\sigma}$, and hence $\tilde{\sigma} \oplus i_*\sigma$, has projective dimension at most 1. By Lemma 1.21 we conclude that $\tilde{\sigma} \oplus i_*\sigma$ is a tilting module.

## 2.7 Recollements induced by ring epimorphisms

In this section we specify some results of Section 2.4 to the particular case of 2-term silting complexes along a recollement induced by a homological ring epimorphism.

Consider a homological ring epimorphism $\lambda \colon R \to S$ and assume it is a *silting* ring epimorphism. Let $\omega$ be a 2-term partial silting complex such that the two subcategories $\mathcal{Y}_\omega$ and $\mathcal{X}_S$ of $\operatorname{Mod} R$ coincide (see Section 1.7). Notice that $\mathbf{D}(S)$ can be described as the full subcategory of $\mathbf{D}(R)$ of objects with cohomologies lying in $\operatorname{Mod} S$ ([6, Lemma 4.6]). On the other hand, by [13, Proposition 4.2], the subcategory $\operatorname{Loc}(\omega)^{\perp_0} = \omega^{\perp_\mathbb{Z}}$ consists precisely of those complexes whose cohomologies lie in $\mathcal{Y}_\omega$. Since $\mathcal{X}_S = \mathcal{Y}_\omega$, we can conclude that $\lambda_*(\mathbf{D}(S)) = i_*(\omega^{\perp_\mathbb{Z}})$ and thus the recollements (1.5) and (2.3) are equivalent and of the form

$$\mathbf{D}(S) \xrightarrow{\lambda_*} \mathbf{D}(R) \xrightarrow{j^*} \operatorname{Loc}(\omega). \qquad (2.5)$$

with $\lambda^* = -\otimes^{\mathbf{L}}_R S$, $\lambda^! = \mathbf{R}\operatorname{Hom}_R(S,-)$, and $j_!$, $j_*$.

In particular, it was shown in [48, Theorem 6.7] that universal localizations in the sense of [66] are silting ring epimorphisms. In the same paper [48], moreover, an explicit construction of the silting object associated to a universal localization is shown. The situation is depicted in Figure 1.

We can then state the following consequence of Theorem 2.10.

**Corollary 2.19.** *Let $\lambda \colon R \to S$ be a silting ring epimorphism which is homological, arising from a 2-term partial silting complex $\omega$, and let $\sigma$ be a 2-term silting complex in $\mathbf{D}(S)$. Assume the following additional hypotheses are satisfied:*

**Figure 1:** A diagram summarizing the situation. The triangles commute.

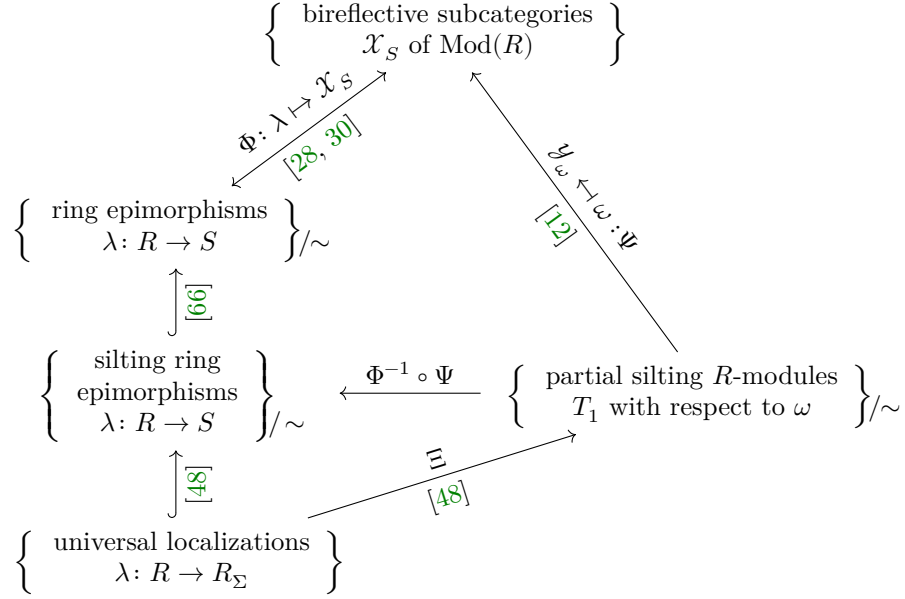

(1) *the class* $(\lambda_*\sigma)^{\perp_{>0}} \cap (j_!\omega)$ *is closed under coproducts in* $\mathbf{D}(R)$,

(2) $\lambda_*\sigma$ *is a 2-term complex of projectives in* $\mathbf{D}(R)$.

*Consider the set* $I = \mathrm{Hom}_{\mathbf{D}(R)}(\lambda_*\sigma, j_!\omega[1])$, *let* $\alpha\colon \lambda_*\sigma^{(I)} \to j_!\omega[1]$ *be the canonically induced map and* $\tilde{\sigma}$ *be the cocone of* $\alpha$:

$$j_!\omega \to \tilde{\sigma} \to \lambda_*\sigma^{(I)} \xrightarrow{\alpha} j_!\omega[1].$$

*Then* $\tilde{\sigma} \oplus \lambda_*\sigma$ *is a 2-term silting complex in* $\mathbf{D}(R)$.

*Proof.* Since $j_!\omega$ and $\lambda_*\sigma$ are 2-term complexes of projectives, hypothesis (2) of Theorem 2.10 is satisfied; thus Theorem 2.10 applies and $\tilde{\sigma} \oplus \lambda_*\sigma$ is a silting complex in $\mathbf{D}(R)$. Finally, the triangle defining $\tilde{\sigma}$ implies that it is a 2-term complex of projectives. □

*Remark* 2.20. Let $\lambda\colon R \to S$ be a homological ring epimorphism such that the projective dimension of $S$ as an $R$-module is not greater than 1. Then

(1) $\omega = \mathrm{Cone}(\lambda)$ is a 2-term partial silting complex and $\lambda$ is the corresponding silting ring epimorphism.

(2) If $\lambda$ is injective, then $\omega$ is quasi-isomorphic to the partial tilting module $T_1 = \mathrm{Coker}\,\lambda$.

(3) Assume $\lambda$ is injective and $R$ hereditary. Then for every tilting $S$-module $T$ satisfying

    (1) $(\lambda_*T)^{\perp_1}$ is closed under coproducts,

the object $\tilde{T} \oplus \lambda_* T$ is a tilting module, where $\tilde{T}$ is constructed as the cocone in $\mathbf{D}(R)$ of the left universal map $\alpha\colon \lambda_* T^{(I)} \to T_1[1]$, for a set of generators $I$ of $\operatorname{Hom}_{\mathbf{D}(R)}(\lambda_* T, T_1[1])$.

*Proof.* (1) Considering a projective resolution $0 \to P_{-1} \xrightarrow{\sigma} P_0 \to S \to 0$ of $S$ as $R$-module, we have the following isomorphisms of triangles in $\mathbf{D}(R)$

$$\begin{array}{ccccccc} R & \longrightarrow & \sigma & \longrightarrow & \omega & \longrightarrow & R[1] \\ \| & & \downarrow \wr & & \downarrow \wr & & \| \\ R & \xrightarrow{\lambda} & S & \longrightarrow & \operatorname{Cone}(\lambda) & \longrightarrow & R[1] \end{array}$$

and we know $H^0(\omega) = \operatorname{Coker}(\lambda)$ is a partial silting module with respect to $\omega$ and with silting ring epimorphism $\lambda$ ([10], see also [48, Example 6.5]).

(2) If $\lambda$ is injective, $H^{-1}(\omega) = \operatorname{Ker}(\lambda) = 0$, hence $T_1 = \operatorname{Coker}(\lambda)$ is partial silting with respect to the injective map $\omega$, and thus partial tilting.

(3) It follows from Corollaries 2.12 and 2.19. □

## 2.8 Recollements induced by triangular matrix rings

We recall here some well-known facts about modules over a ring. Let $R$ be a ring. For idempotent elements $e_1$ and $e_2$, $\operatorname{Hom}_R(e_1 R, e_2 R) \cong e_2 R e_1$. Define the trace $\tau_M(N)$ of a right $R$-module $M$ in $N$ as the sum of the images of all morphisms in $\operatorname{Hom}_R(M, N)$. Note that if $(N_j)_{j \in J}$ is an indexed set of $R$-modules, then $\tau_M(\bigoplus_{j \in J} N_j) \cong \bigoplus_{j \in J} \tau_M(N_j)$ [3, Proposition 8.18]. For $e$ an idempotent element and $P = eR$ the corresponding projective right $R$-module, the two-sided ideal $ReR$ can be computed as the trace $\tau_P(R)$ [42, Remark 2.41].

**Lemma 2.21.** *The canonical epimorphism* $\lambda\colon R \to R/ReR$ *is the universal localization of $R$ at the projective module* $P = eR$, *that means at the zero morphism* $\Sigma = \{0 \to P\}$.

*Proof.* The map $\lambda$ is $\Sigma$-inverting, since $eR \otimes_R (R/ReR) = 0$. It is universal because if $\mu\colon R \to S$ is $\Sigma$-inverting, then $eR \otimes_R S = 0$. Thus $ReR \otimes_R S = 0$, and $ReR \subseteq \operatorname{Ker} \mu$ and $\mu$ factors through $\lambda$. See also [6, Example 4.5(2)]. □

We will use the following result for a ring $R$.

**Theorem 2.22** (Keller). *[34, Theorem 8.5] [6] Let $\mathcal{D}$ be a full subcategory of $\mathbf{D}(R)$ closed under coproducts. If $T$ is a compact generator of $\mathcal{D}$, then there is a differential graded algebra $E = \mathbf{R}\operatorname{Hom}(T, T)$ with homology $H^*(E) \cong \bigoplus_{i \in \mathbb{Z}} \operatorname{Hom}_{\mathcal{D}}(T, T[i])$ such that the functor $- \otimes_E^{\mathbf{L}} T\colon \mathbf{D}(E) \to \mathcal{D}$ is a triangle equivalence. In particular, if $T$ is exceptional, then $E = \operatorname{End}_{\mathcal{D}}(T)$.*

Let now $A$ and $B$ be arbitrary rings and $_A M_B$ be an $A$-$B$-bimodule. Define the matrix ring

$$R = \begin{bmatrix} A & _A M_B \\ 0 & B \end{bmatrix}.$$

We aim at investigating tilting $R$-modules comparing them to tilting $A$- and $B$-modules.

Let $e_A = \begin{bmatrix} 1_A & 0 \\ 0 & 0 \end{bmatrix}$ and $e_B = \begin{bmatrix} 0 & 0 \\ 0 & 1_B \end{bmatrix} = 1_R - e_A$ be the idempotents corresponding respectively to $A$ and $B$.

**Lemma 2.23.** *For a matrix ring $R$ as above, the quotient ring $R/Re_A R$ is isomorphic to $B$.*

*Proof.* We compute the two-sided ideal $Re_A R$, as a right $R$-module:
$$\begin{aligned} Re_A R &\cong \tau_{e_A R}(R) \\ &\cong \tau_{e_A R}(e_A R) \oplus \tau_{e_A R}(e_B R) \\ &= e_A R. \end{aligned}$$

For the last equality, notice that $\tau_{e_A R}(e_A R) = e_A R$, while $\tau_{e_A R}(e_B R) = 0$ since $\operatorname{Hom}_R(e_A R, e_B R) = e_B R e_A = 0$. Thus, the quotient is
$$R/Re_A R \cong (e_A R \oplus e_B R)/e_A R \cong e_B R.$$

Now, as rings, $e_B R \cong e_B R e_B \cong B$. □

**Theorem 2.24.** *Let $A$ and $B$ be algebras over a field $k$ and $R$ be a matrix ring as above. Let $T$ be a tilting $B$-module and $P = e_A R$. Let $\tilde{T}$ be the Bongartz-Trlifaj extension of $T$ and $P$ constructed as in [68, Lemma 6.8], so that*
$$0 \to P \to \tilde{T} \to T^{(I)} \to 0$$
*is a short exact sequence in $R$, with $I = \operatorname{Ext}^1_R(T, P)$. Then $\tilde{T} \oplus T_R$ is a tilting $R$-module.*

*Proof.* We want to apply Corollary 2.19. Let then $\mathcal{U} = \{P\}$; we have already seen that $R_{\mathcal{U}} = R/Re_A R \cong B$ and the canonical epimorphism $\lambda\colon R \to R/Re_A R$ is the universal localization; as a silting ring epimorphism it arises from the partial silting module $P$. Note that $Re_A R$ is an idempotent ideal; moreover, as right $R$-modules, $Re_A R \cong e_A R$ (see the proof of Lemma 2.23), and thus it is projective. For $J = Re_A R$, the exact sequence
$$\operatorname{Tor}^R_{k+1}(R, R/J) \to \operatorname{Tor}^R_{k+1}(R/J, R/J) \to \operatorname{Tor}^R_k(J, R/J).$$
shows that $\operatorname{Tor}^R_k(R/J, R/J) = 0$ for $k \geq 2$, since $J$ and $R$ are projective. Using the well-known identity $\operatorname{Tor}^R_1(R/I, R/I') \cong (I \cap I')/II'$, for any two ideals $I$ and $I'$, we get that also $\operatorname{Tor}^R_1(R/J, R/J) \cong J/J^2 = 0$, since $J$ is idempotent. By [30, Theorem 4.4], the ring epimorphism $\lambda\colon R \to B$ is homological.

Again the proof of Lemma 2.23 shows that $B$ is projective as right $R$-module. Thus, projective $B$-modules are projective as $R$-modules and the projective dimension is preserved, implying hypothesis (2).

Finally, to prove hypothesis (1), use the criterion in Remark 2.11 and show that the recollement induced by $\lambda$ can be extended one step downwards. This is also shown in [7, Example 3.4] and we recall the construction. The recollement has the form

$$\mathbf{D}(B) \xrightarrow{\lambda_*} \mathbf{D}(R) \xrightarrow{j^*} \operatorname{Loc}(P) \cong \mathbf{D}(A) ,$$

with $\lambda^* = -\otimes^{\mathbf{L}}_R B$, $\lambda^! = \mathbf{R}\operatorname{Hom}_R(B,-)$, and $j_!$, $j_*$.

where $A = \mathrm{End}_R(P)$ and the equivalence $\mathrm{Loc}(P) \cong \mathbf{D}(A)$ is provided by Keller's Theorem 2.22. Considering the ideals of $R$ as left $R$-modules, we can prove the dual version of Lemma 2.23, namely, that $R/Re_B R \cong A$. The epimorphism $\mu \colon R \to R/Re_B R \cong A$ is a universal localization and the induced functor $\mu_* \colon \mathbf{D}(A) \to \mathbf{D}(R)$ can be expressed as $\mu_* = \mathbf{R}\mathrm{Hom}_A({}_R A_A, -) \cong - \otimes^{\mathbf{L}}_A {}_A A_R$. As a Hom functor, it coincides with the functor $j_*$ of the recollement (see [7, §2.2.3]); as a tensor functor it has a right adjoint $\mu^! = j^\# = \mathrm{Hom}_R({}_A A_R, -) \colon \mathbf{D}(R) \to \mathbf{D}(A)$. By [7, Proposition 3.2], the recollement can be extended one step downwards. This proves (1) and concludes the proof. □

## 2.9 Silting objects obtained by gluing

We investigate here the conditions under which a silting object in the central category of a recollement is glued from two silting objects. Consider a recollement

$$\mathcal{Y} \xrightleftharpoons[i^!]{\overset{i^*}{\longleftarrow}}_{i_*} \mathcal{D} \xrightleftharpoons[j_*]{\overset{j_!}{\longleftarrow}}_{j^*} \mathcal{X}$$

of well generated triangulated categories.

**Definition 2.25.** Let $\omega$ and $\sigma$ be objects in $\mathcal{X}$ and $\mathcal{Y}$, respectively. We say that a silting object $\tau$ in $\mathcal{D}$ is *glued from $\sigma$ and $\omega$* if there exists a triangle

$$j_!\omega \to \tilde{\sigma} \to i_*\sigma^{(I)} \xrightarrow{\alpha} j_!\omega[1] \qquad (\#)$$

such that $\alpha$ is an $\mathrm{Add}(i_*\sigma)$-precover and $\tilde{\sigma} \oplus i_*\sigma$ is a silting object equivalent to $\tau$.

**Theorem 2.26.** *In the situation above, the silting object $\tau$ is glued from two silting objects $\sigma$ in $\mathcal{Y}$ and $\omega$ in $\mathcal{X}$ if and only if the following conditions hold:*

*(i) $i_* i^* \tau \in \mathrm{Add}(\tau)$,*

*(ii) $j^* \tau$ is partial silting in $\mathcal{X}$.*

*Proof.* Let $\sigma$ and $\omega$ be silting objects respectively in $\mathcal{Y}$ and in $\mathcal{X}$, let $\tilde{\sigma}$ be constructed as in triangle (#), and let $\tau = \tilde{\sigma} \oplus i_*\sigma$. An application of the functor $i^*$ to the triangle (#) yields $i^*\tilde{\sigma} \cong i^* i_* \sigma^{(I)}$, hence $i_* i^* \tau \cong i_* i^* \tilde{\sigma} \oplus i_* i^* i_* \sigma \cong i_* \sigma^{(J)}$ belongs to $\mathrm{Add}\,\tau$ and (i) holds. An application of the functor $j^*$ to the same triangle gives $j^* j_! \omega \cong j^* \tilde{\sigma}$, hence $j^* \tau = j^* \tilde{\sigma} \oplus j^* i_* \sigma \cong j^* j_! \omega \cong \omega$, that is partial silting in $\mathcal{X}$ by hypothesis, proving (ii).

Conversely, suppose $\tau$ is a silting object in $\mathcal{D}$ satisfying (i) and (ii). We claim that $\omega = j^*\tau$ is silting in $\mathcal{X}$. By hypothesis (ii), we only have to check (S3). Suppose then an object $X$ in $\mathcal{X}$ satisfies $\mathrm{Hom}_\mathcal{D}(\omega, X[i]) = 0$ for all $i \in \mathbb{Z}$. Then $\mathrm{Hom}_\mathcal{D}(\tau, j_*X[i]) = 0$ for all $i \in \mathbb{Z}$, hence $j_*X = 0$ since $\tau$ is silting in $\mathcal{D}$, that implies $X = 0$ since $j_*$ is fully faithful, proving the claim. Since $\mathcal{X}$ is well generated, we conclude that $\mathcal{X} = \mathrm{Loc}(\omega)$ and $\mathcal{Y} = \omega^{\perp_\mathbb{Z}}$.

Next, we claim the object $\sigma = i^*\tau$ is silting in $\omega^{\perp_\mathbb{Z}}$. Property (S1) follows since

$$\begin{aligned}\mathrm{Hom}_{\omega^{\perp_\mathbb{Z}}}(\sigma, \sigma[i]) &\cong \mathrm{Hom}_{\omega^{\perp_\mathbb{Z}}}(i^*\tau, i^*\tau[i]) \\ &\cong \mathrm{Hom}_\mathcal{D}(\tau, i_* i^*\tau[i]) = 0\end{aligned}$$

for $i > 0$ by hypothesis (i) and because $\tau$ is silting in $\mathcal{D}$.

For (S2), let $(Y_j)_{j \in J}$ be a family of objects in $\omega^{\perp_\mathbb{Z}}$ lying in $\sigma^{\perp_{>0}}$, and let $Y = \coprod_{j \in J} Y_j$. Then, for $i > 0$,

$$\mathrm{Hom}_{\omega^{\perp_\mathbb{Z}}}(\sigma, Y[i]) \cong \mathrm{Hom}_\mathcal{D}(\tau, i_* Y[i])$$
$$\cong \mathrm{Hom}_\mathcal{D}(\tau, \coprod_{j \in J} i_* Y_j[i]) = 0$$

because $\mathrm{Hom}_\mathcal{D}(\tau, i_* Y_j[i]) \cong \mathrm{Hom}_{\omega^{\perp_\mathbb{Z}}}(\sigma, Y_j[i]) = 0$ for all $i > 0$ and $\tau^{\perp_{>0}}$ is closed under coproducts in $\mathcal{D}$.

To prove (S3), suppose $Y \in \omega^{\perp_\mathbb{Z}}$ satisfies $\mathrm{Hom}_{\omega^{\perp_\mathbb{Z}}}(\sigma, Y[i]) = 0$ for all $i \in \mathbb{Z}$. Then $\mathrm{Hom}_\mathcal{D}(\tau, i_* Y[i]) = 0$ for all $i \in \mathbb{Z}$ and since $\tau$ is a generator, $i_* Y = 0$ and thus $Y = 0$ since $i_*$ is fully faithful.

Consider then the canonical triangle

$$j_! j^* \tau \to \tau \to i_* i^* \tau \xrightarrow{\alpha} j_! \omega[1]$$

induced by the recollement, that can be rewritten as

$$j_! \omega \to \tau \to i_* \sigma \xrightarrow{\alpha} j_! \omega[1].$$

An application of the functor $\mathrm{Hom}_\mathcal{D}(i_* \sigma, -)$ to this triangle yields

$$\mathrm{Hom}_\mathcal{D}(i_* \sigma, i_* \sigma) \xrightarrow{(i_* \sigma, \alpha)} \mathrm{Hom}_\mathcal{D}(i_* \sigma, j_! \omega[1]) \to \mathrm{Hom}_\mathcal{D}(i_* \sigma, \tau[1])$$

where the last term is zero by hypothesis (i) and since $\tau$ is silting. Thus, $\alpha$ is an $\mathrm{Add}(i_* \sigma)$-precover and $\tau$ fits into a triangle of the form (#); moreover $\tau \oplus i_* \sigma$ is silting since $\mathrm{Add}(\tau \oplus i_* \sigma) = \mathrm{Add}(\tau)$ by hypothesis (i). Thus, $\tau \oplus i_* \sigma$ is glued from $\sigma$ and $\omega$ since $\tau$ is equivalent to $\tau \oplus i_* \sigma$. □

*Remark* 2.27. Assume the recollement is of the form (1.5), induced by a homological ring epimorphism $\lambda \colon R \to S$, and let $X$ be the cone of $\lambda$ as a morphism in $\mathbf{D}(R)$. Let $\tau$ be a 2-term silting complex in $\mathbf{D}(R)$ and assume there exist 2-term silting complexes $\sigma$ in $\mathbf{D}(S)$ and $\omega$ in $\mathrm{Loc}(X)$ such that $\tau$ is glued from $\sigma$ and $\omega$. Let $H^0(\tau) = T$. Then the proof of Theorem 2.26 shows $\mathrm{Add}(\lambda^* \tau) = \mathrm{Add}(\sigma)$, in particular $\lambda^* \tau$ is silting in $\mathbf{D}(S)$. Thus, by [20, Theorem 2.2], the condition

(i') $\lambda_* \lambda^* T \in \mathrm{Gen}\, T$

holds (and it is actually equivalent to $\lambda^* \tau$ being silting in $\mathbf{D}(S)$). The conditions (i) and (ii) of Theorem 2.26 are indeed stronger, since the request of $\tau$ being glued from $\sigma$ and $\omega$ implies that $\lambda^* \tau$ is silting.

## 2.10 Other gluing methods

Aim of this section is to compare the previous results, in particular Theorem 2.10 and Theorem 2.16, with other gluing methods appeared in the literature.

**A first gluing procedure**

The method described below appeared in [46] and shows how to glue (classical) silting objects in a way that is compatible with the co-t-structures associated to them. We say an object $M$ in a triangulated category $\mathcal{D}$ is *classical silting* if it satisfies the following two conditions:

- $\text{Hom}_{\mathcal{D}}(M, M[k]) = 0$ for $k > 0$,

- $\mathcal{D} = \text{thick}(M)$.

Note that, by a result of Neeman ([51], see [2, Proposition 4.2]), if $\mathcal{D}$ is a triangulated category with coproducts and $M$ is a silting object in $\mathcal{D}$ which belongs to $\mathcal{D}^c$, then $M$ is classical silting in $\mathcal{D}^c$. Conversely, if $\mathcal{D}$ is compactly generated, then it follows from [56, Example 4.2(1)] that any classical silting object in $\mathcal{D}^c$ is silting in $\mathcal{D}$. Indeed, since $\text{thick}(M) = \mathcal{D}^c$ and $\mathcal{D}$ is compactly generated, there exists a set $\mathcal{S}$ of generators of $\mathcal{D}$ which is contained in $\text{thick}(M)$.

Let $\mathcal{D}$ be a triangulated category. For a co-t-structure $(\mathcal{D}', \mathcal{D}'')$ in $\mathcal{D}$ and an object $Z$ in $\mathcal{D}$, denote by $\beta' Z$ (resp. $\beta'' Z$) a (non-functorial) choice of an object in $\mathcal{D}'$ (resp. in $\mathcal{D}''$) such that

$$\beta' Z \to Z \to \beta'' Z \to \beta' Z[1]$$

is a triangle.

**Theorem 2.28.** *[46, Theorem 3.1] Let*

$$\mathcal{Y} \underset{\substack{\longleftarrow i^* \\ \longrightarrow i_* \\ \longleftarrow i^!}}{} \mathcal{D} \underset{\substack{\longleftarrow j_! \\ \longrightarrow j^* \\ \longleftarrow j_*}}{} \mathcal{X}$$

*be a recollement of triangulated categories. Let $X$ and $Y$ be classical silting objects in $\mathcal{X}$ and $\mathcal{Y}$, respectively. Let $(\mathcal{X}', \mathcal{X}'') = ({}^{\perp_0}(X^{\perp_{>0}}), X^{\perp_{>0}})$ and $(\mathcal{Y}', \mathcal{Y}'') = ({}^{\perp_0}(Y^{\perp_{>0}}), Y^{\perp_{>0}})$ be the co-t-structures associated to $X$ and $Y$ in $\mathcal{X}$ and $\mathcal{Y}$, respectively. Let $\beta'$ and $\beta''$ be a (non-functorial) choice of truncation for the co-t-structure $(\mathcal{Y}', \mathcal{Y}'')$ in $\mathcal{Y}$, and let $\tilde{X}$ be defined by the following triangle*

$$i_* \beta' i^! j_! X \to j_! X \to \tilde{X} \to (i_* \beta' i^! j_!) X[1].$$

*Then the object $Z = \tilde{X} \oplus i_* X$ is silting in $\mathcal{D}$ and the co-t-structure associated to $Z$ is the glued co-t-structure of $(\mathcal{X}', \mathcal{X}'')$ and $(\mathcal{Y}', \mathcal{Y}'')$.*

Assume a recollement of compactly generated triangulated categories of the form (1.3) is given, and assume it restricts to the respective subcategories of compact objects. If $\omega$ and $\sigma$ are compact silting objects in $\mathcal{X}$ and $\mathcal{Y}$ respectively, they glue to a silting object $\tau$ in $\mathcal{D}$ through Theorem 2.10. On the other hand, since they are classical silting in the respective subcategories of compact objects, they glue to a classical silting object in $\mathcal{D}^c$, which is then silting in $\mathcal{D}$. The co-t-structures associated to $\tau$ and $\tilde{\tau}$ coincide, and hence the two silting objects are equivalent.

**A second gluing procedure**

Another gluing procedure was described by Saorín and Zvonareva [65]. In the next definition we introduce what is called 'partial silting' object in [65]: here we are calling it 'weak partial silting' in order to distinguish it from our notion of partial silting (Definition 1.11).

**Definition 2.29.** An object $T$ in a triangulated category $\mathcal{D}$ with coproducts is *weak partial silting* if

(1) $T$ and its positive shifts generate a t-structure, i.e. $(\mathcal{U}_T = {}^{\perp_0}(T^{\perp_{\leq 0}}), T^{\perp_{\leq 0}})$ is a torsion pair;

(2) $\mathcal{U}_T$ is contained in $T^{\perp_{>0}}$.

An object $T$ is partial silting if and only if it is weak partial silting and $T^{\perp_{>0}}$ is cocomplete (see [13, §3.1]).

Note that, in the following, our notation deviates from [65].

**Theorem 2.30.** *[65, Theorem 6.3] Let*

$$\mathcal{X} \xleftarrow{\underset{j_*}{\overset{j^*}{\longleftarrow}}}_{\underset{j^\#}{\longleftarrow}} \mathcal{D} \xleftarrow{\underset{i^!}{\overset{i_*}{\longleftarrow}}}_{\underset{i_\#}{\longleftarrow}} \mathcal{Y}$$

*be a recollement of triangulated categories, let $\omega$ and $\sigma$ be (weak partial) silting objects respectively in $\mathcal{X}$ and $\mathcal{Y}$, let $(\mathcal{X}', \mathcal{X}'')$ and $(\mathcal{Y}', \mathcal{Y}'')$ be the associated t-structures and let $(\mathcal{D}', \mathcal{D}'')$ be the glued t-structure. Suppose the following condition holds:*

($\star$) *there exists a triangle $\tilde{\sigma} \to j_*\omega \to U \to \tilde{\sigma}[1]$ such that $U \in i_*\mathcal{Y}'$ and $\tilde{\sigma} \in {}^{\perp_0}(i_*\mathcal{Y}')$.*

*Then $\tilde{\sigma} \oplus i_*\sigma$ is a (weak partial) silting object with associated t-structure $(\mathcal{D}', \mathcal{D}'')$.*

A criterion for condition ($\star$) to hold is also provided.

**Lemma 2.31.** *[65, Corollary 6.5] Assume $\mathcal{X}$ has coproducts and the recollement of Theorem 2.30 can be extended upwards:*

$$\mathcal{Y} \xleftarrow{\underset{i_*}{\overset{i^*}{\longleftarrow}}}_{\underset{i_\#}{\longleftarrow}}^{\underset{i^!}{\longleftarrow}} \mathcal{D} \xleftarrow{\underset{j^*}{\overset{j_!}{\longleftarrow}}}_{\underset{j^\#}{\longleftarrow}}^{\underset{j_*}{\longleftarrow}} \mathcal{X}.$$

*Let $\omega$ and $\sigma$ be weak partial silting objects in $\mathcal{X}$ and $\mathcal{Y}$ respectively, associated to the t-structures $(\mathcal{X}', \mathcal{X}'')$ and $(\mathcal{Y}', \mathcal{Y}'')$, respectively. If*

($\ddagger$) *the t-structure $(\mathcal{Y}', \mathcal{Y}'')$ has a left adjacent co-t-structure in $\mathcal{Y}$,*

*then condition ($\star$) of Theorem 2.30 holds, so $\tilde{\sigma} \oplus i_*\sigma$ is a weak partial silting object generating the t-structure $(\mathcal{D}', \mathcal{D}'')$ in $\mathcal{D}$ glued with respect to the lower recollement of the ladder.*

*Proof.* Glue the co-t-structures $(^{\perp_0}(\mathcal{Y}'), \mathcal{Y}')$ and $(\mathcal{X}, 0)$ with respect to the upper recollement of the ladder to get a co-t-structure whose coaisle is $i_*(\mathcal{Y}')$. Then the canonical triangle of the torsion pair $(^{\perp_0}(i_*\mathcal{Y}'), i_*\mathcal{Y}')$ is the desired triangle of condition $(\star)$. □

Notice that, by Theorem 1.13, the condition (‡) is always satisfied if $\sigma$ is a partial silting object in a well generated category $\mathcal{Y}$ satisfying dual Brown representability.

Assume a ladder of recollements as above is given, with silting objects $\omega$ in $\mathcal{X}$ and $\sigma$ in $\mathcal{Y}$ such that both Theorem 2.30 and Theorem 2.16 apply. Then the glued silting objects $\tilde{\sigma} \oplus i_*\sigma$ with respect to the procedures of the two theorems are equivalent, since they correspond to the same t-structure.

# The example of Kronecker modules 3

We now apply the construction presented in Corollary 2.19 to the instance of Kronecker modules. More precisely, consider the path algebra $R = kQ$ of the Kronecker quiver

$$Q = \; 1 \xleftarrow[\beta]{\alpha} 2$$

over an algebraically closed field $k$. It is isomorphic to the matrix ring

$$\begin{pmatrix} k & k^2 \\ 0 & k \end{pmatrix}$$

with basis $\{\begin{pmatrix} 1 & 0 \\ 0 & 0 \end{pmatrix}, \begin{pmatrix} 0 & 0 \\ 0 & 1 \end{pmatrix}, \begin{pmatrix} 0 & 1 \\ 0 & 0 \end{pmatrix}, \begin{pmatrix} 0 & x \\ 0 & 0 \end{pmatrix}\}$.

Denote by $P_i$ (respectively $Q_i$), with $i \in \mathbb{N}$, the (finite dimensional) indecomposable preprojective (respectively preinjective) module, indexed such that $\dim_k \operatorname{Hom}_R(P_i, P_{i+1}) = 2$ (resp. $\dim_k \operatorname{Hom}_R(Q_{i+1}, Q_i) = 2$).

The following list gives a complete classification of silting $R$-modules [12, Examples 5.10 and 5.18], together with the corresponding 2-term silting complex, when not quasi-isomorphic:

$$\text{silting non tilting} \begin{cases} 0, \text{ with respect to } P_1[1] \oplus P_2[1] \\ P_1, \text{ with respect to } P_2[1] \oplus P_1 \\ Q_1, \text{ with respect to } P_1[1] \oplus (P_1^2 \to P_2) \end{cases}$$

$$\text{compact tilting} \begin{cases} P_i \oplus P_{i+1}, & i \geq 2 \\ Q_{i+1} \oplus Q_i, & i \geq 1 \\ R = P_1 \oplus P_2 \end{cases}$$

$$\text{non-compact tilting} \begin{cases} R_\mathcal{U} \oplus R_\mathcal{U}/R & 0 \neq \mathcal{U} \in \mathbb{P}_k^1 \\ L \text{ (the Lukas module)} \end{cases}$$

We show that we can retrieve all silting modules, except the Lukas, by gluing silting modules along appropriate recollements of the derived category of Kronecker modules. We summarize here the outcome of the next sections.

## 3.1 Summary

To recover large (non-compact) silting $R$-modules, we consider the recollements induced by a universal localization $R \to R_S$ at a simple regular $R$-module $S$. Such a recollement has the form

$$\mathbf{D}(k[x]) \underset{\substack{\longleftarrow \lambda^* \longrightarrow \\ \longleftarrow \lambda^! \longrightarrow}}{\xrightarrow{\lambda_*}} \mathbf{D}(R) \underset{\substack{\longleftarrow j_! \longrightarrow \\ \longleftarrow j_* \longrightarrow}}{\xrightarrow{j^*}} \mathrm{Loc}(S_\infty)$$

and glued silting modules have always the Prüfer $S_\infty$ as a direct summand. Thus, when we let $S$ run through all simple regular modules, and $T$ run through all (large) silting $R_S$-modules, we recover all large tilting $R$-modules having at least one Prüfer summand, that is, all except for the Lukas. This construction is presented in Section 3.2.

To recover compact silting modules, we use instead the recollement induced by universal localizations $R \to R_{\{P\}}$ at an indecomposable preprojective or preinjective modules. This recollement has the form

$$\mathbf{D}(k) \underset{\substack{\longleftarrow \\ \longleftarrow}}{\xrightarrow{\lambda_*}} \mathbf{D}(R) \underset{\substack{\longleftarrow j_! \longrightarrow \\ \longleftarrow j_* \longrightarrow}}{\xrightarrow{j^*}} \mathbf{D}(k).$$

We compute explicitly the images of $j_!$ and $i_*$ by separating three cases:

- $P$ is projective,
- $P$ is preprojective not projective, or
- $P$ is preinjective.

The first case is studied in Section 3.4 and we further distinguish two subcases:

(1) Localization at the projective $P_1 = e_1 R$. In this case, $\mathrm{Im}\,\lambda_* = \mathrm{Loc}(Q_1)$ and $\mathrm{Im}\,j_! = \mathrm{Loc}(P_1)$; we can glue the silting module $Q_1$ in $\mathrm{Loc}(Q_1)$ and the silting module $P_1$ in $\mathrm{Loc}(P_1)$, to get the glued silting $R$-module $Q_1 \oplus Q_2$. Or, we can glue the silting module $Q_1$ in $\mathrm{Loc}(Q_1)$ and the silting module $0$ (with respect to $P_1 \to 0$) in $\mathrm{Loc}(P_1)$, to get the glued silting $R$-module $Q_1$.

(2) Localization at the projective $P_2 = e_2 R$. Now, $\mathrm{Im}\,\lambda_* = \mathrm{Loc}(P_1)$ and $\mathrm{Im}\,j_! = \mathrm{Loc}(P_2)$; we glue the silting module $P_1$ in $\mathrm{Loc}(P_1)$ and the silting module $0$ (with respect to $P_2 \to 0$) in $\mathrm{Loc}(P_2)$, to get the glued silting $R$-module $P_1$. Or, we can glue the silting module $0$ in $\mathrm{Loc}(P_1)$ (with respect to $P_1 \to 0$) and the silting module $P_2$ in $\mathrm{Loc}(P_2)$, to get the glued silting $R$-module $Q_1$. Finally, we can glue the silting modules $P_1$ in $\mathrm{Loc}(P_1)$ and $P_2$ in $\mathrm{Loc}(P_2)$ to get the silting $R$-module $P_1 \oplus P_2$. The latter case is studied in Section 3.5.

The second case, that is, if $P = P_i$, for $i \geq 3$, is preprojective not projective, is studied in Section 3.5. In this case, $\mathrm{Im}\,\lambda_* = \mathrm{Loc}(P_{i-1})$ and $\mathrm{Im}\,j_! = \mathrm{Loc}(P_i)$. We can glue the silting modules $P_{i-1}$ in $\mathrm{Loc}(P_{i-1})$ and $P_i$ in $\mathrm{Loc}(P_i)$ to get the silting $R$-module $P_{i-1} \oplus P_i$. If $i = 3$, since $P_{i-1} = P_2$ is projective, we can glue the silting module $0$ (with respect to $P_2 \to 0$) in $\mathrm{Loc}(P_2)$ and the silting module $P_3$ in $\mathrm{Loc}(P_3)$, getting the glued silting $R$-module 0 (with respect to $P_1[1] \oplus P_2[1]$).

Finally, the case $P = Q_i$, for $i \geq 1$, is preinjective, is studied in Section 3.6. In this case, $\mathrm{Im}\,\lambda_* = \mathrm{Loc}(Q_{i+1})$ and $\mathrm{Im}\,j_! = \mathrm{Loc}(Q_i)$. We glue the silting modules $Q_{i+1}$ in $\mathrm{Loc}(Q_{i+1})$ and $Q_i$ in $\mathrm{Loc}(Q_i)$ to get the silting $R$-module $Q_{i-1} \oplus Q_i$.

This exhausts the list of compact silting $R$-modules.

## 3.2 Gluing large silting objects

Consider the morphism $\bar{\alpha}\colon P_1 \to P_2$ induced by the arrow $\alpha$ and let the simple regular $S$ be its cokernel. Consider the universal localization $\lambda\colon R \to R_S$ of $R$ at the set $\{\bar{\alpha}\}$. The ring $R_S$ can be regarded as the path algebra $kQ'/\sim$ of the quiver

$$Q' = 1 \underset{\beta}{\overset{\alpha^{-1}}{\underset{\alpha}{\rightleftarrows}}} 2$$

modulo the ideal generated by the relations $\alpha\alpha^{-1} - e_1$ and $\alpha^{-1}\alpha - e_2$, and it is isomorphic to the 2-by-2 matrix ring $M_2(k[x])$ over the polynomial ring, where the isomorphism is given by the map

$$\begin{aligned}
\varphi\colon R_S &\longrightarrow M_2(k[x]) \\
e_1 &\longmapsto \begin{pmatrix} 1 & 0 \\ 0 & 0 \end{pmatrix} \\
e_2 &\longmapsto \begin{pmatrix} 0 & 0 \\ 0 & 1 \end{pmatrix} \\
\alpha &\longmapsto \begin{pmatrix} 0 & 1 \\ 0 & 0 \end{pmatrix} \\
\alpha^{-1} &\longmapsto \begin{pmatrix} 0 & 0 \\ 1 & 0 \end{pmatrix} \\
\beta &\longmapsto \begin{pmatrix} 0 & x \\ 0 & 0 \end{pmatrix}
\end{aligned}$$

extended to the paths in $kQ'$ by setting $\varphi(st) = \varphi(s)\varphi(t)$ for any two paths $s$ and $t$ and by linearity to the whole algebra $kQ'/\sim$. One can check that such map is well-defined and it is a morphism of $k$-algebras. To see that it is bijective note that for any positive integer $n$, $\varphi((\beta\alpha^{-1})^n) = \begin{pmatrix} x^n & 0 \\ 0 & 0 \end{pmatrix}$ and $\varphi((\alpha^{-1}\beta)^n) = \begin{pmatrix} 0 & 0 \\ 0 & x^n \end{pmatrix}$. Moreover, $\varphi(\alpha(\alpha^{-1}\beta)^n) = \begin{pmatrix} 0 & x^n \\ 0 & 0 \end{pmatrix}$ and $\varphi(\alpha^{-1}(\beta\alpha^{-1})^n) = \begin{pmatrix} 0 & 0 \\ x^n & 0 \end{pmatrix}$. This shows that $\varphi$ is surjective; one can also show that the map is injective.

Thus, $R_S$ is Morita equivalent to $k[x]$. The embedding $\lambda_*\colon \mathrm{Mod}\,R_S \to \mathrm{Mod}\,R$ can be explicitly described as follows, taking into account the equivalence between $\mathrm{Mod}\,k[x]$ and the category of representations of the quiver with one vertex and one loop:

$$\begin{array}{ccc}
\mathrm{Mod}\,k[x] & \xrightarrow{\lambda_*} & \mathrm{Mod}\,k(\bullet \rightleftarrows \bullet) \\
\underset{V}{\overset{x}{\circlearrowright}} & \longmapsto & V \underset{x}{\overset{1}{\rightleftarrows}} V
\end{array} \quad (3.1)$$

Since $R$ is a hereditary ring, the universal localization $\lambda$ is homological and there is a recollement of the form (1.5). Moreover, the universal localization $\lambda$ is injective and its cokernel is isomorphic to a sum of two copies of the Prüfer module $S_\infty$ [15, Proposition 1.10]. Thus, the recollement has the form

$$\mathbf{D}(k[x]) \xleftarrow[\phantom{xx}\lambda^!\phantom{xx}]{\overset{\lambda^*}{\longleftarrow}}_{\overset{\lambda_*}{\longrightarrow}} \mathbf{D}(R) \xleftarrow[\phantom{xx}j_*\phantom{xx}]{\overset{j_!}{\longleftarrow}}_{\overset{j^*}{\longrightarrow}} \mathrm{Loc}(S_\infty), \qquad (3.2)$$

where $j_!$ is the natural embedding of $\mathrm{Loc}(S_\infty)$ into $\mathbf{D}(R)$.

By Remark 2.20, $\lambda$ is the ring epimorphism associated to the partial silting $R$-module $S_\infty$ (that is even partial tilting).

Recall now the classification of silting modules over the polynomial ring $k[x]$. Since $k[x]$ is a Dedekind domain, the only silting non-tilting module is the zero module, and [14, Corollary 6.12] tilting modules are in bijection, up to equivalence, with subsets of Max-spec $k[x]$, i.e. sets of maximal ideals, that in turn are in bijection with the elements of $k$. To each such subset $\mathcal{P}$, we associate the corresponding subset $\mathcal{U}_\mathcal{P}$ of simple objects and the tilting module $T_\mathcal{P} = k[x]_{\mathcal{U}_\mathcal{P}} \oplus k[x]_{\mathcal{U}_\mathcal{P}}/k[x]$. We want to glue the tilting objects $T_\mathcal{P}$ in Mod $k[x]$ and $S_\infty$ in $\mathrm{Loc}(S_\infty)$ (or their projective presentation, respectively $\sigma_\mathcal{P}$ and $\omega$). Notice that simple $k[x]$-modules are sent to simple regular $R$-modules by the embedding $\lambda_*$ (see (3.1)) and let $\lambda_*\mathcal{U}_\mathcal{P} = \{\lambda_*(y) \mid y \in \mathcal{U}_\mathcal{P}\}$ be the set of simple regular $R$-modules that consists of the images of all simple $k[x]$-modules in $\mathcal{U}_\mathcal{P}$.

In order to apply Corollary 2.12, we have to determine $\lambda_*(T_\mathcal{P}) = \lambda_*(k[x]_{\mathcal{U}_\mathcal{P}}) \oplus \lambda_*(k[x]_{\mathcal{U}_\mathcal{P}}/k[x])$. For the first summand, observe that iterated universal localizations are again universal localizations [66, Theorem 4.6]. Thus, $\mathrm{Add}(\lambda_*(k[x]_{\mathcal{U}_\mathcal{P}})) \cong \mathrm{Add}(R_{\{S\} \cup \lambda_*\mathcal{U}_\mathcal{P}})$, since $R_S$ and $k[x]$ are Morita equivalent. The second summand is the sum of all the injective envelopes of the $k[x]$-modules belonging to $\mathcal{U}_\mathcal{P}$. The embedding $\lambda_*$ sends them to the corresponding Prüfer $R$-modules:

$$\lambda_*(k[x]_{\mathcal{U}_\mathcal{P}}/k[x]) = \lambda_*\big(\bigoplus_{S' \in \mathcal{U}_\mathcal{P}} S'_\infty\big) \cong \bigoplus_{S' \in \lambda_*\mathcal{U}_\mathcal{P}} S'_\infty \cong R_{\lambda_*\mathcal{U}_\mathcal{P}}/R.$$

This can be directly verified from the description of $\lambda_*$ above; see also [64, Example 1.4]. Consider now hypothesis (2) of Theorem 2.10, that is, $(\lambda_*T_\mathcal{P})^{\perp_{>0}} \cap (j_!S_\infty)^{\perp_{>0}}$ is closed under coproducts. This is satisfied because $\lambda_*T_\mathcal{P} \oplus j_!S_\infty$ is a tilting object in $\mathcal{D}$ (equivalent to $R_{\{S\} \cup \lambda_*\mathcal{U}_\mathcal{P}} \oplus R_{\{S\} \cup \lambda_*\mathcal{U}_\mathcal{P}}/R$).

According to Corollary 2.12, we now have to calculate the set

$$I = \mathrm{Hom}_{\mathbf{D}(R)}(\lambda_*T_\mathcal{P}, j_!S_\infty[1]).$$

Taking into account the equality

$$\mathrm{Add}(\lambda_*T_\mathcal{P}) = \mathrm{Add}\Big(R_{\{S\} \cup \lambda_*\mathcal{U}_\mathcal{P}} \oplus \bigoplus_{S' \in \lambda_*\mathcal{U}_\mathcal{P}} S'_\infty\Big),$$

this amounts to computing the group

$$\mathrm{Ext}^1_R\big(R_{\lambda_*\mathcal{U}_\mathcal{P} \cup \{S\}}, S_\infty\big) \oplus \mathrm{Ext}^1_R\Big(\bigoplus_{S' \in \lambda_*\mathcal{U}_\mathcal{P}} S'_\infty, S_\infty\Big).$$

The summand $R_{\{S\} \cup \lambda_*\mathcal{U}_\mathcal{P}}$ embeds into the localization of $R$ at the set of all simple regular modules, which is isomorphic to the sum of two copies of the generic module $G$ (see [15, Proposition 1.10]); since $\operatorname{Ext}^1_R(G, S_\infty) = 0$ for any Prüfer object (see for instance [60, Proposition 10.1]), then also $\operatorname{Ext}^1_R(R_{\{S\} \cup \lambda_* \mathcal{U}_\mathcal{P}}, S_\infty) = 0$. The second $\operatorname{Ext}^1$-group vanishes because Prüfer modules are $\operatorname{Ext}^1$-orthogonal.

Thus, $I = 0$, the triangle (2.4) is the trivial one

$$S_\infty \to S_\infty \to 0 \to S_\infty[1],$$

and $\tilde{T}_\mathcal{P} = S_\infty$. Thus, the glued silting object (in fact, tilting) is determined, up to equivalence, by

$$\operatorname{Add}(\lambda_* T_\mathcal{P} \oplus \tilde{T}_\mathcal{P}) = \operatorname{Add}(R_{\{S\} \cup \lambda_* \mathcal{U}_\mathcal{P}} \oplus \bigoplus_{S' \in \lambda_* \mathcal{U}_\mathcal{P}} S'_\infty \oplus S_\infty)$$
$$= \operatorname{Add}(R_{\{S\} \cup \lambda_* \mathcal{U}_\mathcal{P}} \oplus R_{\{S\} \cup \lambda_* \mathcal{U}_\mathcal{P}}/R).$$

By letting $\mathcal{P}$ run through through all subsets of Max-spec $k[x]$ we retrieve all large silting $R$-modules having $S_\infty$ as a direct summand. We can then let $S$ run through all simple regular modules and construct a similar recollement localizing at $S$. In this way we get all silting $R$-modules corresponding to a localization at a nonempty set of simple regular modules, and thus all large silting modules except the Lukas.

## 3.3 Gluing compact silting objects

We consider now a different family of recollements of $\mathbf{D}(R)$, namely, the ones induced by compact partial silting objects. We want to glue 2-term silting complexes through Corollary 2.19. We will use the following lemma.

**Lemma 3.1.** *[12, Theorem 5.8] Let $R$ be a hereditary ring and $\lambda \colon R \to S$ a homological ring epimorphism. Then $S \oplus \operatorname{Coker} \lambda$ is a silting $R$-module.*

*Remark* 3.2. If $\operatorname{Ker} \operatorname{Ext}^1_R(S, -)$ is closed under coproducts (e.g., if $S$ is finitely presented as $R$-module), we can deduce the statement of Lemma 3.1 from the gluing procedure of Corollary 2.19. Consider indeed a projective resolution $0 \to P_{-1} \xrightarrow{\sigma} P_0 \to S \to 0$ of $S$ as $R$-module. As in the proof of Remark 2.20, we have the following isomorphisms of triangles in $\mathbf{D}(R)$

$$\begin{array}{ccccccc}
R & \longrightarrow & \sigma & \longrightarrow & \omega & \longrightarrow & R[1] \\
\parallel & & \downarrow \cong & & \downarrow \cong & & \parallel \\
R & \xrightarrow{\lambda} & S & \longrightarrow & \operatorname{Cone}(\lambda) & \longrightarrow & R[1]
\end{array}$$

and the recollement (1.5) induced by $\lambda$ has the form

$$\mathbf{D}(S) \xrightarrow[\longleftarrow]{\overset{\longleftarrow}{\lambda_*}} \mathbf{D}(R) \xrightarrow[\longleftarrow j_*]{\overset{\longleftarrow j_!}{j^*}} \operatorname{Loc}(\omega).$$

Notice that $H^0(\operatorname{Cone} \lambda) = \operatorname{Coker} \lambda$, $H^{-1}(\operatorname{Cone} \lambda) = \operatorname{Ker} \lambda$ and $H^i(\operatorname{Cone} \lambda) = 0$ if $i \neq 0, -1$. We glue the silting objects $S$ in $\mathbf{D}(S)$ and $\omega$ in $\operatorname{Loc}(\omega)$. Notice that $I = \operatorname{Hom}_{\mathbf{D}(R)}(\lambda_* S, j_!\omega[1]) \cong \operatorname{Ext}^1_R(S, \operatorname{Coker} \lambda) = 0$. Indeed, $\operatorname{Ext}^1_R(S, S^{(J)}) \cong \operatorname{Ext}^1_S(S, S^{(J)}) = 0$ for any set $J$ and hence $\operatorname{Add}(S) \subseteq S^{\perp_1}$; but since $R$ is a hereditary ring, the covariant $\operatorname{Ext}^1_R$ functor is right exact and thus $\operatorname{Gen} S \subseteq S^{\perp_1}$. In particular, $\operatorname{Coker} \lambda \in S^{\perp_1}$. Moreover, $\lambda_* S \cong \sigma$ satisfies condition (1) of Corollary 2.19. Indeed, $\sigma^{\perp_{>0}}$ is closed under coproducts since so is $\mathcal{D}_\sigma = \operatorname{Ext}^1_R(S,-)$ by assumption.

Since $I = 0$, the triangle (2.4) of Theorem 2.10 has the form

$$j_!\omega \to j_!\omega \to 0 \to j_!\omega[1],$$

and the glued silting object is $j_!\omega \oplus \sigma$ (note that $\sigma$ is quasi-isomorphic in $\mathbf{D}(R)$ to $\lambda_*(S)$). This is indeed a 2-term complex and its zero-cohomology $H^0(j_!\omega \oplus \sigma) = \operatorname{Coker} \lambda \oplus S$ is a silting object.

In Table 1 we list the compact partial silting modules and describe the ring epimorphisms and the recollements induced. For the sake of completeness we also include the already studied case relative to the non-compact partial silting object $S_\infty$, for a simple regular object $S$.

**Table 1:** The recollements associated to partial silting Kronecker modules (see also [12, Example 5.18]).

| 2-term partial silting complex | Partial silting $R$-module | Ring epimorphism | Recollement |
|---|---|---|---|
| $P_1 \to 0$ <br> $0 \to P_1$ | $0$ <br> $P_1$ | $R \to R/Re_1R$ | $\operatorname{Loc}(Q_1) \equiv \mathbf{D}(R) \equiv \operatorname{Loc}(P_1)$ |
| $P_2 \to 0$ <br> $0 \to P_2$ | $0$ <br> $P_2$ | $R \to R/Re_2R$ | $\operatorname{Loc}(P_1) \equiv \mathbf{D}(R) \equiv \operatorname{Loc}(P_2)$ |
| $P_1^{i+1} \to P_2^i$ | $Q_i, i \geq 1$ | $R \to R_{\{Q_i\}}$ | $\operatorname{Loc}(Q_{i+1}) \equiv \mathbf{D}(R) \equiv \operatorname{Loc}(Q_i)$ |
| $P_1^{i-2} \to P_2^{i-1}$ | $P_i, i \geq 3$ | $R \to R_{\{P_i\}}$ | $\operatorname{Loc}(P_{i-1}) \equiv \mathbf{D}(R) \equiv \operatorname{Loc}(P_i)$ |
| proj. resolution | $S_\infty$ | $R \to R_S$ | $\mathbf{D}(R_S) \equiv \mathbf{D}(R) \equiv \operatorname{Loc}(S_\infty)$ |

## 3.4 Universal localization at a projective module

We consider the two recollements induced by the universal localization at a projective module. Recall that, for an idempotent $e$, defined $P = eR$, the right $R$-module $ReR$ can be computed as the trace $\tau_P(R)$ of the projective module $P$ into $R$. If $e_1$ and $e_2$ denote the idempotent elements of $R$ such that $e_1R = P_1$ and $e_2R = P_2$, we compute the quotient rings $R/Re_1R$ and $R/Re_2R$ as right $R$-modules. We have

$$\begin{aligned} Re_1R &\cong \tau_{P_1}(R) \\ &\cong \tau_{P_1}(P_1) \oplus \tau_{P_1}(P_2) \\ &\cong P_1 \oplus P_1^2 = P_1^3. \end{aligned}$$

Thus, since $R_R = P_1 \oplus P_2$ has dimension vector (1 3) and $P_1$ has dimension vector (0 1), the quotient $R/P_1^3$ has dimension vector (1 0) and it is isomorphic to $Q_1$.

Similarly,
$$\begin{aligned} Re_2R &\cong \tau_{P_2}(R) \\ &\cong \tau_{P_2}(P_1) \oplus \tau_{P_2}(P_2) \\ &\cong 0 \oplus P_2 = P_2. \end{aligned}$$

Since $R_R \cong P_1 \oplus P_2$, then the quotient $R/P_2$ is isomorphic to $P_1$.

Now, the silting object induced by the ring epimorphism $\lambda \colon R \to R/Re_1R$ according to Lemma 3.1 is $R/Re_1R \oplus \operatorname{Coker}\lambda \cong Q_1$; it corresponds to gluing the 2-term partial silting complexes $P_1 \to 0$ in $\operatorname{Loc}(P_1)$ and $P_1^2 \to P_2$ in $\operatorname{Loc}(Q_1)$. Similarly, the silting object induced by the ring epimorphism $\lambda \colon R \to R/Re_2R$ is $R/Re_2R \oplus \operatorname{Coker}\lambda \cong P_1$; it corresponds to gluing the 2-term partial silting complexes $P_2 \to 0$ in $\operatorname{Loc}(P_2)$ and $P_1$ in $\operatorname{Loc}(P_1)$.

For the recollement induced by the ring epimorphism $R \to R/Re_1R$, we may also choose the 2-term silting complexes $0 \to P_1$ in $\operatorname{Loc}(P_1)$ and $P_1^2 \to P_2$ in $\operatorname{Loc}(Q_1)$. In this case,
$$\begin{aligned} \operatorname{Hom}_{\mathbf{D}(R)}(\lambda_*Q_1, j_!P_1[1]) &\cong \operatorname{Ext}^1_R(Q_1, P_1) \\ &\cong D\operatorname{Hom}(P_1, \tau Q_1) \\ &= D\operatorname{Hom}(P_1, Q_3) \end{aligned}$$

has $k$-dimension 2. Indeed, the dimension vectors of $P_1$ and $Q_3$ are, respectively, (0 1) and (3 2), and hence $\operatorname{Hom}_R(P_1, Q_3) \cong \operatorname{Hom}_k(k, k^2)$. Thus, the set $I$ can be taken to have cardinality 2 and the triangle defining $\tilde{\sigma}$ is
$$P_1 \to \tilde{\sigma} \to Q_1^2 \to P_1[1],$$
corresponding to the short exact sequence
$$0 \to P_1 \to Q_2 \to Q_1^2 \to 0.$$
Thus, $\tilde{\sigma} = Q_2$ and the glued silting object is $Q_1 \oplus Q_2$.

Similarly, for the recollement induced by the ring epimorphism $R \to R/Re_2R$, we may choose the 2-term silting complexes $0 \to P_2$ in $\operatorname{Loc}(P_2)$ and $P_1 \to 0$ in $\operatorname{Loc}(P_1)$. Then $\dim_k \operatorname{Hom}(P_1[1], P_2[1]) = 2$, the triangle defining $\tilde{\sigma}$ is
$$P_2 \to \tilde{\sigma} \to P_1^2[1] \to P_2[1],$$
that yields $\tilde{\sigma} \cong Q_1$. The glued silting module is $Q_1$ with respect to the 2-term silting complex $(P_1^2 \to P_2) \oplus P_1[1]$.

## 3.5 Universal localization at a preprojective non projective

Consider now the family of recollements
$$\operatorname{Loc}(P_{i-1}) \xrightarrow{\lambda_*} \mathbf{D}(R) \xrightarrow{\substack{\leftarrow j_! \\ \leftarrow j^* \\ \leftarrow j_*}} \operatorname{Loc}(P_i)$$

induced by the partial silting object $P_i$, for $i \geq 2$.

By [67] and [47, Lemma 4.1] (compare [12, §5]), the universal localization $\lambda\colon R \to R_{\{P_i\}}$ is injective for $i \geq 3$ since $P_i$ is finitely presented with no projective summands. Moreover, the cokernel $\operatorname{Coker}\lambda$ is filtered by $P_i$, and thus it belongs to $\operatorname{Add}(P_i)$. Since $R_{\{P_i\}}$ belongs to $(\operatorname{Coker}\lambda)^{\perp_{0,1}} = \operatorname{Add}(P_{i-1})$, then the glued silting object $R_{\{P_i\}} \oplus \operatorname{Coker}\lambda$ of Lemma 3.1 is equivalent to the compact tilting object $P_{i-1} \oplus P_i$.

For $i = 3$ we may also choose the 2-term silting complexes $0 \to P_3$ in $\operatorname{Loc}(P_3)$ and $P_2 \to 0$ in $\operatorname{Loc}(P_2)$. Then $\dim_k \operatorname{Hom}(P_2[1], P_3[1]) = 2$, the triangle defining $\tilde{\sigma}$ is

$$P_3 \to \tilde{\sigma} \to P_2^2[1] \to P_3[1],$$

that implies $\tilde{\sigma} \cong P_1[1]$. The glued silting module is 0 with respect to the 2-term silting complex $P_1[1] \oplus P_2[1]$.

## 3.6 Universal localization at a preinjective

Finally, consider the recollement

$$\operatorname{Loc}(Q_{i+1}) \xrightarrow{\;\;\lambda_*\;\;} \mathbf{D}(R) \xrightarrow{\;\;j^*\;\;} \operatorname{Loc}(Q_i)$$

induced by the ring epimorphism $\lambda\colon R \to R_{\{Q_i\}}$, for $i \geq 1$. Similarly to the preprojective non-projective case, by Lemma 3.1 we get the silting module $R_{Q_i} \oplus \operatorname{Coker}\lambda$, which is equivalent to $Q_{i+1} \oplus Q_i$.

This exhausts the list of compact silting $R$-modules. Apparently, the only way of getting the Lukas tilting module $L$ by means of gluing is through a trivial recollement. Indeed, this module enjoys the remarkable property that for any direct summand $M$, $\operatorname{Add}(M) = \operatorname{Add}(L)$ (see [4]). Thus, any partial tilting module $M$ that is a direct summand of $L$ generates the same localizing subcategory $\operatorname{Loc}(M) = \operatorname{Loc}(L) = \mathbf{D}(R)$ and one gets the trivial recollement

$$0 \longrightarrow \mathbf{D}(R) \xrightarrow{\;\;j^*\;\;} \operatorname{Loc}(L)\ .$$

The Lukas module can then be obtained by gluing the silting modules $L$ in $\operatorname{Loc}(L)$ and 0 in the zero category.

# Expansions of abelian categories

## 4

In Chapter 3 we employed recollements induced by universal localizations to get the derived category of Kronecker modules from the derived module category of the polynomial ring and eventually get as many silting modules as possible through gluing. Now we want to retrieve silting objects on concealed canonical algebras starting from the classification of silting modules over the Kronecker algebra. Indeed, concealed canonical algebras are related to the Kronecker algebra by iterated universal localizations at one-element sets of simple regular objects. Expansions of abelian categories offer the geometric perspective for the same phenomenon.

## 4.1 Left and right expansions

We recall here briefly the construction of expansions of abelian categories, a concept developed by Chen and Krause [23, 24]. Roughly speaking, an expansion is a fully faithful and exact functor $\mathcal{B} \to \mathcal{A}$ between abelian categories that admits an exact left adjoint and an exact right adjoint. In addition one requires the existence of simple objects $S_\lambda$ and $S_\rho$ in $\mathcal{A}$ such that $S_\lambda^{\perp_{0,1}} = \mathcal{B} = {}^{\perp_{0,1}}S_\rho$, where $\mathcal{B}$ is viewed as a full subcategory of $\mathcal{A}$. In fact, these simple objects are related by an exact sequence $0 \to S_\rho \to S \to S_\lambda \to 0$ in $\mathcal{A}$ such that $S$ is a simple object in $\mathcal{B}$. On the other hand, the left adjoint of $\mathcal{B} \to \mathcal{A}$ identifies $S_\rho$ with $S$, whereas the right adjoint identifies $S_\lambda$ with $S$.

Let $\mathcal{A}$ be an abelian category. Recall that a full subcategory $\mathcal{B}$ of $\mathcal{A}$ is called exact abelian if $\mathcal{B}$ is an abelian category and the inclusion functor is exact. Now let $i\colon \mathcal{B} \to \mathcal{A}$ be a fully faithful and exact functor between abelian categories. It is convenient to identify $\mathcal{B}$ with the essential image of $i$, which means that $\mathcal{B}$ is an exact abelian subcategory of $\mathcal{A}$.

**Definition 4.1.** [24, Definition 3.1.1] A fully faithful end exact functor $i\colon \mathcal{B} \to \mathcal{A}$ between abelian categories is called *left expansion* if the following conditions are satisfied:

(E1) The functor $i\colon \mathcal{B} \to \mathcal{A}$ admits an exact left adjoint.

(E2) The category $^{\perp_{0,1}}\mathcal{B}$ is equivalent to $\operatorname{mod}\Delta$ for some division ring $\Delta$.

(E3) $\operatorname{Ext}^2_{\mathcal{A}}(A,B) = 0$ for all $A, B \in {}^{\perp_{0,1}}\mathcal{B}$.

The functor $\mathcal{B} \to \mathcal{A}$ is called *right expansion* if the dual conditions are satisfied.

Left expansions admit the following characterization.

**Lemma 4.2.** *[24, Lemma 3.1.4] The following are equivalent for an abelian category $\mathcal{A}$ and an exact abelian subcategory $\mathcal{B}$:*

(1) *The inclusion $\mathcal{B} \to \mathcal{A}$ is a left expansion.*

(2) *There exists a simple object $S \in \mathcal{A}$ such that $S^{\perp_{0,1}} = \mathcal{B}$ that satisfies the following properties:*

   a) $\operatorname{Hom}_{\mathcal{A}}(S,A)$ and $\operatorname{Ext}^1_{\mathcal{A}}(S,A)$ *are of finite length over* $\operatorname{End}_{\mathcal{A}}(S)$ *for all $A \in \mathcal{A}$,*

   b) $\operatorname{Ext}^1_{\mathcal{A}}(S,S) = 0$ *and* $\operatorname{Ext}^2_{\mathcal{A}}(S,A) = 0$ *for all $A \in \mathcal{A}$.*

Dually, one can define *right expansions*, which admit a dual characterization.

**Definition 4.3.** A fully faithful and exact functor $i\colon \mathcal{B} \to \mathcal{A}$ between abelian categories is an *expansion* of abelian categories if the functor $i$ is a left and right expansion.

Let $i\colon \mathcal{B} \to A$ be an expansion of abelian categories. We identify $\mathcal{B}$ with the essential image of $i$. We denote by $i_\lambda$ the left adjoint of $i$ and by $i_\rho$ the right adjoint of $i$. We choose an indecomposable object $S_\lambda$ in $^{\perp_{0,1}}\mathcal{B}$ and an indecomposable object $S_\rho$ in $\mathcal{B}^{\perp_{0,1}}$. Thus $^{\perp_{0,1}}\mathcal{B} = \operatorname{add}(S_\lambda)$ and $\mathcal{B}^{\perp_{0,1}} = \operatorname{add}(S_\rho)$ by property (E2).

An expansion $i\colon \mathcal{B} \to \mathcal{A}$ is called *split* if $\mathcal{B}^{\perp_{0,1}} = {}^{\perp_{0,1}}\mathcal{B}$; it is split if and only if ([24, Lemma 3.2.4]) $\mathcal{A}$ decomposes as the disjoint union $\mathcal{B} \amalg \mathcal{C}$ for some Serre subcategory $\mathcal{C}$. If the expansion is non-split, ([24, Lemma 3.2.5]) the object $\bar{S} = i_\lambda(S_\rho)$ is simple in $\mathcal{B}$ and isomorphic to $i_\rho(S_\lambda)$. Moreover, there exists a short exact sequence
$$0 \to S_\rho \to i\bar{S} \to S_\lambda \to 0$$
and the functor $i_\lambda$ (resp. $i_\rho$) induces an equivalence $\mathcal{B}^{\perp_{0,1}} \xrightarrow{\sim} \operatorname{add}(\bar{S})$ (resp. $^{\perp_{0,1}}\mathcal{B} \xrightarrow{\sim} \operatorname{add}(\bar{S})$). In particular, $^{\perp_{0,1}}\mathcal{B} \cong \mathcal{B}^{\perp_{0,1}}$.

An expansion $\mathcal{B} \to \mathcal{A}$ of abelian categories determines a division ring $\Delta$ such that $^{\perp_{0,1}}\mathcal{B}$ and $\mathcal{B}^{\perp_{0,1}}$ are equivalent to $\operatorname{mod}\Delta$: we call $\Delta$ the *associated division ring*.

Let $i\colon \mathcal{B} \to \mathcal{A}$ be an expansion with associated division ring $\Delta$, so that $^{\perp_{0,1}}\mathcal{B} \cong \operatorname{mod}\Delta \cong \mathcal{B}^{\perp_{0,1}}$. Then there are inclusions $j\colon {}^{\perp_{0,1}}\mathcal{B} \hookrightarrow \mathcal{A}$ and $k\colon \mathcal{B}^{\perp_{0,1}} \hookrightarrow \mathcal{A}$ with adjoints $j_\rho$ and $k_\lambda$, yielding a diagram

$$\mathcal{B} \underset{\xleftarrow{i_\rho}}{\overset{\xleftarrow{i_\lambda}}{\underset{\longrightarrow}{\xrightarrow{i}}}} \mathcal{A} \underset{\xleftarrow{k}}{\overset{\xleftarrow{j}}{\underset{\xrightarrow{k_\lambda}}{\xrightarrow{j_\rho}}}} \operatorname{mod}\Delta \ .$$

When one passes from the abelian categories to their derived categories, the diagram above induces a recollement of triangulated categories [24, Proposition 3.3.2]

$$\mathbf{D}^b(\mathcal{B}) \xrightarrow[\overset{\longleftarrow i^*}{\underset{\longleftarrow i^!}{\overset{i_*}{\longrightarrow}}}]{} \mathbf{D}^b(\mathcal{A}) \xrightarrow[\overset{\longleftarrow j_!}{\underset{\longleftarrow j_*}{\overset{j^*}{\longrightarrow}}}]{} \mathbf{D}^b(\mathrm{mod}\,\Delta), \tag{4.1}$$

where $i_*, i^*, i^!, j_!, j_*$ are induced by the corresponding functors $i, i_\lambda, i_\rho, j, k$ between abelian categories, while $j^*$ is the right adjoint of $j_!$, that is isomorphic to the left adjoint of $j_*$.

The following results describe the functor $i\colon \mathcal{B} \to \mathcal{A}$ and its adjoints with respect to simple objects. The left adjoint $i_\lambda$ induces a bijection between the isomorphism classes of simple objects of $\mathcal{A}$ that are different from $S_\lambda$, and the isomorphism classes of simple objects of $\mathcal{B}$.

**Lemma 4.4.** *[24, Lemma 3.4.1] Let $i\colon \mathcal{B} \to \mathcal{A}$ be an expansion of abelian categories.*

(1) *If $S$ is a simple object in $\mathcal{B}$ and $S \not\cong \bar{S}$, then $iS$ is simple in $\mathcal{A}$ and $i_\lambda iS \cong S$.*

(2) *If $S$ is a simple object in $\mathcal{A}$ and $S \not\cong S_\lambda$, then $i_\lambda S$ is simple in $\mathcal{B}$. Moreover, $S \cong ii_\lambda S$ if $S \not\cong S_\rho$.*

## 4.2 Coherent and quasicoherent sheaves

Throughout this section, let $k$ be an algebraically closed field. Weighted projective curves were introduced in [29]; following [9, 41], a noncommutative curve $\mathbb{X}$ can be described axiomatically by a category $\mathcal{H}$ that is regarded as the category $\mathrm{coh}\,\mathbb{X}$ of coherent sheaves over $\mathbb{X}$:

(NC1) $\mathcal{H}$ is small, connected, abelian and every object in $\mathcal{H}$ is noetherian;

(NC2) $\mathcal{H}$ is a $k$-category with finite-dimensional Hom and Ext spaces;

(NC3) $\mathcal{H}$ admits a Serre functor, that is an autoequivalence $\tau$ such that the Serre duality
$$\mathrm{Ext}^1_\mathcal{A}(X,Y) \cong D\,\mathrm{Hom}_\mathcal{A}(Y,\tau X)$$
holds, where $D = \mathrm{Hom}_k(-,k)$.

(NC4) $\mathcal{H}$ contains an object of infinite length.

Notice that (NC3) is equivalent to ask that $\mathcal{A}$ be hereditary and has no nonzero projective objects [23, Proposition 3.4.5].

For a category $\mathcal{H}$ of coherent sheaves there is a positively $H$-graded noetherian ring $R$, with $(H, \leq)$ an ordered abelian group of rank one, such that $\mathcal{H}$ is

$$\mathcal{H} = \frac{\mathrm{mod}^H(R)}{\mathrm{mod}_0^H(R)},$$

the quotient category of finitely generated $H$-graded modules modulo the Serre subcategory of those modules which are finite-dimensional over $k$.

This description can be used to define the category $\overrightarrow{\mathcal{H}} = \operatorname{Qcoh} \mathbb{X}$ of quasi-coherent sheaves as

$$\overrightarrow{\mathcal{H}} = \frac{\operatorname{Mod}^H(R)}{\operatorname{Mod}^H_0(R)},$$

where $\operatorname{Mod}^H_0(R)$ denotes the localizing subcategory of $\operatorname{Mod}^H(R)$ of all locally finite-dimensional modules, that is, direct limits of finite-dimensional modules over $k$. The category $\overrightarrow{\mathcal{H}}$ is hereditary abelian, and a locally noetherian Grothendieck category; every object in $\overrightarrow{\mathcal{H}}$ is a direct limit of objects in $\mathcal{H}$. The full abelian subcategory $\mathcal{H}$ consists of the coherent (= finitely presented = noetherian) objects in $\overrightarrow{\mathcal{H}}$, and we also write $\mathcal{H} = \operatorname{fp}(\overrightarrow{\mathcal{H}}) = \operatorname{coh}(\overrightarrow{\mathcal{H}})$. Every indecomposable coherent sheaf has a local endomorphism ring, and $\mathcal{H}$ is a Krull-Schmidt category. The closure under direct limits of a tube $\mathcal{U}$ in $\overrightarrow{\mathcal{A}}$ will be denoted by $\overrightarrow{\mathcal{U}}$.

If $\mathbb{X}$ is the non-weighted projective line $\mathbb{P}^1_k$ over $k$, then $H = \mathbb{Z}$ and $R = k[x_0, x_1]$ (see [23]).

Let $\mathcal{H}_0$ be the Serre subcategory of $\mathcal{H}$ consisting of all objects of finite length. Then

$$\mathcal{H}_0 = \coprod_{x \in \mathbb{X}} \mathcal{U}_x$$

(for some index set $\mathbb{X}$), where $\mathcal{U}_x$ are connected uniserial categories, called *tubes*. A category $\mathcal{U}$ that is equivalent to a tube $\mathcal{U}_x$, for some weighted projective curve $\mathbb{X}$ and some $x \in \mathbb{X}$, will be called a *tube category*. By [45, Proposition 1.1], every coherent sheaf in $\mathcal{H}$ belongs either to $\mathcal{H}_0$ or to $\mathcal{H}_+ = \operatorname{vect} \mathbb{X}$, that is the class of sheaves not containing any simple sheaf, also called *vector bundles* or *torsionfree sheaves*. Objects in $\mathcal{H}_0$ will be called *torsion sheaves*. We also require any curve to satisfy the following property:

(NC5) $\mathbb{X}$ consists of infinitely many points.

The quotient category $\mathcal{H}/\mathcal{H}_0$ is, by [45, Proposition 3.4], of the form $\mathcal{H}/\mathcal{H}_0 \cong \operatorname{mod}(k(\mathbb{X}))$ for a unique division ring $k(\mathcal{H})$, called the *function field* of $\mathcal{H}$ (or $\mathbb{X}$). If we denote by $\overrightarrow{\mathcal{H}}_0$ the closure under direct limits of $\mathcal{H}_0$ in $\mathcal{H}$, then $\overrightarrow{\mathcal{H}}/\overrightarrow{\mathcal{H}}_0 \cong \operatorname{Mod} k(\mathcal{H})$.

If $\mathcal{H}$ is a category satisfying (NC1) through (NC5), then we call $\mathbb{X}$ a *weighted noncommutative regular projective curve* over $k$. It is shown in [41, Theorem 7.11] that for a weighted noncommutative regular projective curve, for all $x \in \mathbb{X}$ there are (up to isomorphism) precisely $p(x) < \infty$ simple objects in $\mathcal{U}_x$ and for almost all $x$, $p(x) = 1$. The numbers $p(x)$ are called *weights*; points $x \in \mathbb{X}$ are called *exceptional* if $p(x) > 1$, *homogeneous* otherwise. In case $p(x) = 1$ for each $x \in \mathbb{X}$, we say that the curve $\mathbb{X}$ is *non-weighted*. Notice that $\mathcal{H}$ is non-weighted if and only if for any simple object $S$, $\operatorname{Ext}^1_{\mathcal{H}}(S, S) \neq 0$ (equivalently, $\tau S \cong S$). Finally, we say $\mathcal{H}$ is *of genus* 0 if it satisfies the following condition:

(G0) $\mathcal{H}$ admits a tilting object.

## 4.3 A geometric model for the tube category

We describe here a geometric model by Baur, Buan and Marsh [16] for objects in the closure $\overrightarrow{\mathcal{U}}$ under direct limits of a tube category $\mathcal{U}$. This model provides a handy way of visualizing modules and computing Hom and Ext groups in such categories, and gives a characterization of maximal exceptional objects. Recall an object $T$ in an abelian category $\mathcal{A}$ is said to be *exceptional* if $\mathrm{Ext}^i_{\mathcal{A}}(T,T) = 0$ for all integers $i > 0$.

Before talking about tube categories, we present a geometric model for the category $\mathrm{mod}\, kQ$, where $Q$ is the linearly oriented quiver of type $\mathbb{A}_n$:

$$1 \longleftarrow 2 \longleftarrow 3 \longleftarrow \cdots \longleftarrow n-1 \longleftarrow n.$$

Write $S_i$, $1 \leq i \leq n$ for the simple $kQ$-modules. For $i \leq j-2$, let $M_{ij}$ denote the indecomposable $kQ$-module with composition factors $S_{i+1}, \ldots, S_{j-1}$ (starting from the socle). The modules $M_{ij}$ are uniquely determined since $\mathrm{mod}\, kQ$ is a serial category. In this notation, the simple objects are $S_i = M_{i-1\,i+1}$. If $i \geq j-1$, we set $M_{ij} = 0$.

Consider a line segment $\ell_n$ with $n+2$ marked points, numbered from $0$ to $n+1$:

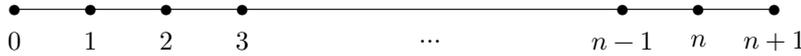

0    1    2    3         $\cdots$         $n-1$    $n$    $n+1$

and associate the module $M_{ij}$ to the (isotopy class of) arcs $[i,j]$ above $\ell_n$ from $i$ to $j$ oriented towards $j$, for $0 \leq i \leq j-2 \leq n$. Recall two arcs are *isotopic* if there exists a continuous transformation which sends one to the other, preserving orientation. This gives a bijection between indecomposable $kQ$-modules and the set $\mathcal{A}(\ell_n)$ of isotopy classes of arcs between marked points of $\ell_n$, above $\ell_n$, which are not isotopic to boundary arcs [16, §3.1]. Moreover, for arcs $[i,j]$ and $[i',j']$, $\mathrm{Ext}^1(M_{ij}, M_{i'j'}) \cong k$ if there is a negative crossing between $[i,j]$ and $[i',j']$, and 0 otherwise. By *negative crossing* we mean the arcs $[i,j]$ and $[i',j']$ cross as in Figure 1. In case $\mathrm{Ext}^1(M_{ij}, M_{i'j'}) \cong k$, then

**Figure 1:** A negative crossing between $[i,j]$ and $[i',j']$.

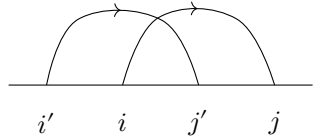

$i'$    $i$    $j'$    $j$

the non-split extension has the form

$$0 \to M_{i'j'} \to M_{i'j} \oplus M_{ij'} \to M_{ij} \to 0$$

if $j' - i \geq 2$, while if $j' = i+1$, it has the form

$$0 \to M_{i'j'} \to M_{i'j} \to M_{ij} \to 0.$$

A geometric representation of extensions is given in Figure 2, where dashed lines represent the direct summands of the middle terms of the short exact sequences above.

The Auslander-Reiten translation $\tau$ moves an arc one step to the left (or gives zero if this is not defined). Nonzero quotients of an indecomposable

**Figure 2:** Non-split extensions.

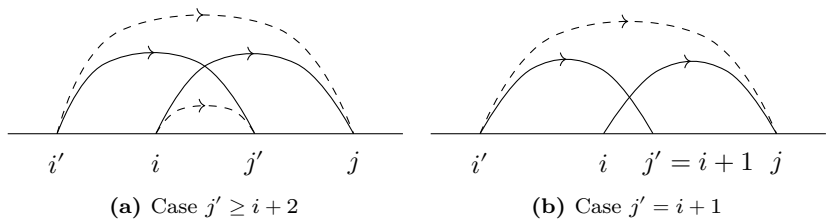

(a) Case $j' \geq i+2$      (b) Case $j' = i+1$

module $M_{ij}$ are modules $M_{i'j}$, with $i \leq i' \leq j-2$, and thus represented by arcs sharing the ending vertex $j$ and with shorter length. Similarly, nonzero subobjects of an indecomposable module $M_{ij}$ are modules $M_{ij'}$, with $i+2 \leq j' \leq j$, and thus represented by arcs sharing the starting vertex $i$ and with shorter length.

Assume now $\mathcal{U}$ is a tube of rank $n$. Denote the simples in $\vec{\mathcal{U}}$ by $S_i$, for $1 \leq 0 \leq n-1$, ordered in such a way that $\tau S_i = S_{i-1}$ for $1 \leq i \leq n-1$ and $\tau S_0 = S_{n-1}$. For $-1 \leq i \leq n-2$ and $j \geq i+2$, denote by $M_{ij}$ the unique indecomposable object in $\vec{\mathcal{U}}$ with simple socle $S_{i+1}$ and length $j-i-1$. For any $k \in \mathbb{Z}$, let $M_{i+kn\ j+kn} = M_{ij}$. With similar notation, denote by $M_{i\infty}$ the Prüfer module with simple socle $S_{i+1}$, for $-1 \leq i \leq n-1$. Note that the indecomposable objects in $\vec{\mathcal{U}}$ are exactly the indecomposables $M_{ij}$ in $\vec{\mathcal{U}}$ and the Prüfer modules $M_{i\infty}$ (see for example [9, Corollary 3.7(4)]).

The geometric model for $\vec{\mathcal{U}}$ extends the geometric model for $\mathcal{U}$ presented above. Consider an annulus $\mathbf{A}_n$ with $n$ marked points on the outer boundary. The points are labeled $0, 1, \ldots, n-1$ and arranged counterclockwise, as in Figure 3. Let $\mathbf{U}_n$ be the universal cover of $\mathbf{A}_n$, that is, a strip $\mathbb{R} \times [0, 1]$ with marked points corresponding to $(k, 0)$, for all $k \in \mathbb{Z}$, together with a continuous surjective map $\pi_n \colon \mathbf{U}_n \to \mathbf{A}_n$ such that for every $x \in \mathbf{A}_n$, there exists an open neighborhood $U$ of $x$, such that $\pi_n^{-1}(U)$ is a union of disjoint open sets in $\mathbf{U}_n$, each of which is mapped homeomorphically onto $U$ by $\pi_n$. The marked points on $\mathbf{U}_n$ are mapped to the corresponding marked points on $\mathbf{A}_n$ (modulo $n$). See Figure 4.

**Figure 3:** The annulus $\mathbf{A}_n$.

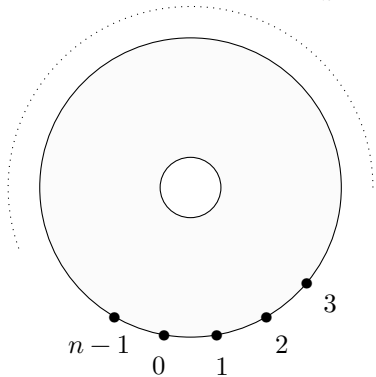

For integers $i, j$, with $j - i \leq 2$, denote by $[i, j]$ the arc in $\mathbf{U}_n$ with starting point $i$ and ending point $j$, oriented from $i$ to $j$. We also allow arcs of the form $[i, \infty]$ with starting point at $i$ and oriented in the positive $x$ direction (see Figure 4). Such arcs will be called *(right) infinite arcs*.

**Figure 4:** The universal cover $\mathbf{U}_n$ with a (right) infinite arc.

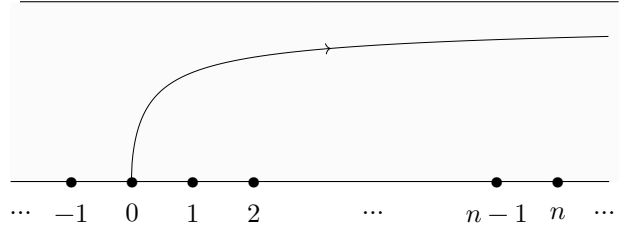

Let $\pi_n([i,j])$ be the corresponding arc in $\mathbf{A}_n$ and let $\vec{\mathcal{A}} = \vec{\mathcal{A}}(\mathbf{A}_n)$ denote the set of (isotopy classes of) such arcs. It contains the set $\mathcal{A} = \mathcal{A}(\mathbf{A}_n)$ of arcs of the form $\pi_n([i,j])$ with $i, j$ finite. The map $\psi \colon \vec{\mathcal{A}} \to \vec{\mathcal{U}}$ sending $\pi_n([i,j])$ to $M_{ij}$ induces a bijection between $\vec{\mathcal{A}}$ and the class of indecomposable objects in $\vec{\mathcal{U}}$ (see [16, §4.2]).

**Figure 5:** A (right) infinite arc in $\mathbf{A}_n$.

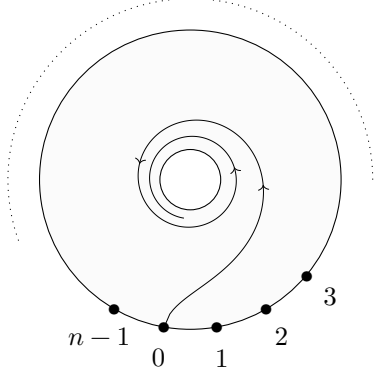

Define a quiver $Q(\mathcal{A})$ with vertices given by the element in $\mathcal{A}$ and arrows
$$\pi_n([i,j]) \longrightarrow \pi([i, j+1]) \quad \text{and}$$
$$\pi_n([i,j]) \longrightarrow \pi([i+i, j]) \quad (\text{if } j > i + 2).$$
The map $\tau \colon \mathcal{A} \to \mathcal{A}$ defined using the formula $\tau(\pi_n([i,j])) = \pi_n([i-1, j-1])$ is a *translate*, that means, it is a bijection and for arcs $\alpha$ and $\beta$ in $\mathcal{A}$, the number of arrows from $\beta$ to $\alpha$ equals the number of arrows from $\tau(\alpha)$ to $\beta$. A quiver together with a translate is called *translation quiver* and the pair $(Q(\mathcal{A}), \tau)$ is called the *(translation) quiver of* $\mathcal{A}(\mathbf{A}_n)$. The restriction of $\psi$ to $\mathcal{A}$ induces an isomorphism between the translation quiver of $\mathcal{A}(\mathbf{A}_n)$ and the Auslander-Reiten quiver of $\mathcal{U}_n$ ([16, Proposition 4.2]).

For arcs $\alpha, \beta$ in $\mathcal{A}(\mathbf{A}_n)$, let $I(\alpha, \beta)$ be the minimum number of intersections between arcs in the isotopy classes of $\alpha$ and $\beta$, not allowing multiple intersections. Similarly, let $I^+(\alpha, \beta)$ (resp. $I^-(\alpha, \beta)$) denote the number of positive (resp. negative) crossings between $\alpha$ and $\beta$ (see Figure 1 for an example of negative crossing). Then the following formula holds for two indecomposable objects $M_{ij}$ and $M_{i'j'}$ in $\overrightarrow{\mathcal{U}}$ [16, Theorem 4.3]:

$$\operatorname{Ext}^1(M_{ij}, M_{i'j'}) \cong \prod_{I^-(\pi_n([i,j]), \pi_n([i',j']))} k. \tag{4.2}$$

Note finally that the map $\psi \colon \pi_n([i,j]) \mapsto M_{ij}$ induces a bijection between maximal collections of mutually noncrossing arcs in $\mathbf{A}_n$ (including right infinite arcs) and maximal exceptional objects in $\overrightarrow{\mathcal{U}}$ [16, Remark 6.3]. Moreover, maximal exceptional objects in $\overrightarrow{\mathcal{U}}$ are exactly exceptional objects with $n$ pairwise nonisomorphic indecomposable direct summands ([21], see also [16, Proposition 5.1]).

## 4.4 Expansions of categories of coherent sheaves

We describe now expansions of categories of coherent sheaves. We start by talking about expansions of tube categories.

Let $\mathcal{U}$ be a tube of rank $n \geq 2$ and fix a simple object $S_\lambda$ of $\mathcal{U}$. This object satisfies conditions (a) and (b) of Lemma 4.2(2), as well as the dual conditions; thus, the inclusion $i \colon S_\lambda^{\perp_{0,1}} \to \mathcal{U}$ is an expansion and $S_\lambda^{\perp_{0,1}}$ is a tube category of rank $n-1$. In Figure 6 is represented the Auslander-Reiten quiver of the tube category $\mathcal{U}$: the vertices represent indecomposable objects and there is an arrow between two vertices if and only if there is an irreducible morphism between them. The indecomposable objects not lying in $S_\lambda^{\perp_{0,1}}$ are represented by circles. The object $\bar{S}$ is the only simple object in $\mathcal{U}$ that is not simple in $S_\lambda^{\perp_{0,1}}$. The same expansion is represented in Figure 7 in the geometric model for the tube category of Section 4.3. The image $i\bar{S}$ of the simple object $\bar{S}$ in $S_\lambda^{\perp_{0,1}}$ is the nontrivial extension of the two simple objects $S_\lambda$ and $S_\rho = \tau S_\lambda$ in $\mathcal{U}$. Considering that $i$ is exact, it is easy to visualize the image under $i$ of any arc. Namely, each time an arc corresponding to an object $M$ passes over $\bar{S}$ in $S_\lambda^{\perp_{0,1}}$, the arc corresponding to $iM$ will pass over both $S_\rho$ and $S_\lambda$ in $\mathcal{U}$.

**Figure 6:** The Auslander-Reiten quiver of a tube category $\mathcal{U}$ of rank $n = 7$. The two dotted lines are identified.

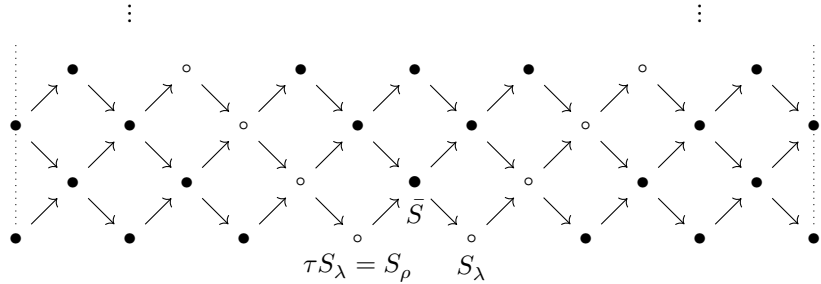

Consider now a category $\mathcal{A} = \operatorname{coh} \mathbb{X}$ of coherent sheaves over a weighted noncommutative regular projective curve of genus 0. Chen and Krause proved the following characterization of expansions of categories of coherent sheaves.

**Figure 7:** The expansion $i$ of tube categories in the geometric model. Note that the dashed lines are identified.

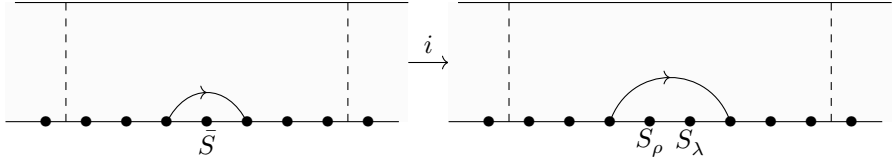

**Lemma 4.5.** *[23, Lemma 6.5.1] Let $\mathcal{A}$ be a category satisfying (NC1)–(NC3) and (G0). For a full subcategory $\mathcal{B}$ of $\mathcal{A}$, the following are equivalent:*

(1) *The inclusion $i\colon \mathcal{B} \to \mathcal{A}$ is an expansion of abelian categories.*

(2) *There exists a simple object $S_\lambda$ in $\mathcal{A}$ such that $\tau S_\lambda \not\cong S_\lambda$ and $S_\lambda^{\perp_{0,1}} = \mathcal{B}$.*

Moreover, if $i\colon \mathcal{B} \to \mathcal{A}$ is a non-split expansion of abelian categories satisfying the conditions of Lemma 4.5, then $i$ restricts to an expansion $\mathcal{B}_0 \to \mathcal{A}_0$ by [24, Proposition 3.4.3]. Furthermore, if we write $\mathcal{A}_0 = \coprod_{x \in \mathbb{X}} \mathcal{A}_x$ and $\mathcal{B}_0 = \coprod_{x \in \mathbb{X}} \mathcal{B}_x$ as decompositions into connected tubes, then there exists $\bar{x} \in \mathbb{X}$ such that $i$ restricts to an expansion $\mathcal{B}_{\bar{x}} \to \mathcal{A}_{\bar{x}}$ and to equivalences $\mathcal{B}_x \xrightarrow{\sim} \mathcal{A}_x$ for all $x \neq \bar{x}$.

We have also the following result.

**Theorem 4.6.** *[23, Theorem 6.5.4] Let $k$ be a field and $\mathcal{B} \to \mathcal{A}$ a non-split expansion of $k$-linear abelian categories with associated division ring $k$. Then $\mathcal{A}$ satisfies (NC1)–(NC5) and (G0) if and only if $\mathcal{B}$ satisfies (NC1)–(NC5) and (G0).*

## 4.5 The classification of tilting sheaves

In [9], Angeleri and Kussin provide the classification of large tilting objects in categories $\overrightarrow{\mathcal{H}} = \operatorname{Qcoh}\mathbb{X}$ of quasicoherent sheaves over weighted noncommutative regular projective curves of genus 0.

Let $\mathcal{U} = \mathcal{U}_x$ be a tube of rank $p > 1$. If $E \in \mathcal{U}$ is an indecomposable exceptional sheaf, then the collection $\mathcal{W}$ of all subquotients of $E$ is called the *wing rooted in $E$*. Notice that the full subcategory $\operatorname{add}(\mathcal{W})$ of $\mathcal{H}$ is equivalent to the category of finite-dimensional representation of the linearly oriented Dynkin quiver $\mathbb{A}_r$ ([43, Chapter 3]), for some $1 \leq r \leq p-1$ which equals the length of the *root $E$*. The set of all simple sheaves in $\mathcal{W}$ is of the form $\{S, \tau^- S, \dots, \tau^{-(r-1)} S\}$ for an exceptional simple sheaf $S$. Such a set of simples is called a *segment* and it forms the *basis* of $\mathcal{W}$. Two wings in $\mathcal{U}$ are *non-adjacent* if their bases are disjoint and their union consists of less than $p$ simples and is not a segment.

Any tilting object $B$ in the category $\operatorname{add}(\mathcal{W})$ has precisely $r$ nonisomorphic indecomposable summands $B_1, \dots, B_r$ (see [62, p. 205]). We call the object $B = B_1 \oplus \cdots \oplus B_r$ a *connected branch* in $\mathcal{W}$: one of the $B_i$ is isomorphic to $E$ and for every $j$ the wing rooted in $B_j$ contains as many summands of $B$ as the length of $B_j$. A sheaf $B$ of finite length is called a *branch sheaf* if it is a multiplicity free direct sum of connected branches in pairwise non-adjacent wings; then $\operatorname{Ext}^1_R(B,B) = 0$.

By [9, Theorem 4.8], tilting sheaves with large torsion part are parametrized by pairs $(B, V)$, where $B$ is a branch sheaf and $V$ a nonempty subset of $\mathbb{X}$. More precisely, every tilting sheaf with large torsion part is of the form $T_{(B,V)} = T_+ \oplus T_0$, where the torsionfree part $T_+$ is $V$-divisible (that is, $\operatorname{Ext}^1_R(S, T_+) = 0$ for every simple object supported in $V$) and the torsion part $T_0$ is given by

$$T_0 = B \oplus \bigoplus_{x \in V} \bigoplus_{j \in \mathcal{R}_x} \tau^j S_x[\infty],$$

where $S_x$ is a simple object in the tube $\mathcal{U}_x$ and where $\mathcal{R}_x = \{x \in \{0, \ldots, p(x) - 1\} \mid \tau^{j+1} S_x \notin \mathcal{W}\}$. We write $B = B_\mathfrak{e} \oplus B_\mathfrak{i}$, where $B_\mathfrak{e}$ is supported in $\mathbb{X} \setminus V$ and $B_\mathfrak{i}$ in $V$. In this case we say $B_\mathfrak{e}$ is *exterior* and $B_\mathfrak{i}$ is *interior* with respect to $V$. In Figure 8 is represented an example.

**Figure 8:** An example of a large tilting sheaf, corresponding to the pair $(B = B_\mathfrak{i} \oplus B_\mathfrak{e}, V)$. The stars $\star$ denote Prüfer summands.

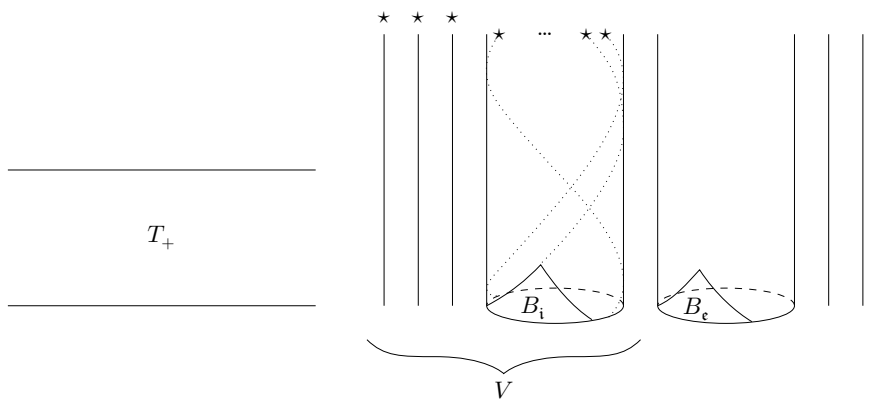

## 4.6 The recollement induced by an expansion

Consider a simple exceptional compact object $S$ in $\mathcal{A}$. Let $\mathcal{B} = S^{\perp_{0,1}}$. Then by Lemma 4.5, the embedding $i\colon \mathcal{B} \to \mathcal{A}$ of the full subcategory $\mathcal{B} = S^{\perp_{0,1}}$ is an expansion. Thus we have a recollement like (4.1)

$$\mathbf{D}^b(\mathcal{B}) \xrightleftharpoons{\longrightarrow} \mathbf{D}^b(\mathcal{A}) \xrightleftharpoons{\longrightarrow} \mathbf{D}^b(\operatorname{mod} k), \qquad (4.3)$$

where $k = \operatorname{End}_\mathcal{A}(S)$.

We may also consider $S$ as a simple exceptional object in the category $\overrightarrow{\mathcal{A}}$ of quasicoherent sheaves. Since $S$ is compact, $S$ is a partial silting (even partial tilting) object in $\overrightarrow{\mathcal{A}}$ and thus we can build the recollement (1.3) associated to $S$

$$S^{\perp_\mathbb{Z}} \underset{\xleftarrow{i^!}}{\overset{\xleftarrow{i^*}}{\underset{i_*}{\longrightarrow}}} \mathbf{D}(\overrightarrow{\mathcal{A}}) \underset{\xleftarrow{j_*}}{\overset{\xleftarrow{j_!}}{\underset{j^*}{\longrightarrow}}} \operatorname{Loc}(S) \cong \mathbf{D}(\operatorname{Mod} k). \qquad (4.4)$$

It is easy to see that the category $\mathcal{B}$ is a subcategory of the triangulated subcategory $S^{\perp_\mathbb{Z}}$ of $\mathbf{D}(\overrightarrow{\mathcal{A}})$: just recall that $\mathcal{B}$ is characterized as the full abelian subcategory of $\mathcal{A}$ of objects $X$ with $\mathrm{Hom}_\mathcal{A}(S, X) = 0 = \mathrm{Ext}^1_\mathcal{A}(S, X)$. Thus, the restriction of $i_*$ to $\mathcal{B}$ gives back the expansion $i$ of abelian categories.

## 4.7 Gluing tilting sheaves

Now we are ready to glue. Assume that $\mathcal{A}$ satisfies (NC1)–(NC5) and (G0). We will need the following results.

**Proposition 4.7.** *[9, Theorem 3.8] If $T \in \overrightarrow{\mathcal{A}}$ is a sheaf satisfying $\mathrm{Ext}^1_{\overrightarrow{\mathcal{A}}}(T, T) = 0$, then $T$ is a direct sum of its torsionfree part $T_0$ and its torsion part $T_+$.*

**Lemma 4.8.** *Let $\mathcal{D}$ be a triangulated category.*

*(1) If $A \to B \to C \oplus D \xrightarrow{\alpha} A[1]$ is a triangle in $\mathcal{D}$ and $\mathrm{Hom}_\mathcal{D}(D, A[1]) = 0$, then $B \cong B' \oplus D$, where $B'$ is the cocone of the restriction $\alpha'$ of $\alpha$ to $C$.*

*(2) Dually, if $A \xrightarrow{\alpha} B \oplus C \to D \to A[1]$ is a triangle in $\mathcal{D}$ and $\mathrm{Hom}_\mathcal{D}(A, C) = 0$, then $D \cong D' \oplus C$, where $D'$ is the cone of the map $\alpha' \colon A \to B$ induced by $\alpha$.*

*Proof.* (1) The map $\alpha$ decomposes as $\alpha' \oplus 0 \colon C \oplus D \to A \oplus 0$ and we have the following isomorphism of triangles

$$\begin{array}{ccccccc}
A & \longrightarrow & B & \longrightarrow & C \oplus D & \xrightarrow{\alpha} & A[1] \\
\Big\| & & \Big\downarrow{\scriptstyle \mathbb{R}} & & \Big\| & & \Big\| \\
A \oplus 0 & \longrightarrow & B' \oplus D & \longrightarrow & C \oplus D & \xrightarrow{\alpha' \oplus 0} & A \oplus 0[1].
\end{array}$$

(2) is dual. □

### 4.7.1 A summary

We summarize here the outcome of the next sections, which are devoted to prove the following

**Proposition 4.9.** *Let $\overrightarrow{\mathcal{A}}$ be a category of quasicoherent sheaves over a weighted noncommutative regular projective curve $\mathbb{X}$ of genus $1$, with nontrivial weights. Then every large tilting sheaf $T$ is obtained by gluing appropriate tilting sheaves along a recollement of the form*

$$S^{\perp_\mathbb{Z}} \begin{array}{c} \xleftarrow{\ i^* \ } \\ \xrightarrow{\ i_* \ } \\ \xleftarrow{\ i^! \ } \end{array} \mathbf{D}(\overrightarrow{\mathcal{A}}) \begin{array}{c} \xleftarrow{\ j_! \ } \\ \xrightarrow{\ j^* \ } \\ \xleftarrow{\ j_* \ } \end{array} \mathrm{Loc}(S) \cong \mathbf{D}(\mathrm{Mod}\, k)$$

*for some partial tilting sheaf $S$ in $\overrightarrow{\mathcal{A}}$.*

In Section 4.7 we glue the tilting objects $S$ in $\mathrm{Loc}(S)$ and $T = T_{(B,V)}$ in $S^{\perp_\mathbb{Z}}$. The glued tilting object is (in the equivalence class of) $T_{(i_*B \oplus \tilde{S}), V}$, where $\tilde{S}$ is the indecomposable torsion sheaf of minimal length with the properties of having $S_\lambda = j_! S$ as simple socle and such that $\mathrm{Ext}^1_{\overrightarrow{\mathcal{A}}}(i_* T, \tilde{S}) = 0$. Actually, we just have to check that $\mathrm{Ext}^1_{\overrightarrow{\mathcal{A}}}(B, \tilde{S})$ vanishes for the summand $B$ of $T$ belonging to the tube $\mathcal{U}_{\bar{x}}$ in which $S_\lambda$ lies (see Figure 9).

The computation is divided in the two cases: $\bar{x}$ belongs to $V$ (Subsection 4.8.1), or $\bar{x}$ does not belong to $V$ (Subsection 4.8.2).

Next, in Section 4.9 we glue with respect to Theorem 2.16. Here, we have to distinguish some cases:

(1) If either

  - $\bar{x}$ belongs to $V$ (computed in Subsection 4.9.1), or
  - $\bar{x}$ does not belong to $V$ but $S_\rho = j_*S$ lies inside a wing corresponding to some indecomposable branch summand of $i_*T$ (see Figure 11, computed in Subsection 4.9.2),

  then the glued tilting object is (equivalent to) $T_{(i_*B \oplus \tilde{S}, V)}$, where $\tilde{S}$ is the indecomposable torsion sheaf of minimal length with the properties of having $S_\rho$ as simple top and such that $\mathrm{Ext}^1_{\overline{\mathcal{A}}}(\tilde{S}, i_*T) = 0$. Actually, we just have to check that $\mathrm{Ext}^1_{\overline{\mathcal{A}}}(\tilde{S}, B)$ vanishes for the summand $B$ of $T$ belonging to the tube $\mathcal{U}_{\bar{x}}$ in which $S_\lambda$ lies (see Figure 10).

(2) If $\bar{x}$ does not belong to $V$ and $\tau S_\rho$ belongs to $(i_*B_{\bar{x}})^{\perp_{0,1}}$, where $B_{\bar{x}}$ is the branch component of $T$ lying in $\mathcal{U}_{\bar{x}}$ (see Figure 12). In this case, the torsionfree part of $i_*T$ is modified, while its torsion part is not. The computation is done in Subsection 4.9.3 and the glued tilting sheaf is equivalent to $T_{(B,V)}$.

(3) If $\tau S_\rho$ lies inside a wing corresponding to some indecomposable branch summand of $i_*T$, but $S_\rho$ does not. In this case it is not clear what is the outcome of the gluing procedure.

Finally, in Section 4.10 we put together pieces and show that, even if our understanding of the outcome is not complete, every tilting sheaf in $\overrightarrow{\mathcal{A}}$ can be obtained by gluing appropriate tilting objects in $\overrightarrow{\mathcal{B}}$ and $\mathrm{Loc}(S)$.

## 4.8 Gluing, the first procedure

Consider the simple partial tilting object $S$ in $\mathrm{Loc}(S)$ and a tilting object

$$T = T_{(B,V)} = T_+ \oplus T_0$$

in $\overrightarrow{\mathcal{B}} \subseteq \mathbf{D}(\overrightarrow{\mathcal{B}})$. Let $\bar{x}$ be the point of $\mathbb{X}$ such that the simple object $S_\lambda = j_!S$ belongs to the tube $\mathcal{U}_{\bar{x}}$; let $\bar{r}$ be the rank of $\mathcal{U}_{\bar{x}}$. Denote by $\mathcal{U}'_{\bar{x}}$ the corresponding tube in $\overrightarrow{\mathcal{B}}$, which has rank $\bar{r} - 1$ (see Section 4.4), such that $i_*\mathcal{U}'_{\bar{x}} \subseteq \mathcal{U}_{\bar{x}}$. Write

$$T_0 = \bigoplus_{x \in \mathbb{X}} B_x = \bigoplus_{x \in \mathbb{X} \setminus \{\bar{x}\}} B_x \oplus B_{\bar{x}},$$

where $B_x$ is the direct sum of the branch and Prüfer modules supported in $x$, for each $x \in \mathbb{X}$. Let $B' = \bigoplus_{x \in \mathbb{X} \setminus \{\bar{x}\}} B_x$. We have to determine the set

$$\begin{aligned} I &= \mathrm{Hom}_{\mathbf{D}(\overline{\mathcal{A}})}(i_*T, S_\lambda[1]) \\ &\cong \mathrm{Ext}^1_{\overline{\mathcal{A}}}(i_*T, S_\lambda) \\ &\cong \mathrm{Ext}^1_{\overline{\mathcal{A}}}(i_*T_+, S_\lambda) \oplus \mathrm{Ext}^1_{\overline{\mathcal{A}}}(i_*B_{\bar{x}}, S_\lambda) \oplus \mathrm{Ext}^1_{\overline{\mathcal{A}}}(i_*B', S_\lambda). \end{aligned}$$

Here, the only possibly nonzero summand is $\operatorname{Ext}^1_{\mathcal{A}}(i_*B_{\bar{x}}, S_\lambda)$, because there are no nonsplit extensions between the tube $\mathcal{U}_{\bar{x}}$ and vector bundles or a tube $\vec{\mathcal{U}}_x$, with $x \neq \bar{x}$. Recall the tubes $\vec{\mathcal{U}}_x$ are serial categories, thus each nonzero indecomposable has exactly one simple subobject, its *socle*. From the Auslander-Reiten formula

$$D\operatorname{Ext}^1_{\mathcal{A}}(B, S_\lambda) \cong \operatorname{Hom}_{\overline{\mathcal{A}}}(\tau^- S_\lambda, B)$$

it follows that for an indecomposable object $B$ in $\vec{\mathcal{U}}_{\bar{x}}$, $\operatorname{Ext}^1_R(B, S_\lambda) \neq 0$ if and only if the simple socle of $B$ is $\tau^-(S_\lambda)$. This is also clear from the geometric model. More generally, by Equation 4.2 the dimension of $\operatorname{Ext}^1_{\mathcal{A}}(B, S_\lambda)$ equals the number of (nonisomorphic) indecomposable direct summands of $B$ with simple socle $\tau^-(S_\lambda)$. Notice also that if $B$ is indecomposable of length $1 \leq \ell \leq \infty$ with simple socle $\tau^-(S_\lambda)$, then the unique sheaf $E$ fitting in the short exact sequence

$$0 \to S_\lambda \to E \to B \to 0$$

is the one of length $\ell + 1$ and simple subobject $S_\lambda$.

The triangle (2.4) has the form

$$S_\lambda \to \tilde{T} \to i_*T^{(I)} \to S_\lambda[1]$$

and it is induced by the short exact sequence

$$0 \to S_\lambda \to \tilde{T} \to i_*T^{(I)} \to 0$$

in $\vec{\mathcal{A}}$. Up to multiplicity of the summands, $\tilde{T}$ may be written as

$$\tilde{T} = \tilde{T}_+ \oplus \tilde{B}$$
$$= i_*T_+ \oplus i_*B' \oplus \tilde{B}_{\bar{x}},$$

where $\tilde{B}_{\bar{x}}$ is the universal extension of $S_\lambda$ by $i_*B_{\bar{x}}$:

$$0 \to S_\lambda \to \tilde{B}_{\bar{x}} \to i_*B_{\bar{x}}^{(I)} \to 0.$$

Next, we determine the object $i_*B_{\bar{x}} \oplus \tilde{B}_{\bar{x}}$.

Recall that $\tilde{B}_{\bar{x}}$ has the following properties:

- $S_\lambda$ is a subobject of $\tilde{B}_{\bar{x}}$,
- $\operatorname{Ext}^1_{\mathcal{A}}(i_*B_{\bar{x}}, \tilde{B}_{\bar{x}}) = 0$ and $\operatorname{Ext}^1_{\mathcal{A}}(\tilde{B}_{\bar{x}}, i_*B_{\bar{x}}) = 0$.
- $i_*B_{\bar{x}}$ is a quotient of $\tilde{B}_{\bar{x}}$.

This translates into the following set of properties of the geometric model, assuming $S_\lambda$ is the simple object labeled $i$ in $\mathcal{U}_{\bar{x}}$ and identifying objects in $\vec{\mathcal{U}}_{\bar{x}}$ with their representation in the geometric model:

- $\tilde{B}_{\bar{x}}$ has at least one arc starting in $i - 1$,
- no arc of $i_*B_{\bar{x}}$ intersects any arc of $\tilde{B}_{\bar{x}}$,
- all arcs of $i_*B_{\bar{x}}$ share the ending vertex with some arc of $\tilde{B}_{\bar{x}}$.

Recall that objects with simple socle $S_\lambda$ are not in the image of $i_*$. Thus, $i_*B_{\bar{x}}$ has no arcs starting in $i - 1$ and therefore $\tilde{B}_{\bar{x}}$ has at least one arc, starting in $i - 1$ that is not isotopic to any arc of $i_*B_{\bar{x}}$.

### 4.8.1 Case $\bar{x} \in V$

If $\bar{x}$ belongs to $V$, then $B_{\bar{x}}$ is maximal exceptional in $\vec{\mathcal{B}}$ and thus $i_* B_{\bar{x}}$ has exactly $\bar{r} - 1$ nonisomorphic indecomposable summands; the sum $i_* B_{\bar{x}} \oplus \tilde{B}_{\bar{x}}$ has at least $\bar{r}$ nonisomorphic indecomposable summands, and thus exactly $\bar{r}$. Hence, the only indecomposable summand of $\tilde{B}_{\bar{x}}$ that is not isomorphic to any indecomposable summand of $i_* B_{\bar{x}}$ is the unique one with simple socle $S_\lambda$ and $\mathrm{Ext}^1$-orthogonal to $i_* B_{\bar{x}}$. See Figure 9 (a).

**Figure 9:** The component $B_{\bar{x}}$ of a tilting sheaf, with $\mathcal{U}_{\bar{x}}$ of rank 8. Here the expansion $i \colon \mathcal{B} \to \mathcal{A}$ sends the simple $M_{13}$ in $\mathcal{B}$ to $M_{13}$ (of length 2) in $\mathcal{A}$. On the right is represented $i_* B_{\bar{x}} \oplus \tilde{B}_{\bar{x}}$. The dashed arc represents the unique indecomposable summand of $\tilde{B}_{\bar{x}}$ that is not isomorphic to any summand of $i_* B_{\bar{x}}$. The dotted arc represents $S_\lambda$ and it is not a summand of $i_* B_{\bar{x}} \oplus \tilde{B}_{\bar{x}}$.

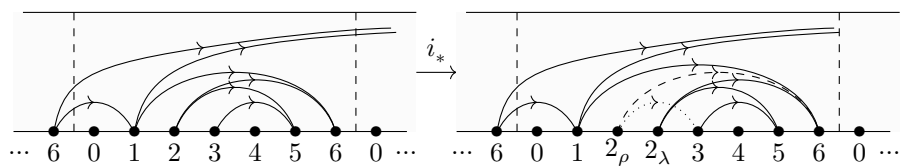

(a) Case $\bar{x} \in V$.

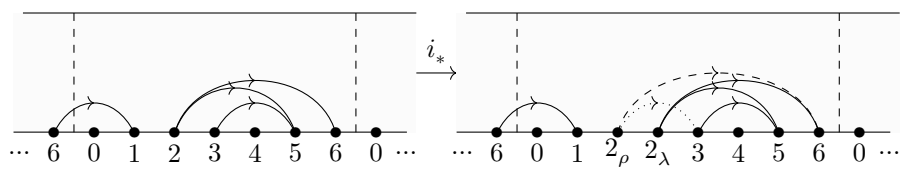

(b) Case $\bar{x} \in \mathbb{X} \setminus V$.

### 4.8.2 Case $\bar{x} \notin V$

Similarly, if $\bar{x}$ does not belong to $V$, then $B_{\bar{x}}$ is the direct sum of connected branches. There exists at most one connected branch $C$ which is a direct summand of $B_{\bar{x}}$ and such that $\mathrm{Ext}^1_{\overline{\mathcal{A}}}(i_* C, S_\lambda) \neq 0$. If there is no such $C$, then $\tilde{B}_{\bar{x}} \cong S_\lambda$; assume such $C$ exists and let $E$ be its root and $\ell$ the length of $E$; denote by $\mathcal{W}_C$ the wing rooted in $i_* E$ if $\mathrm{Ext}_{\overline{\mathcal{A}}}(i_* E, S_\lambda) = 0$, or rooted in the unique nontrivial extension of $S_\lambda$ by $i_* E$ if $\mathrm{Ext}_{\overline{\mathcal{A}}}(i_* E, S_\lambda) \neq 0$. Then

$$\mathrm{Ext}^1_{\overline{\mathcal{A}}}(i_* B_{\bar{x}}, S_\lambda) \cong \mathrm{Ext}_{\overline{\mathcal{A}}}(i_* C, S_\lambda)$$
$$\cong \mathrm{Ext}_{\mathcal{W}_C}(i_* C, S_\lambda).$$

Thus, we determine the universal extension

$$0 \to S_\lambda \to \tilde{C} \to i_* C^{(I)} \to 0$$

in $\mathcal{W}_C$. In $\mathcal{W}_C$ maximal exceptional objects have $\ell + 1$ nonisomorphic indecomposable summands. Since $i_* C$ has $\ell$ nonisomorphic indecomposable summands in $\mathcal{W}_C$, none of them with socle $S_\lambda$, and since $\tilde{C}$ has at least one indecomposable summand with socle $S_\lambda$, we can conclude as above that $i_*C \oplus \tilde{C}$ has exactly one indecomposable object with simple socle $S_\lambda$, which is the unique indecomposable summand of $\tilde{C}$ not appearing as a summand in $i_*C$. See Figure 9 (b).

In both cases ($x \in V$ or $x \notin V$) the glued silting object has exactly one indecomposable direct summand $\tilde{S}$ not isomorphic to a direct summand of $i_*T$; $\tilde{S}$ has $S_\lambda$ as simple socle and it is the only summand of $\tilde{T} \oplus i_*T$ with this property. When starting from the tilting object $T = T_{(B,V)}$ in $\vec{\mathcal{B}}$, the glued tilting object $\tilde{T} \oplus i_*T$ is the one corresponding to the pair $(i_*B \oplus \tilde{S}, V)$ in $\vec{\mathcal{A}}$.

If we let $S_\lambda$ run though all simple exceptional objects in $\vec{\mathcal{A}}$, and let $T$ run through all tilting objects in $S_\lambda^{\perp_\mathbb{Z}}$ concentrated in degree 0, we obtain a large part of the classification of tilting objects in $\vec{\mathcal{A}}$. To get a more complete picture of the classification, however, we need to employ the dual construction of Theorem 2.16 (and its Corollary 2.18).

## 4.9 Gluing on the other side

Consider again the recollement (4.4) induced by an expansion $i \colon \mathcal{B} \to \mathcal{A}$ of categories of coherent sheaves. Assume $\vec{\mathcal{A}} = \operatorname{Qcoh} \mathbb{X}$ and $\vec{\mathcal{B}} = \operatorname{Qcoh} \mathbb{X}'$. We first investigate the functor $j_*$, which plays a role in Theorem 2.16. To this aim, let $S_\lambda = j_! S$ as in previous subsection, where $S$ is the chosen simple exceptional compact object. Recall that $i^! S_\lambda = \bar{S}$, and notice that there is an isomorphism of triangles

$$\begin{array}{ccccccc}
i_* i^! S_\lambda & \longrightarrow & S_\lambda & \longrightarrow & j_* j^* S_\lambda & \longrightarrow & i_* i^! S_\lambda[1] \\
\downarrow \mathrel{\rotatebox{90}{$\simeq$}} & & \| & & \downarrow \mathrel{\rotatebox{90}{$\simeq$}} & & \downarrow \mathrel{\rotatebox{90}{$\simeq$}} \\
i_* \bar{S} & \longrightarrow & S_\lambda & \longrightarrow & S_\rho[1] & \longrightarrow & i_* \bar{S}[1],
\end{array}$$

the first one being the canonical triangle given by the recollement (4.4), the second one induced by the short exact sequence $0 \to S_\rho \to i_* \bar{S} \to S_\lambda \to 0$. Thus, $j_* S \cong j_* j^* j_! S = j_* j^* S_\lambda \cong S[1]$.

Consider then the simple partial tilting object $S[-1]$ in $\operatorname{Loc}(S)$ and the tilting object

$$T = T_{(B,V)} = T_+ \oplus T_0$$

in $\mathbf{D}(\vec{\mathcal{B}})$. Let $S_\rho = j_* S[-1]$ and note that $S_\rho = \tau S_\lambda$, where $S_\lambda = i_! S$. Let $\bar{x}$ be the point on $\mathbb{X}$ such that the simple object $S_\rho$ belongs to the tube $\mathcal{U}_{\bar{x}}$ of rank $\bar{r}$ in $\vec{\mathcal{A}}$. Write

$$T_0 = \bigoplus_{x \in \mathbb{X}} B_x = \bigoplus_{x \in \mathbb{X} \setminus \{\bar{x}\}} B_x \oplus B_{\bar{x}},$$

where $B_x$ is the direct sum of the branch and Prüfer modules supported in $x$. Let $B' = \bigoplus_{x \in \mathbb{X} \setminus \{\bar{x}\}} B_x$. We have to determine a set of generators $I$ of

$$\operatorname{Hom}_{\mathbf{D}(\vec{\mathcal{A}})}(S_\rho, i_*T[1]) \cong \operatorname{Ext}^1_{\vec{\mathcal{A}}}(S_\rho, i_*T)$$
$$\cong \operatorname{Ext}^1_{\vec{\mathcal{A}}}(S_\rho, i_*T_+) \oplus \operatorname{Ext}^1_{\vec{\mathcal{A}}}(S_\rho, i_*B_{\bar{x}}) \oplus \operatorname{Ext}^1_{\vec{\mathcal{A}}}(S_\rho, i_*B').$$

Theorem 2.16 requires the set $I$ to be finite. It is clear that $\operatorname{Ext}^1_{\overline{\mathcal{A}}}(S_\rho, i_*B') = 0$; we have to make sure the two remaining summands are finitely generated.

Since $i_*B_{\bar{x}}$ has at most $\bar{r}-1$ summands and $S_\rho$ is simple, then $\operatorname{Ext}^1_{\overline{\mathcal{A}}}(S_\rho, i_*B_{\bar{x}})$ has dimension at most $\bar{r}-1$. We investigate then $\operatorname{Ext}^1_{\overline{\mathcal{A}}}(S_\rho, i_*T_+)$.

The object $T_+$ is $V$-divisible in $\overrightarrow{\mathcal{B}}$ by [9, Theorem 4.8]; hence, when $\bar{x} \in V$,

$$\operatorname{Ext}^1_{\overline{\mathcal{A}}}(S_\rho, i_*T_+) = \operatorname{Hom}_{\mathbf{D}(\overline{\mathcal{A}})}(S_\rho, i_*T_+[1])$$
$$\cong \operatorname{Hom}_{\mathbf{D}(\overline{\mathcal{B}})}(i^*S_\rho, T_+[1])$$
$$\cong \operatorname{Ext}^1_{\overline{\mathcal{B}}}(i^*S_\rho, T_+) = 0.$$

Assume then $\bar{x} \notin V$ and, as in Section 4.5, write

$$T_0 = B_{\mathfrak{i}} \oplus B_{\mathfrak{e}} \oplus \bigoplus_{x \in V} \bigoplus_{j \in \mathcal{R}_x} \tau^j S_x[\infty].$$

We denote by $\overrightarrow{\mathcal{B}}'$ the full exact subcategory $(B_{\mathfrak{e}} \oplus \tau^- B_{\mathfrak{i}})^{\perp_{0,1}} \subseteq \overrightarrow{\mathcal{B}}$. Note that $\overrightarrow{\mathcal{B}}' = \operatorname{Qcoh} \mathbb{X}'$, where $\mathbb{X}'$ is a noncommutative regular projective curve with reduced weights ([9, Lemma 4.16(1)]).

Choose a canonical configuration $T'_{\operatorname{can}}$ in $\overrightarrow{\mathcal{B}}'$, that is a compact torsionfree tilting sheaf whose endomorphism algebra is canonical in the sense of [63] (see [9, §5], [44]). By [9, Theorem 5.8], the torsionfree summand $T_+$ is determined by the short exact sequence

$$0 \to T'_{\operatorname{can}} \to T_+ \to \bigoplus_{x \in V} \mathfrak{S}_x \to 0 \tag{4.5}$$

in $\overrightarrow{\mathcal{B}}'$, where $\mathfrak{S}_x$ is a direct sum of Prüfer sheaves supported in $x$ such that $\operatorname{Add}(\mathfrak{S}_x) = \operatorname{Add}(\bigoplus_{j \in \mathcal{R}_x} \tau^j S_x[\infty])$.

Notice that $\operatorname{Hom}_{\overline{\mathcal{A}}}(S_\rho, \bigoplus_{x \in V} \mathfrak{S}_x) = 0 = \operatorname{Ext}^1_{\overline{\mathcal{A}}}(S_\rho, \bigoplus_{x \in V} \mathfrak{S}_x)$ since $\bar{x} \notin V$. Thus, the exact sequence

$$\operatorname{Hom}_{\overline{\mathcal{A}}}(S_\rho, \bigoplus_{x \in V} \mathfrak{S}_x) \to \operatorname{Ext}^1_{\overline{\mathcal{A}}}(S_\rho, i_*T'_{\operatorname{can}}) \to \operatorname{Ext}^1_{\overline{\mathcal{A}}}(S_\rho, T_+) \to \operatorname{Ext}^1_{\overline{\mathcal{A}}}(S_\rho, \bigoplus_{x \in V} \mathfrak{S}_x)$$

yields $\operatorname{Ext}^1_{\overline{\mathcal{A}}}(S_\rho, T_+) \cong \operatorname{Ext}^1_{\overline{\mathcal{A}}}(S_\rho, T'_{\operatorname{can}}) \cong \operatorname{Ext}^1_{\overline{\mathcal{B}}}(i^*S_\rho, T'_{\operatorname{can}})$. The vector bundle $T'_{\operatorname{can}}$ is finitely presented in $\overrightarrow{\mathcal{B}}'$, hence $\operatorname{Ext}^1_{\overline{\mathcal{B}}'}(i^*S_\rho, T'_{\operatorname{can}})$ is finitely generated and so is $\operatorname{Ext}^1_{\overline{\mathcal{B}}}(i^*S_\rho, T'_{\operatorname{can}})$, because $\overrightarrow{\mathcal{B}}'$ is an exact subcategory of $\overrightarrow{\mathcal{B}}$.

Thus, the set of generators $I$ of $\operatorname{Ext}^1_{\overline{\mathcal{A}}}(S_\rho, i_*T_+) \oplus \operatorname{Ext}^1_{\overline{\mathcal{A}}}(S_\rho, i_*B_{\bar{x}})$ can be chosen to be finite and the hypotheses of Corollary 2.18 are satisfied.

We construct then $\tilde{T}$ as the cocone of the right universal map $\alpha$:

$$i_*T^I \to \tilde{T} \to S_\rho \xrightarrow{\alpha} i_*T^I[1].$$

Such a triangle induces the short exact sequence

$$0 \to i_*(T_+ \oplus B_{\bar{x}} \oplus B')^I \to \tilde{T} \to S_\rho \to 0 \tag{4.6}$$

in $\overrightarrow{\mathcal{A}}$.

In view of Proposition 4.7 and since $\operatorname{Hom}_{\overline{\mathcal{A}}}(\tilde{T}, \tilde{T}) = 0$, we can write $\tilde{T} = \tilde{T}_+ \oplus \tilde{T}_0 = \tilde{T}_+ \oplus \tilde{B}_{\bar{x}} \oplus \tilde{B}'$ as the direct sum of its torsionfree part $\tilde{T}_+$, the torsion summand $\tilde{B}_{\bar{x}}$ belonging to $\vec{\mathcal{U}}_{\bar{x}}$ and the sum $\tilde{B}'$ of torsion summands in tubes $\vec{\mathcal{U}}_x$ with $x \neq \bar{x}$.

We determine the torsion part $\tilde{T}_0 \oplus i_*T_0$ of $\tilde{T} \oplus i_*T$. The component $\tilde{B}' \oplus i_*B'$ belonging to tubes $\vec{\mathcal{U}}_x$ with $x \neq \bar{x}$ belongs to $\mathrm{Add}(i_*B')$: indeed, since $\mathrm{Ext}^1_{\mathcal{A}}(S_\rho, i_*B') = 0$, then $i_*B'^I$ appears as a summand of $\tilde{T}$ (see Lemma 4.8) and there is no indecomposable summand supported in points $x \neq \bar{x}$ being not isomorphic to a summand of $i_*B'$.

We concentrate then on the component $\tilde{B}_{\bar{x}} \oplus i_*B_{\bar{x}}$ lying in the tube $\vec{\mathcal{U}}_{\bar{x}}$. The short exact sequence (4.6) implies that $\tilde{B}_{\bar{x}}$ has the following properties:

- $S_\rho$ is a quotient of $\tilde{B}_{\bar{x}}$,
- $\mathrm{Ext}^1_{\mathcal{A}}(\tilde{B}_{\bar{x}}, i_*B_{\bar{x}}) = 0$ and $\mathrm{Ext}^1_{\mathcal{A}}(i_*B_{\bar{x}}, \tilde{B}_{\bar{x}}) = 0$.
- $i_*B_{\bar{x}}$ is a subobject of $\tilde{B}_{\bar{x}}$.

We distinguish now some cases in which we can see what happens.

### 4.9.1 Case $\bar{x} \in V$

In case $\bar{x} \in V$, it holds additionally $\mathrm{Ext}^1_{\mathcal{A}}(S_\rho, i_*T_+) = 0$ (see Section 4.5). Thus, by Lemma 4.8 and since $\mathrm{Ext}^1_{\mathcal{A}}(S_{\rho, i_*B'}) = 0$, $\tilde{T}$ is just determined by the extension of $B_{\bar{x}}$ by $S_\rho$ in $\vec{\mathcal{U}}_{\bar{x}}$.

A maximality argument shows that there is exactly one indecomposable summand $\tilde{S}$ of $\tilde{B}_{\bar{x}}$ (up to multiplicity) with the three properties above (see also Figure 10 (a)). The glued tilting object in this case is equivalent to the tilting object $T_{(i_*B \oplus \tilde{S}, V)}$ in $\vec{\mathcal{A}}$.

**Figure 10:** The component $\tilde{B}_{\bar{x}}$ of a tilting sheaf when gluing through a right universal map. The expansion $i\colon \mathcal{B} \to \mathcal{A}$ sends the simple $M_{46}$ in $\mathcal{B}$ to $M_{46}$ (of length 2) in $\mathcal{A}$. The dashed arc represents the unique indecomposable summand $\tilde{S}$ of $\tilde{B}_{\bar{x}}$ that is not isomorphic to any summand of $i_*B_{\bar{x}}$. The dotted arc represents $S_\rho$ and it is not a summand of $i_*B_{\bar{x}} \oplus \tilde{B}_{\bar{x}}$.

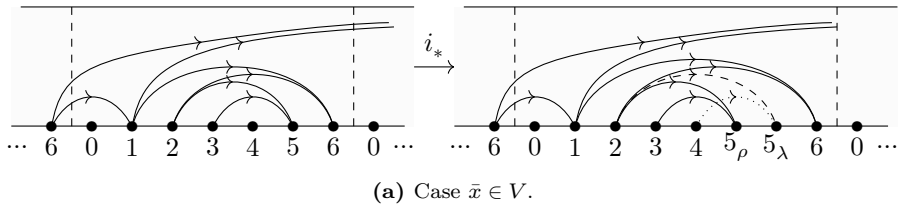

(a) Case $\bar{x} \in V$.

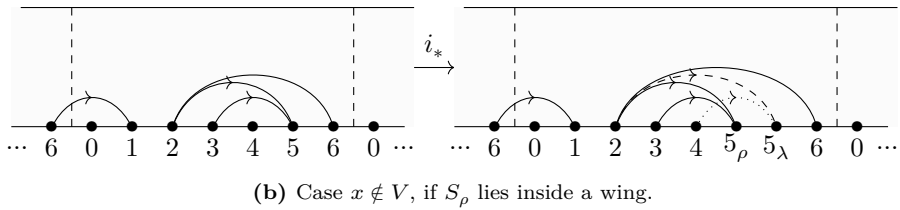

(b) Case $x \notin V$, if $S_\rho$ lies inside a wing.

### 4.9.2 If $\bar{x} \notin V$, first case

**Figure 11:** Case $\bar{x} \notin V$ and $S_\rho$ belongs to the wing $\mathcal{W}_{\bar{x}}$. The solid line is the root of the wing, the dotted lines are the possible arcs for $S_\rho$.

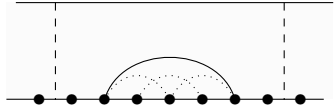

Let $\mathcal{W}_{\bar{x}}$ be the union of wings corresponding to $B_{\bar{x}}$ in $\overrightarrow{\mathcal{A}}$ (notice that $i_* B_{\bar{x}}$ is no more a branch module in $\overrightarrow{\mathcal{A}}$, but we can still associate a union of wings to it). In case $S_\rho$ belongs to $\mathcal{W}_{\bar{x}}$, then

$$\operatorname{Ext}^1_{\overrightarrow{\mathcal{A}}}(S_\rho, i_* T_+) \cong \operatorname{Hom}_{\mathbf{D}(\overrightarrow{\mathcal{A}})}(S_\rho, i_* T_+[1]) \cong \operatorname{Ext}^1_{\overrightarrow{\mathcal{B}}}(i^* S_\rho, T_+) = 0,$$

since $\bar{S} = i^* S_\rho$ belongs to the wing $\mathcal{W}_{\bar{x}}$ associated to $B_{\bar{x}}$ in $\overrightarrow{\mathcal{B}}$ and thus $\bar{S}$ and $T_+$ are direct summands of the same tilting object in $\overrightarrow{\mathcal{B}}$ [9].

Thus, $\operatorname{Ext}^1_{\overrightarrow{\mathcal{A}}}(S_\rho, i_* T) \cong \operatorname{Ext}^1_{\overrightarrow{\mathcal{A}}}(S_\rho, i_* B_{\bar{x}})$. This means the glued tilting object is determined by the extensions of $i_* B_{\bar{x}}$ by $S_\rho$, which can be computed within the abelian category $\operatorname{add}(\mathcal{W}_{\bar{x}})$. If we denote by $\tilde{S}$ the unique indecomposable object in $\mathcal{W}_{\bar{x}}$ having simple top $S_\rho$ and no nontrivial extensions with $i_* B_{\bar{x}}$, then the glued tilting object $\tilde{T} \oplus i_* T$ is equivalent to the tilting object $T_{(i_* B \oplus \tilde{S}, V)}$ in $\overrightarrow{\mathcal{A}}$. The situation is depicted in Figure 10 (b).

### 4.9.3 If $\bar{x} \notin V$, second case

**Figure 12:** Case $\bar{x} \notin V$ and $\tau S_\rho$ belongs to $(i_* B_{\bar{x}})^{\perp_{0,1}}$. The solid line is the root of the wing, the dotted lines are the possible arcs for $S_\rho$.

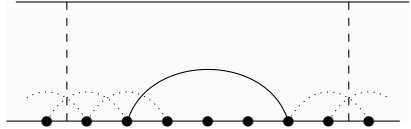

Assume then $\bar{x} \notin V$ and $\tau S_\rho$ belongs to $(i_* B_{\bar{x}})^{\perp_{0,1}}$. Then,

$$0 = \operatorname{Hom}_{\overrightarrow{\mathcal{A}}}(i_* B_{\bar{x}}, \tau S_\rho) \cong D\operatorname{Ext}^1_{\overrightarrow{\mathcal{A}}}(S_\rho, i_* B_{\bar{x}}).$$

Thus, $\operatorname{Ext}^1_{\overrightarrow{\mathcal{A}}}(S_\rho, i_* T) \cong \operatorname{Ext}^1_{\overrightarrow{\mathcal{A}}}(S_\rho, i_* T_+)$. Moreover,

$$D\operatorname{Ext}^1_{\overrightarrow{\mathcal{A}}}(S_\rho, i_* T_+) \cong \operatorname{Hom}_{\overrightarrow{\mathcal{A}}}(i_* T_+, \tau S_\rho).$$

An application of $\operatorname{Hom}_{\overrightarrow{\mathcal{A}}}(-, \tau S_\rho)$ to the short exact sequence (4.5) defining $T_+$ shows that $\operatorname{Hom}_{\overrightarrow{\mathcal{A}}}(i_* T_+, \tau S_\rho) \cong \operatorname{Hom}_{\overrightarrow{\mathcal{A}}}(i_* T'_{\operatorname{can}}, \tau S_\rho)$, since $\operatorname{Hom}_{\overrightarrow{\mathcal{A}}}(\bigoplus_{x \in V} \mathfrak{S}_x, \tau S_\rho[k]) = 0$ for each $k \in \mathbb{Z}$. But

$$\operatorname{Hom}_{\overrightarrow{\mathcal{A}}}(i_* T'_{\operatorname{can}}, \tau S_\rho) \cong \operatorname{Hom}_{\overrightarrow{\mathcal{B}}}(T'_{\operatorname{can}}, i^! \tau S_\rho)$$
$$= \operatorname{Hom}_{\overrightarrow{\mathcal{B}'}}(T'_{\operatorname{can}}, \tau S_\rho)$$

(recall $i^!$ is the identity on simple objects different from $S_\rho$ and $S_\lambda$). This is nonzero, since $T'_{\text{can}}$ is tilting in $\vec{\mathcal{B}}'$ ([44, Proposition 5.4]) and $\text{Ext}^1_{\vec{\mathcal{B}}'}(T'_{\text{can}}, \tau S_\rho) = 0$. We conclude that $\text{Ext}^1_{\vec{\mathcal{A}}}(S_\rho, T_+)$ is nonzero.

Notice that the extensions of $T_+$ by $S_\rho$ do not contribute to the torsion part of $\tilde{T}$: for any non-split short exact sequence

$$0 \to T_+ \to E \xrightarrow{g} S_\rho \to 0,$$

the sheaf $E$ does not have any nonzero torsion subobject. Indeed, assume $E_0$ is a torsion subobject of $E$ not isomorphic to $S_\rho$. Then there is a diagram as follows, with exact rows and where vertical arrows are inclusions:

$$\begin{array}{ccccccccc} 0 & \longrightarrow & \ker g' & \longrightarrow & E_0 & \xrightarrow{g'} & S_\rho & \longrightarrow & 0 \\ & & \downarrow & & \downarrow & & \| & & \\ 0 & \longrightarrow & T_+ & \longrightarrow & E & \xrightarrow{g} & S_\rho & \longrightarrow & 0. \end{array}$$

Since $\vec{\mathcal{A}}_0$ is an abelian subcategory of $\vec{\mathcal{A}}$, $\ker g'$ is a torsion subobject of the torsion subobject $T_+$, and thus it is zero. This means $E_0 \cong S_\rho$, but in this case the sequence would split.

We conclude that in this case, $\text{Add}(\tilde{B}_{\bar{x}} \oplus B_{\bar{x}}) = \text{Add}(B_{\bar{x}})$ and the glued object is the tilting object $T_{(i_*B, V)}$ in $\vec{\mathcal{A}}$.

### 4.9.4 If $\bar{x} \notin V$, third case

Consider now the case that $\tau S_\rho$ belongs to $\mathcal{W}_{\bar{x}}$ but $S_\rho$ does not, see Figure 13. In this case, $\tau S_\rho$ does not belong to $B_{\bar{x}}^{\perp_{0,1}}$ as $\text{Ext}_{\vec{\mathcal{A}}}(S_\rho, B_{\bar{x}}) \neq 0$, but we cannot show that $\text{Ext}_{\vec{\mathcal{A}}}(S_\rho, T_+) = 0$ as in Subsection 4.9.2. I am not certain whether or not $\tilde{T}$ provides a new summand in $\vec{\mathcal{U}}_{\bar{x}}$.

**Figure 13:** Case $\bar{x} \in V$, with $\tau S_\rho \in \mathcal{W}_{\bar{x}}$ and $S_\rho \notin \mathcal{W}_{\bar{x}}$. The solid line is the root of the wing $\mathcal{W}_{\bar{x}}$, the dotted line is $S_\rho$.

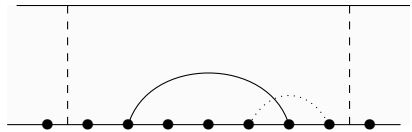

## 4.10 Conclusion

Even if our understanding of the behavior of glued tilting sheaves is not complete, we can retrieve the whole classification of large tilting sheaves over a weighted (with nontrivial weights) noncommutative regular projective curves of genus zero, by gluing along the recollement induced by an expansion and known the classification on curves with smaller weights.

Consider then a category $\vec{\mathcal{A}} = \text{Qcoh}\,\mathbb{X}$ as above, and assume there is at least one non-homogeneous tube $\vec{\mathcal{U}}_x$ of rank $r \geq 2$. Consider a large tilting object $T = T_{(B,V)}$ in $\vec{\mathcal{A}}$ and denote by $B_x$ the summand of $T$ supported in $\vec{\mathcal{U}}_x$. Our goal is to appropriately choose a simple object $S_\lambda$ in $\mathcal{U}_x$ (and hence $S_\rho = \tau S_\lambda$) and a tilting sheaf $T'$ in $\vec{\mathcal{B}} = \operatorname{Qcoh} \mathbb{X}' = S^{\perp_{0,1}}$ such that $T$ results from gluing $T'$ and $S_\lambda$ or $S_\rho$ along the recollement

$$\mathbf{D}(\vec{\mathcal{B}}) \xleftarrow{\;\;i^*\;\;}_{\xleftarrow{\;\;i^!\;\;}} \xrightarrow{i_*} \mathbf{D}(\vec{\mathcal{A}}) \xleftarrow{\;\;j_!\;\;}_{\xleftarrow{\;\;j_*\;\;}} \xrightarrow{j^*} \operatorname{Loc}(S_\lambda).$$

with respect to either the method of Theorem 2.10 or of Theorem 2.16.

### 4.10.1 Case 1. $B_x = 0$

Assume first that $B_x = 0$ (and thus necessarily $x \notin V$). In this case, we can choose $S_\lambda$ to be any simple object in $\mathcal{U}_x$, and $T'$ the tilting object corresponding to the pair $(B, V)$ in $\vec{\mathcal{B}}$. Then gluing with respect to Theorem 2.16, as in 4.9, gives the desired object $T$.

### 4.10.2 Case 2. $B_x$ has a nonzero branch component

Assume now that $B_x$ is nonzero and it has a branch summand. In this case, $B_x$ has at least one simple direct summand $S_1$. Write $B_x = B_1 \oplus S_1$, so that $B_1$ has no summands isomorphic to $S_1$. Then exactly one of the following conditions is satisfied:

(a) $\operatorname{Hom}_{\vec{\mathcal{A}}}(B_1, S_1) \neq 0$;

(b) $\operatorname{Hom}_{\vec{\mathcal{A}}}(S_1, B_1) \neq 0$;

(c) $\operatorname{Hom}_{\vec{\mathcal{A}}}(B_1, S_1) = 0 = \operatorname{Hom}_{\vec{\mathcal{A}}}(S_1, B_1)$.

In case (a), choose $S_\lambda = S_1$ and let $T'$ be the tilting object corresponding to the pair $(B_1, V)$ in $\vec{\mathcal{B}}$. Notice that $B_1$ is a branch sheaf in $\vec{\mathcal{B}}$ and it is isomorphic to $i^* B_x$. Glue with respect to Theorem 2.10, as in Section 4.7, to get back the tilting sheaf $T$. See Figure 14.

**Figure 14:** Case 2(a). To the left is represented the branch component $B_x$ in $\mathcal{A}$; to the right is its image $i^* B_x$. Note that the image $i^* S_\lambda$ of the simple summand represented by the dashed line is zero.

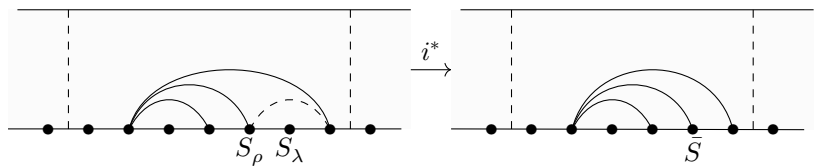

In case (b), choose $S_\lambda = \tau^- S_1$ (hence $S_1 = S_\rho$) and let $T'$ be the tilting object corresponding to the pair $(B_1, V)$ in $\vec{\mathcal{B}}$. The object $B_1$ is a branch sheaf in $\vec{\mathcal{B}}$ and it is isomorphic to $i^! B$. Glue with respect to Theorem 2.16, as in Section 4.9, to get back the tilting sheaf $T$. See Figure 15.

Case (c) occurs exactly when $S_1$ is the root of the connected branch containing it. Here the choice is the same as case (a). See Figure 16

**Figure 15:** Case 2(b). To the left is represented the branch component $B_x$ in $\mathcal{A}$; to the right is its image $i^!B_x$. Note that the image $i^!S_\rho$ of the simple summand represented by the dashed line is zero.

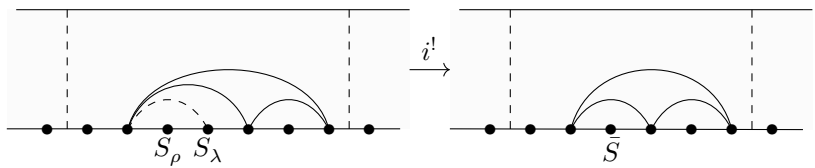

**Figure 16:** Case 2(c). As for case (a), the object to the left is the image $i^*B_x$ of the branch component $B_x$ in $\mathcal{A}$. Note that the image $i^*S_\lambda$ of the simple summand represented by the dashed line is zero.

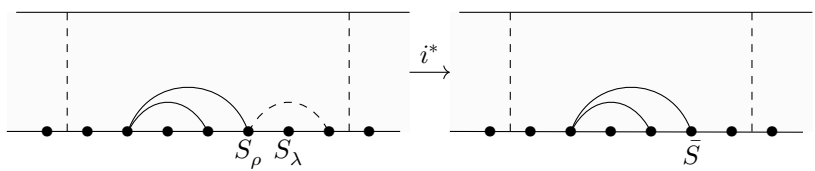

### 4.10.3   Case 3. $B_x$ is a sum of Prüfer sheaves

In this case, $x$ belongs necessarily to $V$. Pick any simple $S_\lambda$ in $\mathcal{U}_x$ and the tilting sheaf $T'$ corresponding to $(B, V)$ in $\vec{\mathcal{B}}$, whose component in $\vec{\mathcal{U}}'_x$ is again the sum of all Prüfer sheaves the tube. Gluing $S_\lambda$ and $T'$ with respect to Theorem 2.10, as in Section 4.7, gives back the tilting sheaf $T$. See Figure 17.

**Figure 17:** Case 3. An object $B_x$ in $\vec{\mathcal{A}}$ that is a direct sum of all Prüfer objects supported in $x$ is sent to the object $i^*B$ that is the direct sum in $\vec{\mathcal{B}}$ of all Prüfer supported in $x$. The Prüfer with simple socle $S_\lambda$ is sent to zero.

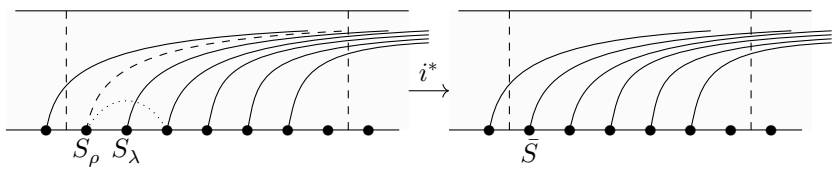

# Bibliography

[1]  T. Adachi, O. Iyama, and I. Reiten. "$\tau$-tilting theory". *Compositio Mathematica* 150.3 (2014), pp. 415–452.

[2]  T. Aihara and O. Iyama. "Silting mutation in triangulated categories". *Journal of the London Mathematical Society* 85.3 (2012), pp. 633–668.

[3]  F. W. Anderson and K. R. Fuller. *Rings and categories of modules*. Vol. 13. Springer Science & Business Media, 2012.

[4]  L. Angeleri Hügel. "Infinite dimensional tilting theory". Ed. by D. J. Benson, H. Krause, and A. Skowronski. EMS European Mathematical Society, 2013, pp. 1–37.

[5]  L. Angeleri Hügel. "Silting objects". *ArXiv e-prints* (2018). To appear in Bulletin of the London Math. Soc. eprint: `1809.02815`.

[6]  L. Angeleri Hügel, S. Koenig, and Q. Liu. "Recollements and tilting objects". *Journal of Pure and Applied Algebra* 215.4 (2011), pp. 420–438.

[7]  L. Angeleri Hügel, S. Koenig, Q. Liu, and D. Yang. "Ladders and simplicity of derived module categories". *Journal of Algebra* 472 (2017), pp. 15–66.

[8]  L. Angeleri Hügel and D. Kussin. "Tilting and cotilting modules over concealed canonical algebras". *Mathematische Zeitschrift* (2015), pp. 1–30.

[9]  L. Angeleri Hügel and D. Kussin. "Large tilting sheaves over weighted noncommutative regular projective curves". *Documenta Mathematica* 22 (2017), pp. 67–134.

[10] L. Angeleri Hügel, F. Marks, J. Šťovíček, R. Takahashi, and J. Vitória. "Silting and flat ring epimorphisms". 2019.

[11] L. Angeleri Hügel, F. Marks, and J. Vitória. "Silting modules". *International Mathematics Research Notices* 2016.4 (2015), pp. 1251–1284.

[12] L. Angeleri Hügel, F. Marks, and J. Vitória. "Silting modules and ring epimorphisms". *Advances in Mathematics* 303 (2016), pp. 1044–1076.

[13] L. Angeleri Hügel, F. Marks, and J. Vitória. "Partial silting objects and smashing subcategories". *ArXiv e-prints* (2019). arXiv: `1902.05817`.

[14] L. Angeleri Hügel and J. Sánchez. "Tilting modules arising from ring epimorphisms". *Algebras and Representation Theory* 14.2 (2011), pp. 217–246.

[15] L. Angeleri Hügel and J. Sánchez. "Tilting modules over tame hereditary algebras". *Journal für die reine und angewandte Mathematik (Crelles Journal)* 2013.682 (2013), pp. 1–48.

[16] K. Baur, A. B. Buan, and R. J. Marsh. "Torsion pairs and rigid objects in tubes". *Algebras and Representation Theory* 17.2 (2014), pp. 565–591.

[17] A. A. Beilinson, J. Bernstein, and P. Deligne. "Faisceaux pervers". *Astérisque* 100 (1982).

[18] A. Beligiannis and I. Reiten. *Homological and homotopical aspects of torsion theories*. American Mathematical Soc., 2007.